\documentclass[bj,11pt,preprint]{imsart}

\RequirePackage{amsthm,amsmath, enumerate}
\RequirePackage[colorlinks,citecolor=blue,urlcolor=blue]{hyperref}
\usepackage{mathrsfs}
\usepackage[margin=1in]{geometry}
\usepackage[english]{babel}
\usepackage{multicol}
\usepackage{pdfpages}
\usepackage{wrapfig}
\usepackage{tabularx}
\usepackage{caption}
\usepackage{bm}
\usepackage{mathrsfs}
\usepackage{amsmath}
\usepackage{enumitem}
\usepackage{pdfpages}
\usepackage{color}
\usepackage{float}
\usepackage{graphicx}
\usepackage{csquotes}
\usepackage{hyperref}
\usepackage{amsmath}
\usepackage{mathrsfs}
\usepackage[utf8]{inputenc}
\usepackage[english]{babel}
\usepackage{multicol}
\usepackage{pdfpages}
\usepackage{wrapfig}
\usepackage{color}
\usepackage{bbm}
\usepackage{tabularx}
\usepackage{caption}
\usepackage{bm}
\usepackage{mathrsfs}
\usepackage{amssymb}
\usepackage{enumitem}
\usepackage{pdfpages}
\usepackage{float}
\usepackage{bigints}
\usepackage{graphicx}
\usepackage{csquotes}
\usepackage{hyperref}
\usepackage{amsmath}
\usepackage{ulem}
\usepackage{amssymb}
\usepackage[overload]{empheq}

\usepackage{bbm}
\usepackage{multirow}

\startlocaldefs
\numberwithin{equation}{section}
\theoremstyle{plain}

\theoremstyle{remark}
\newtheorem{rem}{Remark}
\theoremstyle{plain}

\theoremstyle{plain}

\theoremstyle{plain}

\theoremstyle{plain}

\theoremstyle{plain}

\theoremstyle{plain}

\endlocaldefs

\newtheorem{theo}{Theorem}
\newtheorem{lem}{Lemma}

\numberwithin{equation}{section}
\numberwithin{thm}{section}
\numberwithin{cor}{section}
\numberwithin{lemma}{section}
\numberwithin{proposition}{section}
\numberwithin{rem}{section}
\numberwithin{example}{section}
\numberwithin{res}{section}

\newtheorem{innercustomgeneric}{\customgenericname}
\providecommand{\customgenericname}{}
\newcommand{\newcustomtheorem}[2]{%
  \newenvironment{#1}[1]
  {%
   \renewcommand\customgenericname{#2}%
   \renewcommand\theinnercustomgeneric{##1}%
   \innercustomgeneric
  }
  {\endinnercustomgeneric}
}

\newcustomtheorem{customthm}{Theorem}
\newcustomtheorem{customlemma}{Lemma}

\begin{document}

\begin{frontmatter}

\linespread{1}
\title{On Second order correctness of Bootstrap in Logistic Regression}%Pebble: Second Order Theory of Bootstrap in Logistic regression}
\runtitle{SOC of Bootstrap in Logistic Regression}

\begin{aug}
\author{\fnms{Debraj} \snm{Das}$^{\text{a}\thanksref{T1}}$\ead[label=c]{rajdas@iitk.ac.in}}
\and
\author{\fnms{Priyam} \snm{Das}$^{\text{b}}$\ead[label=d]{priyam\_das@hms.harvard.edu}}

\thankstext{T1}{Research partially supported by DST Inspire fellowship DST/INSPIRE/04/2018/001290}

\address{$^{\text{a}}$Department of Mathematics and Statistics, Indian Institute of Technology Kanpur,
India.
\printead{c}}
\address{$^{\text{b}}$Department of Biomedical Informatics, Harvard Medical School, Boston, USA.\\
\printead{d}}

\runauthor{Das D. and Das P.}

\end{aug}

\linespread{1.25}

\begin{abstract} $\;$In the fields of clinical trial, biomedical surveys, marketing, banking, with dichotomous response variable, the logistic regression is considered as an alternative convenient approach to linear regression. In this paper, we develop a novel perturbation Bootstrap technique for approximating the distribution of the maximum likelihood estimator (MLE) of the regression parameter vector. We establish second order correctness for the proposed Bootstrap method which results in improved inference performance  compared to that based on asymptotic normality. The main challenge in establishing second order correctness remains in the fact that the response variable being binary, the resulting MLE has a lattice structure. We show that the direct Bootstrapping approach fails even after studentization. We adopt smoothing technique developed in Lahiri (1993) to ensure that the smoothed studentized version of the MLE has a density. Similar smoothing strategy is employed to the Bootstrap version to achieve second order correct approximation. Good finite-sample properties of the proposed Bootstrap method is shown using simulation experiments. The proposed method is used to find the confidence intervals of the coefficients of the covariates on a dataset in the field of healthcare operations decision.
\end{abstract}


\end{frontmatter}

\section{Introduction}\label{sec:intro}
Logistic regression is one of the most widely used regression techniques when the response variable is binary. The use of the `logit' function as a statistical tool dates back to Berkson (1944), followed by Cox (1958), who popularized it in the field of regression. Following those seminal works, numerous applications of logistic regression can be found in different fields, from banking sectors to epidemiology, clinical trials, biomedical surveys, among others (Hosmer, Lemeshow and Sturdivant (2013)).
%An overview on different applications of logistic regression can be found in
 The logistic regression model is given as follows. Suppose $y$ denotes the binary response variable and the value of $y$ depends on the $p$ independent variables $\bm{x} = (x_1,\ldots,x_p)^\prime$. Instead of capturing this dependence by modelling $y$ directly on the covariates, in logistic regression, log-odds corresponding to the success of $y$, denoted by $p(\mathbf{x})=P(y=1)$, is modeled as a linear function of the covariates.
% $x_1,\dots,x_p$. 
%Let us look into the model more mathematically. 
% Suppose $y_1,\dots,y_n$ are $n$ independent observed values of the response variable. $x_{ij}$ denotes the $i$-th value of the $j$-th independent variable \bl{(or, regressor)} and $\mathbf{x}_i=(x_{i1},\dots,x_{ip})$ is the $i$th design vector where $i=1,\dots, n$ and $j=1,\dots,p$. Let $\mathbf{P}(y_i=1)=p(\mathbf{x}_i)$, i.e. $p(\mathbf{x}_i)$ is the probability of success of $y_i$. 
The odds ratio for the event $\{y=1\}$ is given by $odd(\bm{x}) = \dfrac{p(\bm{x})}{1-p(\bm{x})}$. 
The logistic regression model is given by 
\begin{align}\label{eqn:model}
\text{logit} (p(\bm{x}))=\log\bigg[\dfrac{p(\bm{x})}{1-p(\bm{x})}\bigg]=\bm{x}^{\prime}\bm{\beta},
\end{align}
where $\bm{\beta}=(\beta_1,\ldots, \beta_p)$ is the $p$-dimensional vector of regression parameters. In convention, the maximum likelihood estimator (MLE) of $\bm{\beta}$ is used for the purpose of inference. 
%In this paper we also stick to the MLE of $\bm{\beta}$. %In this article we assume the design vectors to be fixed. 
%\bl{(LETS PUT IT LATER, POSSIBLY ITS TRIVIAL, NO NEED TO STATE?, NOT NEEDED IN INTRODUCTION ``It is routine to consider the maximum likelihood estimator (MLE) of $\boldsymbol{\beta}$ for the purpose of inference. In this paper we also stick to the MLE'')}.
For a given sample $\{(\mathbf{x_i},y_i)\}_{i=1}^n$, the likelihood is given by
\begin{align*}
L(\bm{\beta}|y_1,\dots,y_n, \bm{x}_1,\dots,\bm{x}_n)=\prod_{i=1}^{n}p(\bm{x}_i)^{y_i}(1-p(\bm{x}_i))^{1-y_i},    
\end{align*}
where $p(\bm{\beta}|\bm{x}_i)=\dfrac{e^{\bm{x}_i^{\prime}\bm{\beta}}}{1+e^{\bm{x}_i^{\prime}\bm{\beta}}}$.
% as implied from the model (\ref{eqn:model}). 
The MLE $\hat{\bm{\beta}}_n$ of $\boldsymbol{\beta}$ is defined as the maximizer of $L(\bm{\beta}|y_1,\dots,y_n, \bm{x}_1,\dots,\bm{x}_n)$, which is obtained by solving
\begin{align}\label{eqn:def}
\sum_{i=1}^{n}(y_i-p(\bm{\beta}|\bm{x}_i))\bm{x}_i=0.
\end{align}
%(since the second derivative of the left hand side (LHS) of (\ref{eqn:def}) with respect to $\bm{\beta}$ is negative definite.) We assume that the design vectors $\bm{x}_1,\dots,\bm{x}_n$ are fixed.
In order to find confidence intervals for different regression coefficients or to test whether a certain covariate is of importance or not, it is required to find a good approximatation of the distribution of $\hat{\bm{\beta}}_n$. $\hat{\bm{\beta}}_n$ being the MLE, the distribution of $\hat{\bm{\beta}}_n$ is approximately normal under certain regularity conditions. 
% This asymptotic distribution can be used for inferences regarding the regression parameter $\bm{\beta}$.
Asymptotic normality as well as other large sample properties of $\hat{\bm{\beta}}_n$ have been studied extensively in the literature (cf. Haberman (1974), McFadden (1974), Amemiya (1976), Gourieroux and Monfort (1981), Fahrmeir and Kaufmann (1985)).

%\bl{(ADD ONE LINE ON HOW Bootstrap IS BETTER THAN PREVIOUS APPROACH)}
As an alternative to asymptotic normality, Efron (1979) proposed the Bootstrap approximation which has been shown to work in wide class of models, specially in case of multiple linear regression. %\sout{Bootstrap inference is an alternative to the asymptotic inference for underlying parameter of interest. In several scenarios, the Bootstrap inference has also been shown to be more accurate than} \sout{the case in linear regression}.
In the last few decades, several variants of Bootstrap have been developed in linear regression. Depending on whether the covariates are non-random or random in linear regression setup, Freedman (1981) proposed the residual Bootstrap or the paired Bootstrap. A few other variants of Bootstrap methods in linear regression setup are the wild Bootstrap (cf. Liu (1988), Mammen (1993)), the weighted Bootstrap (Lahiri (1992), Barbe and Bertail (2012)) and the perturbation Bootstrap (Das and lahiri (2019)). 
%\sout{For the understanding of these Bootstrap methods and their higher order properties see for example Liu (1988), Lahiri (1992), Mammen (1993), Chatterjee and Bose (2005), Das and Lahiri (2019) and references therein.} 
Using similar mechanism of the residual and the paired Bootstrap, Moulton and Zeger (1989, 1991) developed the standardized Pearson residual resampling and the observation vector resampling Bootstrap methods in generalized linear models (GLM). Lee (1990) considered the logistic regression model and showed that the conditional distribution of these resample based Bootstrap estimators for the given data are close to the distribution of the original estimator in almost sure sense. Claeskens et al. (2003) proposed a couple of Bootstrap methods for logistic regression in univariate case, namely `Linear one-step Bootstrap' and `Quadratic one-step Bootstrap'. `Linear one-step Bootstrap' was developed  following the linearization principle proposed in Davison et al. (1986), whereas, `Quadratic one-step Bootstrap' was constructed based on the quadratic approximation of the estimators as discussed in Ghosh (1994). The validity of these two Bootstrap methods for approximating the underlying distribution in almost sure sense was established in Claeskens et al. (2003). They also developed a finite sample bias correction of logistic regression estimator using their quadratic one-step Bootstrap method. 

 %(\gr{Comment: limitation is that these are not SOC. But we should not mention that since we are not showing this.}) \sout{However, to the best of our knowledge, none of the related existing methods including the above-mentioned ones have studied the error rate of the \gr{Bootstrap} approximations} (\rd{what approximation?}).
In order to have an explicit understanding about the sample size requirement for practical implementations of any asymptotically valid method, it is essential to study the error rate of the approximation. The Bootstrap methods in linear regression have been shown to achieve second order correctness (SOC), i.e. having the error rate $o(n^{-1/2})$. In order to draw more accurate inference results compared to that based on asymptotic normal distribution, SOC is essential. An elaborate description on the results on SOC of residual for generalized and perturbation Bootstrap methods in linear regression can be found in Lahiri (1992), Barbe and Bertail (2012) and Das and Lahiri (2019) and references their in. 
% In the same spirit, it is expected for any Bootstrap method to achieve SOC in logistic regression as well. 
However, to the best of our knowledge, for none of the existing Bootstrap methods for logistic regression in the literature, SOC has been explored. In this paper, we propose \underline{Pe}rtur\underline{b}ation \underline{B}ootstrap in \underline{L}ogistic R\underline{e}gression (PEBBLE) as an alternative of the normal approximation approach. Whenever the underlying estimator is a minimizer of certain objective function, perturbation Bootstrap simply produces a Bootstrap version of the estimator by finding the minimizer of a random objective function, suitably developed by perturbing the original objective function using some non-negative random variables. 
%For details of the construction of perturbation Bootstrapped logistic regression estimator, see Section \ref{sec:despb}.}
We show that the perturbation Bootstrap attains SOC in approximating the distribution of $\hat{\bm{\beta}}_n$. For the sake of comparison with the proposed Bootstrap method, we also find the error rate for the normal approximation of the studentized version of the distribution of $\hat{\bm{\beta}}_n$ which comes out to be of $O(n^{-1/2}\log n)$. The extra ``$\log n$'' term in the error rate appears due to the underlying lattice structure. Therefore, the inference based on our Bootstrap method is more accurate than that based on the asymptotic normality. %\sout{One of the critical aspects of the proposed approach remains in choosing the perturbation factors appropriately. Afterwards, (\ref{eqn:def}) needs to be changed based on these perturbation factors and hence the Bootstrap estimator in our case can be defined as the solution of the changed equation.} 
% See section \ref{sec:despb} for details on the construction of the Bootstrap version of the logistic regression estimator.

%$n$ independent copies $G_1^*,\ldots, G_n^*$ of a non-negative and non-degenerate random variable $G^*$ having mean $\mu_{G^*}$, $var(G^*)=\mu_{G^*}^2$ and $\mathbf{E(G^*-\mu_{G^*})^3}=\mu_{G^*}^3$. These characteristics of $G^*$ are assumed to be true throughout this paper. Any additional condition on $G^*$ will be stated in respective theorems. One immediate choice of the distribution of $G^*$ is Beta($1/2, 3/2$). Other choices can be found easily by investigating the generalized beta family of distribution.'')} See Section 2 for details on the construction of the Bootstrap version of the logistic regression estimator.

In order to establish SOC for the proposed method, we start with studentization of $\sqrt{n}(\hat{\bm{\beta}}_n-\boldsymbol{\beta})$ and its perturbation Bootstrap version. We show that unlike in the case of multiple linear regression, here SOC cannot be achieved only by studentization of $\sqrt{n}(\hat{\bm{\beta}}_n-\bm{\beta})$ due to the lattice nature of the distribution of the logistic regression estimator $\hat{\bm{\beta}}_n$, in general. The lattice nature of the distribution is induced by the binary nature of the response variable. It is a common practice to establish SOC by comparing the Edgeworth expansions in original and Bootstrap case (cf. Hall (1992)). However the usual Edgeworth expansion does not exist when the underlying setup is lattice. Therefore, correction terms are required to take care of the lattice nature. For example, one can compare Theorem 20.8 and corollary 23.2 in Bhattacharya and Rao (1986) [hereafter referred to as BR(86)] to learn the correction terms required in the Edgeworth expansions whenever the underlying structure is lattice. In general, these correction terms cannot be approximated with an error of  $o(n^{-1/2})$, which makes SOC unachievable even with studentization. As a remedy we adopt the novel smoothing technique developed in Lahiri (1993). % The smoothing technique was shown to work for residual Bootstrap to achieve SOC in estimating the mean of independent and identical lattice random variables (Lahiri (1993)).
 First, this smoothing technique is applied to transform the lattice nature of the distribution of the studentized version to make it absolutely continuous. Thus the resulting correction terms do not appear in the underlying Edgeworth expansion. Further we use the same smoothing technique for the Bootstrap version and establish SOC by comparing the Edgeworth expansions across the original and the Bootstrap cases. Moreover, an interesting property of the smoothing is that it has negligible effect on the asymptotic variance of $\hat{\bm{\beta}}_n$ and therefore it is not required to incorporate the effect of the smoothing in the form of the studentization. In order to prove the results, we establish the Edgeworth expansion of a smoothed version of a sequence of sample
means of independent random vectors even if they are not identically distributed (cf. Lemma 3). Lemma 3 may be of independent interest for establishing SOC of Bootstrap in other related problems. %To compare the error in Bootstrap approximation with that of normal approximation, we also establish the rate of convergence of the studentized $\hat{\bm{\beta}}_n$ to normality.

The rest of the paper is organized as follows. The perturbation Bootstrap version of the logistic regression estimator is described in Section \ref{sec:despb}. Main results including theoretical properties of the Bootstrap along with normal approximation are stated in Section \ref{sec:main}. In Section \ref{sec:simulation}, finite-sample performance of PEBBLE is evaluated comparing with other related existing methods by simulation experiments. Section \ref{sec:realdata} gives an illustration of PEBBLE in healthcare operations decision dataset.  Auxiliary lemmas and the proof of the theorems are presented in Section \ref{sec:proofs}. Finally, we conclude on the proposed methodology in Section \ref{sec_con}.

% Proof of Lemma \ref{lem:3}, \ref{lem:10} and \ref{lem:11} are relegated to the Supplementary material file. The explicit forms of the Bootstrap confidence intervals are also presented in the Supplementary material file.  

\section{Description of PEBBLE}\label{sec:despb}
In this section, we define the Perturbation Bootstrapped version of the logistic regression estimator. Let $G_1^*,\ldots, G_n^*$ be $n$ independent copies of a non-negative and non-degenerate random variable $G^*$ with mean $\mu_{G^*}$, $Var(G^*)=\mu_{G^*}^2$ and $\mathbf{E}(G^*-\mu_{G^*})^3=\mu_{G^*}^3$. These quantities serve as perturbing random quantities in the construction of the perturbation Bootstrap version of the logistic regression estimator. We define the Bootstrap version as the minimizer of a carefully constructed objective function which involves the observed values $y_1,\dots,y_n$ as well as the estimated probability of successes $\hat{p}(\bm{x}_i)=\dfrac{e^{\bm{x}_i^{\prime}\hat{\bm{\beta}}_n}}{1+e^{\bm{x}_i^{\prime}\hat{\bm{\beta}}_n}}$, $i=1,\dots, n$. Formally, the perturbation Bootstrapped logistic regression estimator $\bm{\hat{\beta}}_n^*$ is defined as
\begin{align*}
\hat{\bm{\beta}}_n^* = \operatorname*{arg\,max}_{\bm{t}}\Bigg[\sum_{i=1}^{n}\Big\{(y_i - \hat{p}(\bm{x}_i))\bm{x}_i^{\prime}\bm{t}\Big\}(G^*_i-\mu_{G^*}) + \mu_{G^*}\sum_{i=1}^{n}\Big\{\hat{p}(\bm{x}_i)(\bm{x}_i^{\prime}\bm{t})-\log (1+e^{\bm{x}_i^{\prime}\bm{t}})\Big\}\Bigg].
\end{align*}
In other words, $\hat{\bm{\beta}}_n^*$ is the solution of the equation
\begin{align}\label{eqn:defpb}
\sum_{i=1}^{n}\big(y_i-\hat{p}(\bm{x}_i)\big)\bm{x}_i(G_i^*-\mu_{G^*})\mu_{G}^{*-1}+\sum_{i=1}^{n}\big(\hat{p}(\bm{x}_i)-p(\bm{t}|\bm{x}_i)\big)\bm{x}_i=0,
\end{align}
since the derivative of the LHS of (\ref{eqn:defpb}) with respect to $\bm{t}$ is negative definite. If Bootstrap equation (\ref{eqn:defpb}) is compared to the original equation (\ref{eqn:def}), it is easy to note that the second part of the LHS of (\ref{eqn:defpb}) is the estimated version of the LHS of (\ref{eqn:def}). The Bootstrap randomness is coming from the first part of the LHS in (\ref{eqn:defpb}), i.e., $\sum_{i=1}^{n}\big(y_i-\hat{p}(\bm{x}_i)\big)\bm{x}_i(G_i^*-\mu_{G^*})\mu_{G}^{*-1}$. Also, the first part is the main contributing term in the asymptotic expansion of the studentized version of $\hat{\bm{\beta}}_n^*$. One immediate choice for the distribution of $G^*$ is Beta($1/2, 3/2$) since the required conditions of $G^*$ are satisfied for this distribution. Other choices can be found in Liu (1988), Mammen (1993) and Das et al. (2019). The moment characteristics of $G^*$ are assumed to be true for the rest of this paper. Further any additional assumption on $G^*$ will be stated in respective theorems.

%\section{Assumptions}
\section{Main Results}\label{sec:main}
In this section, we describe the theoretical results of Bootstrap as well as the normal approximation. In \ref{rate_normal} we state a Berry-Esseen type theorem for a studentized version of the logistic regression estimator $\hat{\bm{\beta}}_n$. 
%\sout{In Theorem \ref{thm:nor} we derive the error rate of the normal approximation which is generally used for inference in logistic regression.}
In \ref{rate_boot} we explore the effectiveness of Bootstrap in approximating the distribution of the studentized version. Theorem \ref{thm:negboot} shows that SOC is not achievable solely by studentization even when $p=1$. As a remedy, we introduce a smoothing in the studentization and show that proposed Bootstrap method achieves SOC.

Before exploring the rate of normal approximation, first we define the class of sets that we would consider in the following theorems. For any natural number $m$, the class of sets $\mathcal{A}_m$ is the collection of Borel subsets of $\mathcal{R}^m$ satisfying 
\begin{equation*}
\sup\limits_{B \in  \mathcal{A}_m}\; \mathbf{\Phi}((\delta B)^{\epsilon}) = O(\epsilon) \;\;\; \text{as}\; \;\epsilon \downarrow 0.
\end{equation*} 
Here $\mathbf{\Phi}$ denotes the normal distribution with mean $\mathbf{0}$ and dispersion matrix being the identity matrix. We are going to use the class $\mathcal{A}_p$ for the uniform asymptotic results on normal and Bootstrap approximations. $\mathbf{P_*}$ denotes the conditional Bootstrap probability of $G^{*}$ given data $\{y_1,\dots,y_n\}$.
\subsection{Rate of Normal Approximation}\label{rate_normal}
%From the definition (\ref{eqn: def}), we have that$\sum_{i=1}^{n}(y_i-\hat{p}(\mathbf{x}_i))\mathbf{x}_i=0$. Now using Taylor's expansion of $\hat{\boldsymbol{\beta}}_n$ around the true value $\boldsymbol{\beta}$, it is easy to see that the asymptotic variance of $\sqrt{n}\hat{\boldsymbol{\beta}}_n$ is $\mathbf{L}_n^{-1}$ where $\mathbf{L}_n=n^{-1}\sum_{i=1}^{n}\mathbf{x}_i\mathbf{x}_i^{\prime}e^{\mathbf{x}_i^{\prime}\boldsymbol{\beta}}(1+e^{\mathbf{x}_i^{\prime}\boldsymbol{\beta}})^{-2}$. Hence the studentized version of $\hat{\boldsymbol{\beta}}_n$ is defined as
%\begin{equation*}\label{eqn:stu1}
%\mathbf{H}_n=\sqrt{n}\hat{\Sigma}_n^{1/2}\big(\hat{\boldsymbol{\beta}}_n-\boldsymbol{\beta}\big),
%\end{equation*} 
%where $\hat{\mathbf{\Sigma}}_n$ is an estimator of $\mathbf{L}_n$.
%An immediate choice of $\hat{\mathbf{\Sigma}}_n$ is $\hat{\mathbf{L}}_n=n^{-1}\sum_{i=1}^{n}\mathbf{x}_i\mathbf{x}_i^{\prime}e^{\mathbf{x}_i^{\prime}\hat{\boldsymbol{\beta}}_n}\big(1+e^{\mathbf{x}_i^{\prime}\hat{\boldsymbol{\beta}}_n}\big)^{-2}$. To find the rate of convergence in normal approximation, we are going to use $\mathbf{H}_n$ wih $\hat{\mathbf{\Sigma}}_n=\hat{\mathbf{L}}_n$; although any $\sqrt{n}-$consistent estimator $\hat{\mathbf{\Sigma}}_n$ of $\mathbf{L}_n$ will work.

In this sub-section we explore the rate of normal approximation of suitable studentized version of the logistic regression estimator $\hat{\bm{\beta}}_n$, uniformly over the class of sets $\mathcal{A}_p$. From the definition (\ref{eqn:def}) of $\hat{\bm{\beta}}_n$, we have that
$\sum_{i=1}^{n}(y_i-\hat{p}(\bm{x}_i))\bm{x}_i=0$. Now using Taylor's expansion of $\sqrt{n}\big(\hat{\boldsymbol{\beta}}_n-\boldsymbol{\beta}\big)$, it is easy to see that the asymptotic variance of $\sqrt{n}\big(\hat{\boldsymbol{\beta}}_n-\boldsymbol{\beta}\big)$ is $\bm{L}_n^{-1}$ where $\bm{L}_n=n^{-1}\sum_{i=1}^{n}\bm{x}_i\bm{x}_i^{\prime}e^{\bm{x}_i^{\prime}\boldsymbol{\beta}}(1+e^{\bm{x}_i^{\prime}\boldsymbol{\beta}})^{-2}$. An estimator of $\bm{L}_n$ can be obtained by replacing $\boldsymbol{\beta}$ by $\hat{\boldsymbol{\beta}}_n$ in the form of $\bm{L}_n$. Hence we can define the studentized version of $\hat{\boldsymbol{\beta}}_n$ as
\begin{equation*}\label{eqn:stu1}
\tilde{\mathbf{H}}_n=\sqrt{n}\hat{\bm{L}}_n^{1/2}\big(\hat{\bm{\beta}}_n-\bm{\beta}\big),
\end{equation*} 
where $\hat{\bm{L}}_n=n^{-1}\sum_{i=1}^{n}\bm{x}_i\bm{x}_i^{\prime}e^{\bm{x}_i^{\prime}\hat{\boldsymbol{\beta}}_n}\big(1+e^{\bm{x}_i^{\prime}\hat{\bm{\beta}}_n}\big)^{-2}$. Other studentized versions can be constructed by considering other estimators of $\bm{L}_n$. For details of the construction of different studentized versions, one can look into Lahiri (1994). The result on normal approximation will hold for other studentized versions also as long as it involves the estimator of $\bm{L}_n$ which is $\sqrt{n}-$consistent.

Berry-Esseen theorem states that the error in normal approximation for the distribution of the mean of a sequence of independent random variables is $O(n^{-1/2})$, provided the average third absolute moment is bounded (cf. Theorem 12.4 in BR(86)). Note that there is an extra multiplicative ``$\log n$'' term besides the usual $n^{-1/2}$ term in the error rate of the normal approximation which is due to the error incurred in Taylor's approximation of $\sqrt{n}(\hat{\boldsymbol{\beta}}_n-\boldsymbol{\beta})$. Since the underlying setup in logistic regression has lattice nature, in general, this error cannot be corrected by higher order approximations, like Edgeworth expansions. Further one important tool in deriving the error rate in normal approximation, and later for deriving the higher order result for the Bootstrap is to find the rate of convergence of $\hat{\boldsymbol{\beta}}_n$ to $\boldsymbol{\beta}$. To this end, we state our first theorem as follows.

\begin{theo}\label{thm:nor}
Suppose $n^{-1}\sum_{i=1}^{n}\|\bm{x}_i\|^3=O(1)$ and $\bm{L}_n\rightarrow \mathbf{L}$ as $n\rightarrow \infty$ where $\bm{L}$ is a pd matrix. Then
\begin{itemize}
\item[\emph{(a)}] there exists a positive constant $C_0$ such that when $n>C_0$ we have
\begin{equation*}
\mathbf{P}\Big(\hat{\boldsymbol{\beta}}_n \; solves  \;(\ref{eqn:def})\; and \;\|\hat{\boldsymbol{\beta}}_n - \boldsymbol{\beta}\|\leq C_0n^{-1/2}(log n)^{1/2}\Big) = 1 - o\big(n^{-1/2}\big).
\end{equation*} 
\item[\emph{(b)}] we have
\begin{align*}
\sup\limits_{\bm{B} \in  \mathcal{A}_p}\big| \mathbf{P}\big(\tilde{\mathbf{H}}_n\in \bm{B}\big)-\mathbf{\Phi}(\bm{B})\big|=O\big(n^{-1/2}\log n\big).
\end{align*}
\end{itemize}
\end{theo}
The proof of Theorem \ref{thm:nor} is presented in Section \ref{sec:proofs}. Theorem \ref{thm:nor} shows that the normal approximation of the distribution of $\tilde{\mathbf{H}}_n$, the studentized logistic regression estimator, has near optimal Berry-Esseen rate. However the rate can be improved significantly by Bootstrap and an application of a smoothing, as described in \ref{rate_boot}.

\subsection{Rate of Bootstrap Approximation}\label{rate_boot}
In this sub-section, we extensively study the rate of Bootstrap approximation for the distribution of the logistic regression estimator. To that end, before exploring the rate of convergence of Bootstrap we need to define the suitable studentized versions in both original and the Bootstrap setting. Similar to the original case, the asymptotic variance of the Bootstrapped logistic regression estimator $\hat{\bm{\beta}}_n^*$ is needed to be found to define the studentized version in the Bootstrap setting. Using Taylor's expansion, from (\ref{eqn:defpb}) it is easy to see that the asymptotic variance of $\sqrt{n}\big(\hat{\bm{\beta}}_n^*-\hat{\bm{\beta}}_n\big)$ is $\hat{\bm{L}}_n^{-1}\hat{\bm{M}}_n\hat{\bm{L}}_n^{-1}$ where $\hat{\bm{L}}_n=n^{-1}\sum_{i=1}^{n}\bm{x}_i\bm{x}_i^{\prime}e^{\bm{x}_i^{\prime}\hat{\bm{\beta}}_n}(1+e^{\bm{x}_i^{\prime}\hat{\bm{\beta}}_n})^{-2}$ and $\hat{\bm{M}}_n=n^{-1}\sum_{i=1}^{n}\big(y_i-\hat{p}(\bm{x}_i)\big)^2\bm{x}_i\bm{x}_i^\prime$. Therefore the studentized version in Bootstrap setting can be defined as 
$$\mathbf{H}_n^*=\sqrt{n}\bm{\hat{M}}_n^{*-1/2}\bm{L}_n^{*}\big(\hat{\bm{\beta}}_n^*-\hat{\bm{\beta}}_n\big),$$
where $\bm{L}_n^*=n^{-1}\sum_{i=1}^{n}\bm{x}_i\bm{x}_i^{\prime}e^{\bm{x}_i^{\prime}\hat{\bm{\beta}}_n^*}\big(1+e^{\bm{x}_i^{\prime}\hat{\bm{\beta}}_n^*}\big)^{-2}$ and $\hat{\bm{M}}_n^*=n^{-1}\sum_{i=1}^{n}\big(y_i-\hat{p}(\bm{x}_i)\big)^2\bm{x}_i\bm{x}_i^\prime\mu_{G^*}^{-2}(G_i^*-\mu_{G^*})^2$. Analogously, we define the original studentized version as $$\mathbf{H}_n=\sqrt{n}\hat{\bm{M}}_n^{-1/2}\hat{\bm{L}}_n\big(\hat{\bm{\beta}}_n-\bm{\beta}\big),$$
which will be used for investigating SOC of Bootstrap for rest of this section. In the next theorem we show that $\mathbf{H}_n^*$ fails to be SOC in approximating the distribution of $\mathbf{H}_n$ even when $p=1$.

%The natural way of defining the Bootstrap version of the original pivot is by replacing $\boldsymbol{\beta}$ by $\hat{\boldsymbol{\beta}}_n$ and $\hat{\boldsymbol{\beta}}_n$ by $\hat{\boldsymbol{\beta}}_n^*$. Therefore, the Bootstrap version of $\mathbf{H}_n$ is defined as 
%$$\mathbf{H}_n^*=\sqrt{n}\mathbf{L}_n^{*1/2}\big(\hat{\boldsymbol{\beta}}_n^*-\hat{\boldsymbol{\beta}}_n\big),$$
%where $\mathbf{L}_n^*=n^{-1}\sum_{i=1}^{n}\mathbf{x}_i\mathbf{x}_i^{\prime}e^{\mathbf{x}_i^{\prime}\hat{\boldsymbol{\beta}}_n^*}\big(1+e^{\mathbf{x}_i^{\prime}\hat{\boldsymbol{\beta}}_n^*}\big)^{-2}$ with $\hat{\boldsymbol{\beta}}_n^*$ being the solution of (\ref{eqn:defpb}). In the next result we show that $\mathbf{H}_n^*$ fails to be SOC in approximating the distribution of $\mathbf{H}_n$ even when $p=1$.

\begin{theo}\label{thm:negboot}
Suppose $p=1$ and denote the only covariate by $x$ in the model (\ref{eqn:model}). Let $x_1,\dots,x_n$ be the observed values of $x$ and $\beta$ be the true value of the regression parameter. Define, $\mu_n=n^{-1}\sum_{i=1}^{n}x_ip(\beta|x_i)$. Assume the following conditions hold:
\begin{enumerate}[label=(C.\arabic*)]
    \item 
$x_1,\dots,x_n$ are non random and are all integers. 
\item $x_{i_1},\dots, x_{i_m} =1$ where $\{i_1,\dots,i_m\}\subseteq \{1,\dots,n\}$ with $m \geq (\log n)^2$.
\item $\max\{|x_i|:i=1,\dots,n\}=O(1)$ and $\liminf_{n\rightarrow \infty}\big[n^{-1}\sum_{i=1}^{n}|x_i|^6\big]>0$.
%\item $\mathbf{L}_n\rightarrow \mathbf{L}$ as $n\rightarrow \infty$ where $\mathbf{L}$ is a pd matrix. (\rd{may be dropped})
\item $\sqrt{n}|\mu_n|< M_1$ for $n\geq M
_1$ where $M_1$ is a positive constant. 
\item The distribution of $G^*$ has an absolutely continuous component with respect to Lebesgue measure and $\mathbf{E}G^{*4}<\infty$.
\end{enumerate} 
Then there exist an interval $\bm{B}_n$ and a positive constant $M_2$ (does not depend on $n$) such that
\begin{align*}
\lim_{n\rightarrow \infty}\mathbf{P}\Big(\sqrt{n}\big|\mathbf{P}_*\big(\mathbf{H}_n^*\in \bm{B}_n\big)- \mathbf{P}\big(\mathbf{H}_n\in \bm{B}_n\big)\big|\geq M_2\Big)=1
\end{align*}
\end{theo}
The proof of Theorem \ref{thm:negboot} is presented in Section \ref{sec:proofs}. Theorem \ref{thm:negboot} shows that unlike in the case of multiple linear regression, in general the Bootstrap cannot achieve SOC even with studentization. Now we further look into the form of the set $\bm{B}_n$. $\bm{B}_n$ is of the form $f_n(\bm{E}_n\times \mathcal{R})$ with $\bm{E}_n=(-\infty,z_n]$ and $z_n=\Big(\dfrac{3}{4n}-\mu_n\Big)$. $f_n(\cdot)$ is a continuous function which is obtained from the Taylor expansion of $\mathbf{H}_n$. Since $\bm{E}_n\times \mathcal{R}$ is a convex subset of $\mathcal{R}^2$, it is also a connected set. Since $f_n(\cdot)$ is a continuous function, $\bm{B}_n$ is a connected subset of $\mathcal{R}$ and hence is an interval. 

Now, we define the smoothed versions of $\mathbf{H}_n$ and $\mathbf{H}_n^*$ which are necessary in achieving SOC by the Bootstrap for general $p$. Note that the primary reason behind Bootstrap's failure is the lattice nature of the distribution of $\sqrt{n}(\hat{\boldsymbol{\beta}}_n-\boldsymbol{\beta})$. Hence if one can somehow smooth the distribution  $\sqrt{n}(\hat{\boldsymbol{\beta}}_n-\boldsymbol{\beta})$, or more generally the distribution of $\mathbf{H}_n$, so that the smoothed version has density with respect to Lebesgue measure, then the Bootstrap may be shown to achieve SOC by employing theory of Edgeworth expansions. To that end, suppose $Z$ is a $p-$dimensional standard normal random vector, independent of $y_1,\dots,y_n$. Define the smoothed version of $\mathbf{H}_n$ as
\begin{align}\label{eqn:mop}
\check{\mathbf{H}}_n=\mathbf{H}_n+\hat{\bm{M}}_n^{-1/2}b_nZ,
\end{align}
where $\{b_n\}_{n\geq 1}$ is a suitable sequence such that it has negligible effect on the variance of $\sqrt{n}(\hat{\boldsymbol{\beta}}_n-\boldsymbol{\beta})$ and hence on the studentization factor. See Theorem \ref{thm:posboot} for the conditions on $\{b_n\}_{n\geq 1}$.
To define the smoothed studentized version in Bootstrap setting, consider another $p-$dimensional standard normal vector by $Z^*$ which is independent of $y_1,\dots,y_n$, $G_1^*,\dots, G_n^*$ and $Z$. Define the smoothed version of $\mathbf{H}_n^*$ as
\begin{align}\label{eqn:msp}
\check{\mathbf{H}}_n^*=\mathbf{H}_n^*+\hat{\bm{M}}_n^{*-1/2}b_nZ^*.
\end{align}
The following theorem can be distinguished as the main theorem of this section as it shows that the smoothing does the trick for Bootstrap to achieve SOC. Thus the inference on $\boldsymbol{\beta}$ based on the Bootstrap after smoothing is much more accurate than the normal approximation. To state the main theorem, define $W_i = \Big(y_i\bm{x}^\prime_i,\big[y_i^2-\mathbf{E}y_i^2\big]\bm{z}^\prime_i\Big)^\prime$ where $y_i=(y_i-p(\bm{\beta}|\bm{x}_i))$ and $z_i=(x_{i1}^2,x_{i1}x_{i2},\dots,x_{i1}x_{ip},x_{i2}^2,x_{i2}x_{i3}$ $,\dots,x_{i2}x_{ip},\dots,x_{ip}^2)^{\prime}$ with $\bm{x}_i=(x_{i1},\dots,x_{ip})^\prime$, $i\in\{1,\dots,n\}$.

\begin{theo}\label{thm:posboot}
Suppose $n^{-1}\sum_{i=1}^{n}\|\mathbf{x}_i\|^6=O(1)$ and the matrix $n^{-1}\sum_{i=1}^{n}Var(W_i)$ converges to some positive definite matrix as $n\rightarrow \infty$. Also choose the sequence $\{b_n\}_{n\geq 1}$ such that $b_n=O(n^{-d})$ and $n^{-1/p_1}\log n = o(b_n^2)$ where $d>0$ is a constant and $p_1=\max\{p+1,4\}$. Then
\begin{itemize}
\item[\emph{(a)}] there exists two positive constant $C_2$ such that when $n>C_2$ we have
\begin{equation*}
\mathbf{P}_*\Big(\hat{\boldsymbol{\beta}}_n^* \; solves  \;(\ref{eqn:defpb})\; \text{and} \;\|\hat{\boldsymbol{\beta}}_n^* - \hat{\boldsymbol{\beta}}_n\|\leq C_2.n^{-1/2}.(log n)^{1/2}\Big) = 1 - o_p\big(n^{-1/2}\big).
\end{equation*} 
\item[\emph{(b)}] we have
\begin{align*}
\sup\limits_{\bm{B} \in  \mathcal{A}_p}\big|\mathbf{P}_*\big(\check{\mathbf{H}}_n^*\in \bm{B}\big)- \mathbf{P}\big(\check{\mathbf{H}}_n\in \bm{B}\big)\big|=o_p\big(n^{-1/2}\big).
\end{align*}
\end{itemize}

\end{theo}
The proof of Theorem \ref{thm:posboot} is presented in Section \ref{sec:proofs}. Theorem \ref{thm:posboot} shows that SOC of PEBBLE can be achieved by a simple smoothing in the studentized pivotal quantities. As a result, much more accurate inference on $\bm{\beta}$ can be drawn based on Bootstrap than that based on normal approximation specially when $n$ is not large enough compared to $p$. The finite sample simulation results presented in Table \ref{sim_table} also confirms this fact.
\begin{rem}
The class of sets $\mathcal{A}_p$ used to state the uniform asymptotic results is somewhat abstract. Note that there are two major reasons behind considering this class. The first reason is to obtain asymptotic normality or to obtain valid Edgeworth expansions for the normalized part of the underlying pivot and the second one is to bound the remainder term by required small magnitude with sufficiently large probability (or Bootstrap probability). A natural choice for $\mathcal{A}$ is the collection of all Borel measurable convex subsets of $\mathcal{R}^p$, due to Theorem 3.1 in BR(86).   
\end{rem}

\begin{rem}
The results on Bootstrap approximation presented in Theorem \ref{thm:posboot}, may be established in almost sure sense also. In that case the only additional requirement is to have $n^{-1}\sum_{i=1}^{n}$ $\|\bm{x}_i\|^{12}=O(1)$, since $y_1,\dots,y_n$ can take either $0$ or $1$. Actually an almost sure version of part (a) of Theorem \ref{thm:posboot} is necessary to establish Theorem \ref{thm:negboot}. Note that the requirement for almost sure version is met under the assumptions of Theorem \ref{thm:negboot}.
\end{rem}

\begin{rem}
Note that the random quantities $Z$ and $Z^*$ respectively, introduced in (\ref{eqn:mop}) and (\ref{eqn:msp}), are essential in achieving SOC of the Bootstrap. $Z$ and $Z^*$ both are assumed to be distributed as $N(\bm{0},\mathbb{I}_p)$, $\mathbb{I}_p$ being the $p\times p$ identity matrix.  However, Theorem \ref{thm:posboot} remains to be true if we replace $\mathbb{I}_p$ by any diagonal matrix, i.e., Theorem \ref{thm:posboot} is true even if we only assume that the components of $Z$ (and of $Z^*$) are independent and have normal distributions. 
\end{rem}

\section{Simulation Study}\label{sec:simulation}
In this section, we compare the performance of PEBBLE with other existing methods via simulation experiments. For comparative study, we consider the Normal approximation, Pearson Residual Resampling Bootstrap (PRRB, Moulton and Zeger (1991)), One-Step Bootstrap (OSB)  and Quadratic Bootstrap (QB) (Claeskens et al.(2003)). We consider $$\boldsymbol{b} = (1,.5,-2,-0.75, 1.5,-1, 1.85, -1.6).$$ Note that $\boldsymbol{b}$ has length 8. For the scenarios where $p \leq 8$, we take the true parameter vector $\boldsymbol{\beta}$ to be the first $p$-many elements of $\boldsymbol{b}$. The covariate vector $X$ is generated from multivariate normal distribution with mean $\mathbf{0}$ and variance $\Sigma = \{\sigma_{ij}\}_{p\times p}$ where $\sigma_{ij} = 0.5^{|i-j|}$. Now, in order to access the performance of all the methods for various dimensional coefficient vectors and sample sizes, we consider the following cases $(n,p) = (30,3), (50,3), (50,4), (100,3),$ $(100,4), (100,6), (200,3), (200,4), (200,6)$ and $(200,8)$. 

In PEBBLE, we take $p_1 = \max\{p+1,4\}$, $b_n^2 = n^{-\frac{1}{p_1+1}}$. Both $Z$ and $Z^*$ are drawn from independent multivariate normal distribution with mean $\mathbf{0}$ and variance $\frac{1}{4}\mathbb{I}_p$. $G_i^*$ is genrated from $Beta(\frac{1}{2},\frac{3}{2})$. Further details regarding the forms of the confidence sets for PEBBLE is provided in the Supplementary Material Section 2. PEBBLE is implemented in \textbf{R}. Other methods namely Normal approximation, PRRB, OSB and QB are also implemented in \textbf{R}. For the experiment, we consider 1000 Bootstrap iterations. In order to find coverage, such experiment is repeated 1000 times for each $(n,p)$ scenario. In Table \ref{sim_table}, we note down the empirical coverage of lower $90\%$ confidence region of $\bm{\beta}$, upper, middle and lower $90\%$ Confidence intervals (CIs) corresponding to the minimum and maximum components of $\bm{\beta}$. We also note down the average over empirical coverages of upper, middle and lower $90\%$ CI corresponding to all components of $\bm{\beta}$. Average widths of $90\%$ CI corresponding to all applicable cases are also noted in parenthesis. It is noted that in general, PEBBLE performs better than other methods; specifically, for lower $n:p$ scenarios (small sample size, high dimension), i.e., cases corresponding to $(n,p) =(30,3), (50,4),(100,6),(200,8)$ in our study. For example, for $(n,p) =(100,6),(200,8)$ it is noted that PEBBLE outperforms other methods by a big margin. As $n$ increases for fixed $p$, performance of PEBBLE is noted to improve and the widths of CIs tend to decrease, as expected. PEBBLE performs better in comparatively bigger margin than other methods. It is also noted that for all the simulation scenarios, the average coverage over all coordinates is much closer to $0.90$ for PEBBLE compared to other methods. We observe that for relatively smaller $n:p$ scenarios, the PEBBLE CIs are a little wider compared to other methods, but, as $n$ increases (for fixed $p$), PEBBLE CI widths become closer to those observed for other methods.  

\begin{table}[]
\def\arraystretch{0.8}
\resizebox{1\columnwidth}{!}{%
\begin{tabular}{llcccccccccc}
\hline
\multicolumn{1}{c}{(n,p)} & Methods & \begin{tabular}[c]{@{}c@{}}$\beta$\\ (lower)\end{tabular} & \begin{tabular}[c]{@{}c@{}}$\beta_{min}$\\ middle (width)\end{tabular} & \begin{tabular}[c]{@{}c@{}}$\beta_{min}$\\ upper\end{tabular} & \begin{tabular}[c]{@{}c@{}}$\beta_{min}$\\ lower\end{tabular} & \begin{tabular}[c]{@{}c@{}}$\beta_{max}$\\ middle (width)\end{tabular} & \begin{tabular}[c]{@{}c@{}}$\beta_{max}$\\ upper\end{tabular} & \begin{tabular}[c]{@{}c@{}}$\beta_{max}$\\ lower\end{tabular} & \begin{tabular}[c]{@{}c@{}}$\beta$ avg.\\ middle (width)\end{tabular} & \begin{tabular}[c]{@{}c@{}}$\beta$ avg.\\ upper\end{tabular} & \begin{tabular}[c]{@{}c@{}}$\beta$ avg.\\ lower\end{tabular} \\ \hline
\multicolumn{1}{c}{\multirow{5}{*}{(30,3)}} & PEBBLE & 0.916 & 0.885 (2.82) & 0.861 & 0.918 & 0.936 (3.95) & 0.928 & 0.914 & 0.900 (3.09) & 0.888 & 0.913 \\
\multicolumn{1}{c}{} & Normal & 0.952 & 0.947 (2.31) & 0.956 & 0.896 & 0.964 (2.86) & 0.909 & 0.993 & 0.958 (2.42) & 0.939 & 0.935 \\
\multicolumn{1}{c}{} & PRRB & 0.946 & 0.916 (2.17) & 0.926 & 0.873 & 0.943 (2.66) & 0.914 & 0.932 & 0.930 (2.27) & 0.915 & 0.905 \\
\multicolumn{1}{c}{} & OSB & 0.953 & 0.942 (2.34) & 0.940 & 0.889 & 0.930 (2.67) & 0.911 & 0.939 & 0.930 (2.38) & 0.916 & 0.921 \\
\multicolumn{1}{c}{} & QB & 0.976 & 0.952 (2.49) & 0.950 & 0.924 & 0.958 (3.07) & 0.936 & 0.965 & 0.936 (2.53) & 0.920 & 0.940 \\ \hline
\multicolumn{1}{c}{\multirow{5}{*}{(50,3)}} & PEBBLE & 0.888 & 0.891 (2.07) & 0.878 & 0.923 & 0.909 (2.89) & 0.925 & 0.895 & 0.904 (2.20) & 0.901 & 0.912 \\
 & Normal & 0.937 & 0.927 (1.76) & 0.924 & 0.901 & 0.948 (2.18) & 0.906 & 0.971 & 0.936 (1.80) & 0.920 & 0.930 \\
 & PRRB & 0.917 & 0.892 (1.68) & 0.896 & 0.880 & 0.912 (2.06) & 0.899 & 0.933 & 0.902 (1.71) & 0.896 & 0.909 \\
 & OSB & 0.925 & 0.911 (1.79) & 0.907 & 0.885 & 0.905 (2.04) & 0.903 & 0.928 & 0.913 (1.77) & 0.908 & 0.910 \\
 & QB & 0.932 & 0.915 (1.86) & 0.904 & 0.913 & 0.916 (2.12) & 0.904 & 0.935 & 0.922 (1.84) & 0.914 & 0.922 \\ \hline
\multicolumn{1}{c}{\multirow{5}{*}{(50,4)}} & PEBBLE & 0.909 & 0.901 (2.92) & 0.877 & 0.936 & 0.909 (3.87) & 0.936 & 0.879 & 0.902 (2.71) & 0.897 & 0.910 \\
\multicolumn{1}{c}{} & Normal & 0.931 & 0.926 (2.14) & 0.952 & 0.902 & 0.951 (2.62) & 0.906 & 0.985 & 0.939 (2.03) & 0.926 & 0.926 \\
\multicolumn{1}{c}{} & PRRB & 0.928 & 0.899 (1.99) & 0.933 & 0.860 & 0.938 (2.42) & 0.899 & 0.949 & 0.906 (1.88) & 0.906 & 0.894 \\
\multicolumn{1}{c}{} & OSB & 0.958 & 0.928 (2.20) & 0.943 & 0.920 & 0.937 (2.44) & 0.908 & 0.952 & 0.928 (2.03) & 0.926 & 0.919 \\
\multicolumn{1}{c}{} & QB & 0.954 & 0.924 (2.11) & 0.931 & 0.915 & 0.926 (2.40) & 0.891 & 0.954 & 0.924 (1.99) & 0.923 & 0.912 \\ \hline
\multirow{5}{*}{(100,3)} & PEBBLE & 0.880 & 0.877 (1.19) & 0.878 & 0.896 & 0.896 (1.69) & 0.912 & 0.891 & 0.887 (1.35) & 0.894 & 0.894 \\
 & Normal & 0.926 & 0.912 (1.08) & 0.909 & 0.904 & 0.918 (1.40) & 0.911 & 0.901 & 0.913 (1.18) & 0.903 & 0.903 \\
 & PRRB & 0.905 & 0.901 (1.08) & 0.907 & 0.901 & 0.912 (1.39) & 0.916 & 0.891 & 0.901 (1.18) & 0.902 & 0.898 \\
 & OSB & 0.906 & 0.897 (1.09) & 0.900 & 0.899 & 0.896 (1.39) & 0.915 & 0.877 & 0.897 (1.18) & 0.900 & 0.894 \\
 & QB & 0.899 & 0.897 (1.08) & 0.889 & 0.900 & 0.880 (1.33) & 0.907 & 0.873 & 0.894 (1.17) & 0.895 & 0.895 \\ \hline
\multirow{5}{*}{(100,4)} & PEBBLE & 0.885 & 0.907 (1.79) & 0.891 & 0.927 & 0.900 (2.24) & 0.920 & 0.880 & 0.898 (1.71) & 0.899 & 0.902 \\
 & Normal & 0.928 & 0.917 (1.39) & 0.924 & 0.903 & 0.942 (1.65) & 0.912 & 0.929 & 0.916 (1.35) & 0.910 & 0.904 \\
 & PRRB & 0.901 & 0.889 (1.35) & 0.892 & 0.900 & 0.896 (1.60) & 0.905 & 0.881 & 0.887 (1.32) & 0.893 & 0.887 \\
 & OSB & 0.915 & 0.904 (1.41) & 0.918 & 0.900 & 0.914 (1.63) & 0.915 & 0.899 & 0.904 (1.36) & 0.906 & 0.900 \\
 & QB & 0.940 & 0.920 (1.49) & 0.934 & 0.902 & 0.943 (1.86) & 0.937 & 0.926 & 0.912 (1.42) & 0.913 & 0.903 \\ \hline
\multicolumn{1}{c}{\multirow{5}{*}{(100,6)}} & PEBBLE & 0.931 & 0.910 (1.77) & 0.880 & 0.917 & 0.907 (2.79) & 0.929 & 0.868 & 0.906 (2.08) & 0.908 & 0.902 \\
\multicolumn{1}{c}{} & Normal & 0.857 & 0.874 (1.23) & 0.883 & 0.871 & 0.903 (1.68) & 0.882 & 0.937 & 0.871 (1.34) & 0.877 & 0.887 \\
\multicolumn{1}{c}{} & PRRB & 0.849 & 0.854 (1.22) & 0.878 & 0.870 & 0.884 (1.66) & 0.869 & 0.914 & 0.848 (1.33) & 0.866 & 0.874 \\
\multicolumn{1}{c}{} & OSB & 0.933 & 0.797 (1.29) & 0.848 & 0.831 & 0.832 (1.66) & 0.845 & 0.872 & 0.791 (1.37) & 0.837 & 0.846 \\
\multicolumn{1}{c}{} & QB & 0.953 & 0.819 (1.37) & 0.865 & 0.838 & 0.863 (1.84) & 0.857 & 0.902 & 0.807 (1.44) & 0.848 & 0.854 \\ \hline
\multirow{5}{*}{(200,3)} & PEBBLE & 0.891 & 0.906 (0.86) & 0.897 & 0.905 & 0.918 (1.21) & 0.908 & 0.915 & 0.903 (1.01) & 0.896 & 0.906 \\
 & Normal & 0.905 & 0.904 (0.78) & 0.902 & 0.910 & 0.910 (1.03) & 0.936 & 0.879 & 0.902 (0.89) & 0.912 & 0.894 \\
 & PRRB & 0.902 & 0.900 (0.77) & 0.896 & 0.904 & 0.899 (1.02) & 0.930 & 0.874 & 0.893 (0.88) & 0.904 & 0.892 \\
 & OSB & 0.905 & 0.902 (0.78) & 0.900 & 0.917 & 0.897 (1.01) & 0.935 & 0.870 & 0.895 (0.88) & 0.910 & 0.893 \\
 & QB & 0.867 & 0.890 (0.75) & 0.889 & 0.913 & 0.868 (0.93) & 0.924 & 0.842 & 0.871 (0.83) & 0.893 & 0.877 \\ \hline
\multirow{5}{*}{(200,4)} & PEBBLE & 0.872 & 0.898 (1.08) & 0.890 & 0.908 & 0.912 (1.54) & 0.922 & 0.893 & 0.900 (1.11) & 0.900 & 0.905 \\
 & Normal & 0.919 & 0.917 (0.89) & 0.902 & 0.917 & 0.910 (1.18) & 0.918 & 0.893 & 0.906 (0.92) & 0.905 & 0.902 \\
 & PRRB & 0.899 & 0.908 (0.88) & 0.891 & 0.915 & 0.891 (1.15) & 0.916 & 0.876 & 0.892 (0.91) & 0.896 & 0.890 \\
 & OSB & 0.905 & 0.911 (0.89) & 0.897 & 0.914 & 0.901 (1.16) & 0.925 & 0.880 & 0.900 (0.92) & 0.905 & 0.898 \\
 & QB & 0.926 & 0.924 (0.93) & 0.905 & 0.923 & 0.921 (1.23) & 0.930 & 0.892 & 0.917 (0.97) & 0.912 & 0.907 \\ \hline
\multirow{5}{*}{(200,6)} & PEBBLE & 0.927 & 0.915 (1.32) & 0.890 & 0.930 & 0.921 (1.79) & 0.933 & 0.875 & 0.913 (1.59) & 0.908 & 0.906 \\
 & Normal & 0.794 & 0.833 (0.89) & 0.855 & 0.868 & 0.892 (1.17) & 0.915 & 0.862 & 0.847 (1.01) & 0.863 & 0.868 \\
 & PRRB & 0.791 & 0.829 (0.90) & 0.860 & 0.872 & 0.872 (1.18) & 0.911 & 0.859 & 0.840 (1.02) & 0.860 & 0.865 \\
 & OSB & 0.904 & 0.751 (0.92) & 0.813 & 0.842 & 0.794 (1.18) & 0.893 & 0.780 & 0.741 (1.03) & 0.814 & 0.818 \\
 & QB & 0.902 & 0.738 (0.88) & 0.804 & 0.837 & 0.784 (1.15) & 0.893 & 0.768 & 0.736 (1.01) & 0.814 & 0.814 \\ \hline
\multicolumn{1}{c}{\multirow{5}{*}{(200,8)}} & PEBBLE & 0.841 & 0.869 (1.75) & 0.837 & 0.948 & 0.866 (2.28) & 0.965 & 0.776 & 0.851 (1.94) & 0.866 & 0.877 \\
\multicolumn{1}{c}{} & Normal & 0.405 & 0.679 (0.94) & 0.886 & 0.676 & 0.734 (1.19) & 0.696 & 0.961 & 0.688 (1.00) & 0.778 & 0.800 \\
\multicolumn{1}{c}{} & PRRB & 0.496 & 0.679 (0.98) & 0.887 & 0.673 & 0.731 (1.23) & 0.701 & 0.953 & 0.691 (1.03) & 0.780 & 0.803 \\
\multicolumn{1}{c}{} & OSB & 0.861 & 0.468 (0.97) & 0.810 & 0.571 & 0.569 (1.17) & 0.634 & 0.843 & 0.486 (1.00) & 0.680 & 0.714 \\
\multicolumn{1}{c}{} & QB & 0.852 & 0.470 (0.98) & 0.805 & 0.575 & 0.551 (1.15) & 0.637 & 0.837 & 0.480 (0.99) & 0.680 & 0.713 \\ \hline
\end{tabular}}
\caption{Comparative performance study of the proposed method \underline{Pe}rtur\underline{b}ation \underline{B}ootstrap in \underline{L}ogistic R\underline{e}gression (PEBBLE) and other existing methods Normal approximation (Normal), Pearson Residual Resampling Bootstrap (PRRB), One-Step Bootstrap (OSB) and Quadratic Bootstrap (QB). All considered coverage analysis is based on $90\%$ confidence intervals (CI) and average is noted over 1000 experiments, results for each experiment is evaluated based on 1000 Bootstrap iterations. We consider the average coverages based on lower CI of norm of $\bm{\beta}$ (column 1), upper, lower and middle CI of the minimum absolute value of $\bm{\beta}$ (column 2,3,4), upper, lower and middle CI of the maximum absolute value of the $\bm{\beta}$ (column 5,6,7), upper, lower and middle CI of the all components of $\bm{\beta}$, on average (column 8,9,10). The average width of the middle CI corresponding to the min, max and average components are provided in parenthesis in columns 2,5,8 respectively.}
\label{sim_table}
\end{table}

\section{Application to Healthcare Operations Decision}\label{sec:realdata}
Vaginal delivery is the most common type of birth. However due to several medical reasons, with advancement of medical procedures, caesarian delivery is often considered as an alternative way for delivery. Recently a few studies showed how the recommended type of delivery may depend on various clinical aspects of the mother including age, blood pressure and heart problem (Rydahl et al. (2019), Amorim et al. (2017), Pieper (2012)). We consider a dataset about caesarian section results of 80 pregnant women along with several important related clinical covariates. The dataset is avialable in the following link \footnote{\url{https://archive.ics.uci.edu/ml/datasets/Caesarian+Section+Classification+Dataset}}. We regress the type of delivery (caesarian or not) on several related covariates namely age, delivery number, delivery time, blood pressure and presence of heart problem. Delivery time can take three values 0 (timely), 1 (premature) and 2 (latecomer). Blood pressure is denoted by 0, 1, 2 for the cases low, normal and high respectively. The covariate presence of heart problem is also binary, 0 denoting apt behaviour and 1 denoting its inept condition. We perform a logistic regression and corresponding CIs are computed using PEBBLE and in Table \ref{app_table} we note down the results. It is noted that although $90\%$ CIs for all the covariates contain zero, however, the $90\%$ CI for heart problem belong to the positive quadrant mostly; also the upper $90\%$ CI completely belongs to the positive quadrant, which implies women with heart problems tend to have caesarian procedure, coinciding with the findings in Yap et al. (2008) and Blaci et al. (2011).

\begin{table}[]\def\arraystretch{1}
\resizebox{0.7\columnwidth}{!}{%
\begin{tabular}{lcccc}
\hline
Variables & $\hat{\beta}$ & \begin{tabular}[c]{@{}c@{}}$90\%$ CI\\ (mid)\end{tabular} & \begin{tabular}[c]{@{}c@{}}$90\%$ CI\\ (upper)\end{tabular} & \begin{tabular}[c]{@{}c@{}}$90\%$ CI\\ (lower)\end{tabular} \\ \hline
%Intercept & 0.107 & (-8.568, 4.264) & \textgreater -7.199 & \textless 2.804 \\
Age & -0.010 & (-0.151, 0.300) & \textgreater -0.100 & \textless 0.237 \\
Delivery number & 0.263 & (-0.544, 0.740) & \textgreater -0.398 & \textless 0.601 \\
Delivery time & -0.427 & (-0.643, 0.466) & \textgreater -0.521 & \textless 0.348 \\
Blood pressure & -0.251 & (-0.709, 0.680) & \textgreater -0.548 & \textless 0.531 \\
Heart problem & 1.702 & (-0.139, 2.327) & \textgreater 0.145 & \textless 2.105 \\ \hline
\end{tabular}}
\caption{Real Data Analysis : The estimated coefficients and corresponding middle, upper and lower $90\%$ CIs are noted for all the covariates; the type of delivery is the dependent variable, which takes values 1 or 0 based on if the delivery was caesarian or not.}
\label{app_table}
\end{table}

\section{Proof of the Results} \label{sec:proofs}

\subsection{Notations}
Before going to the proofs we are going to define few notations. Suppose, $\mathbf{\Phi_V}$ and $\mathbf{\phi_V}$ respectively denote the normal distribution and its density with mean $\mathbf{0}$ and covariance matrix $\bm{V}$. We will write $\mathbf{\Phi_V} = \mathbf{\Phi}$ and $\phi_{\mathbf{V}}=\phi$ when the dispersion matrix $\mathbf{V}$ is the identity matrix.
$C, C_1, C_2,\cdots$ denote generic constants that do not depend on the variables like $n, x$, and so on. $\bm{\nu}_1$, $\bm{\nu}_2$ denote vectors in $\mathscr{R}^p$, sometimes with some specific structures (as mentioned in the proofs). $(\mathbf{e_1},\ldots,\mathbf{e_p})'$ denote the standard basis of $\mathcal{R}^p$. For a non-negative integral vector $\bm{\alpha} = (\alpha_1, \alpha_2,\ldots,\alpha_l)'$ and a function $f = (f_1,f_2,\ldots,f_l):\ \mathcal{R}^l\ \rightarrow \ \mathcal{R}^l$, $l\geq 1$, let $|\bm{\alpha}| = \alpha_1 +\ldots+ \alpha_l$, $\bm{\alpha}! = \alpha_1!\ldots \alpha_l!$, $f^{\bm{\alpha}} = (f_1^{\alpha_1})\ldots(f_l^{\alpha_l})$, $D^{\bm{\alpha}}f_1 = D_1^{\alpha_1}\cdots D_l^{\alpha_l}f_1$, where $D_jf_1$ denotes the partial derivative of $f_1$ with respect to the $j$th component of $\bm{\alpha}$, $1\leq j \leq l$. We will write $D^{\bm{\alpha}}=D$ if $\bm{\alpha}$ has all the component equal to 1. For $\mathbf{t} =(t_1,t_2,\cdots t_l)'\in \mathcal{R}^l$ and $\mathbf{\alpha}$ as above, define $\mathbf{t}^{\bm{\alpha}} = t_1^{\alpha_1}\cdots t_l^{\alpha_l}$. For any two vectors $\bm{\alpha}, \bm{\beta} \in \mathcal{R}^k$, $\bm{\alpha} \leq \bm{\beta}$ means that each of the component of $\bm{\alpha}$ is smaller than that of $\bm{\beta}$. For a set $\bm{A}$ and real constants $a_1,a_2$, $a_1\bm{A}+a_2=\{a_1y+a_2:y\in \bm{A}\}$, $\partial A$ is the boundary of $A$ and $A^{\epsilon}$ denotes the $\epsilon-$neighbourhood of $A$ for any $\epsilon>0$. $\mathcal{N}$ is the set of natural numbers. $C(\cdot),C_1(\cdot),\dots$ denote generic constants which depend on only their arguments. Given two probability measures $P_1$ and $P_2$ defined on the same space $(\Omega,\mathcal{F})$, $P_1*P_2$ defines the measure on $(\Omega,\mathcal{F})$ by convolution of $P_1$ \& $P_2$ and $\|P_1-P_2\|=|P_1-P_2|(\Omega)$, $|P_1-P_2|$ being the total variation of $(P_1-P_2)$. For a function $g:\mathcal{R}^k\rightarrow \mathcal{R}^m$ with $g=(g_1,\dots,g_m)^\prime$, $$Grad [g(\bm{x})]=\Big(\Big(\frac{\partial g_i(\bm{x})}{\partial x_j}\Big)\Big)_{m\times k}.$$

%Before moving to the proofs of the theorems, we are going to state few essential Lemmas. Also we are going to present the proof of Theorem \ref{thm:negboot} at last, since some proof steps of Theorem \ref{thm:posboot} will be essential in proving Theorem \ref{thm:negboot}. For Lemma \ref{lem:2} below, define $\xi_{1,n,s}(\bm{t})=\Big(1+\sum_{i=1}^{s-2}n^{-r/2}\tilde{P}_r\big(i\bm{t}:\{\bar{\chi}_{\nu,n}\}\big)\Big)\exp\Big\{-\bm{t}^\prime \bm{E}_n \bm{t}/2\Big\}$ where $\bm{E}_n=n^{-1}\sum_{i=1}^{n}Var(Y_i)$ and $\bar{\chi}_{\nu,n}$ is the average $\nu$th cumulant of $Y_1,\dots,Y_n$. Define $\bar{\rho}_l=n^{-1}\sum_{i=1}^{n}$ $\mathbf{E}\|Y_i\|^{l}$, the average $l$th absolute moment of $\{Y_1,\dots,Y_n\}$. The polynomials $\tilde{P}_r\big(\bm{z}:\{\bar{\chi}_{\nu,n}\}\big)$ are defined on the pages of $51-53$ of Bhattacharya and Rao (1986). Define the identity $$\xi_{1,n,s}(\bm{t})\Big(\sum_{j=0}^{\infty}(-\|\bm{t}\|^2b_n^2)^j/j!\Big)=\xi_{n,s}(\bm{t})+ o\big(n^{-(s-2)/2}\big),$$ uniformly in $\|\bm{t}\|<1$, where $c_n$ is defined in Lemma \ref{lem:2}. $\psi_{n,s}(\cdot)$ is the Fourier inverse of $\xi_{n,s}(\cdot)$.

Before moving to the proofs of the main theorems, we state some auxiliary lemmas. The proofs of lemma \ref{lem:3}, \ref{lem:10} and \ref{lem:11} are relegated to the Supplementary material file to save space. Also we are going to present the proof of Theorem \ref{thm:negboot} at last, since some proof steps of Theorem \ref{thm:posboot} will be essential in proving Theorem \ref{thm:negboot}.

\subsection{Auxiliary Lemmas}
\begin{lem}\label{lem:1}
Suppose $Y_1,\dots,Y_n$ are zero mean independent r.v.s with $\mathbf{E}(|Y_i|^t)< \infty$ for $i = 1,\dots,n$ and $S_n = \sum_{i = 1}^{n}Y_i$. Let $\sum_{i = 1}^{n}\mathbf{E}(|Y_i|^t) = \sigma_t$, $c_t^{(1)}=\big(1+\frac{2}{t}\big)^t$ and $c_t^{(2)}=2(2+t)^{-1}e^{-t}$. Then, for any $t\geq 2$ and $x>0$,
\begin{equation*}
P[|S_n|>x]\leq c_t^{(1)}\sigma_t x^{-t} + exp(-c_t^{(2)}x^2/\sigma_2)
\end{equation*}
\end{lem}
%Proof of Lemma 8.1. 
\textbf{Proof of Lemma \ref{lem:1}}.
This inequality was proved in Fuk and Nagaev (1971).

\begin{lem}\label{lem:2}
For any $t>0$, $\dfrac{1-N(t)}{n(t)}\leq \dfrac{1}{t}$ wher $N(\cdot)$ and $n(\cdot)$ respectively denote the cdf and pdf of real valued standard normal rv. 
\end{lem}
Proof of Lemma \ref{lem:2}: This inequality is proved in Birnbaum (1942).\\

\begin{lem}\label{lem:3}
Suppose $Y_1,\dots,Y_n$ are mean zero independent random vectors in $\mathcal{R}^k$ with $\bm{E}_n=n^{-1}$ $ \sum_{i=1}^{n}Var(Y_i)$ converging to some positive definite matrix $V$. Let $s\geq 3$ be an integer and $\bar{\rho}_{s
+\delta} =O(1)$ for some $\delta>0$. Additionally assume $Z$ to be a $N(\bm{0},\bm{I}_k)$ random vector which is independent of $Y_1,\dots,Y_n$ and the sequence $\{c_n\}_{n\geq 1}$ to be such that $c_n=O(n^{-d})$ \& $n^{-(s-2)/\tilde{k}}\log n = o(c_n^2)$ where $\tilde{k}=\max\{k+1,s+1\}$ \& $d>0$ is a constant. Then for any Borel set $B$ of $\mathcal{R}^k$,
\begin{align}\label{eqn:1}
\Big|\mathbf{P}\big(\sqrt{n}\bar{Y}+c_nZ\in B\big)- \int_{B}\psi_{n,s}(x)dx\Big|=o\Big(n^{-(s-2)/2}\Big),
\end{align}
where $\psi_{n,s}(\cdot)$ is defined above.
\end{lem}

Proof of Lemma \ref{lem:3}. See Section 1 of Supplementary material file.

\begin{lem}\label{lem:4}
Suppose all the assumptions of Lemma \ref{lem:2} are true. Define $d_n=n^{-1/2}c_n$ and $A_{\delta}=\{x\in \mathcal{R}^k:\|\bm{x}\|<\delta\}$ for some $\delta>0$. Let $H:\mathcal{R}^k\rightarrow \mathcal{R}^m$ ($k\geq m$) has continuous partial derivatives of all orders on $A_{\delta}$ and $Grad[H(\bm{0})]$ is of full row rank. Then for any Borel set $B$ of $\mathcal{R}^m$ we have
\begin{align}\label{eqn:20}
\Big|\mathbf{P}\Big(\sqrt{n}\big(H(\bar{Y}_n+d_nZ)-H(\bm{0}))\in B\Big)- \int_{B}\check{\psi}_{n,s}(x)dx\Big|=o\Big(n^{-(s-2)/2}\Big),
\end{align}
where $\check{\psi}_{n,s}(\bm{x})=\Big[1+\sum_{r=1}^{s-2}n^{-r/2}a_{1,r}(Q_n,\bm{x})\phi_{\check{\bm{M}}_n}(\bm{x})\Big]\Big[\sum_{j=1}^{m_1-1}c_n^{2j}a_{2,j}(\bm{x})\Big]$ with $m_1=\inf \big\{j:c_n^{2j}=o\big(n^{-(s-2)/2}\big)\big\}$ and $Q_n$ being the distribution of $\sqrt{n}\bar{Y}_n$. $a_{1,r}(Q_n,\cdot)$, $r\in \{1,\dots,(s-2)\}$, are polynomials whose coefficients are continuous functions of first $s$ average cumulants of $\{Y_1,\dots, Y_n\}$. $a_{2,j}(\cdot)$, $j\in\{1,\dots,(m-1)\}$, are polynomials whose coefficients are continuous functions of partial derivatives of $H$ of order $(s-1)$ or less. $\check{\bm{M}}_n=\bar{\bm{B}}\bm{E}_n\bar{\bm{B}}^\prime$ with $\bar{\bm{B}}=
Grad [H(\bm{0})]$ and $\bm{E}_n=n^{-1}\sum_{i=1}^{n}Var(Y_i)$.
\end{lem}
Proof of Lemma \ref{lem:4}. This follows exactly through the same line of the proof of Lemma 3.2 in Lahiri (1989).

\begin{lem}\label{lem:5}
Let $Y_1,\dots,Y_n$ be mean zero independent random vectors in $\mathcal{R}^k$ with $n^{-1}\sum_{i=1}^{n}E\|Y_i\|^3 =O(1)$. Suppose $T_n^2=E_n^{-1}$ where $E_n=n^{-1}\sum_{i=1}^{n}Var(Y_i)$ is the average positive definite covariance matrix and $E_n$ converges to some positive definite matrix as $n\rightarrow \infty$. Then for any Borel subset $B$ of $\mathcal{R}^k$ we have
\begin{align*}
\Big| \mathbf{P}\Big(n^{-1/2}T_n\sum_{i=1}^{n}Y_i\in B\Big)-\mathbf{\Phi}(B)\Big|\leq C_{22}(k)n^{-1/2}\rho_3 + 2\;\mathbf{\Phi}\Big((\partial B)^{C_{22}(k)\rho_3 n^{-1/2}}\Big),
\end{align*}
where $\rho_3=n^{-1}\sum_{i=1}^{n}E\|T_nY_i\|^3$.
\end{lem}
Proof of Lemma \ref{lem:5}. This is a direct consequence of part (a) of corollary 24.3 in BR(86).

\begin{lem}\label{lem:6}
Suppose $A,B$ are matrices such that $(A-aI)$ and $(B-aI)$ are positive semi-definite matrices of same order, for some $a>0$. For some $r>0$, $A^r, B^r$ are defined in the usual way. Then for any $0<r<1$, we have $$\|A^r-B^r\|\leq ra^{r-1}\|A-B\|.$$
\end{lem}
Proof of Lemma \ref{lem:6}. More general version of this lemma is stated as corollary (X.46) in Bhatia (1996).

\begin{lem}\label{lem:7}
Suppose all the assumptions of Lemma \ref{lem:4} are true and $\check{\bm{M}}_n=I_m$, the $m\times m$ identity matrix. Define $\hat{H}_n = \Big[\sqrt{n}\big(H(\bar{Y}_n+d_nZ)-H(\bm{0}))\Big] + R_n$ where $\mathbf{P}\Big(\|R_n\|=o\big(n^{-(s-2)/2}\big)\Big)=1-o\big(n^{-(s-2)/2}\big)$ and $s$ is as defined in Lemma \ref{lem:3}.  Then we have
\begin{align}\label{eqn:21}
\sup_{B\in \mathcal{A}_m}\Big|\mathbf{P}\Big(\hat{H}_n\in B\Big)- \int_{B}\check{\psi}_{n,s}(x)dx\Big|=o\Big(n^{-(s-2)/2}\Big),
\end{align}
where the class of sets $\mathcal{A}_m$ is as defined in section \ref{sec:main}.
\end{lem}
Proof of Lemma \ref{lem:7}. Recall the definition of $(\partial B)^{\epsilon}$ which is given in section \ref{sec:main}. For some $B\subseteq \mathcal{R}^m$ and $\delta>0$, define $B^{n,s, \delta}=(\partial B)^{\delta n^{-(s-2)/2}}$. Hence using Lemma \ref{lem:4}, for any $B\in \mathcal{A}_m$ we have
\begin{align}\label{eqn:22}
&\Big|\mathbf{P}\Big(\hat{H}_n\in B\Big)- \int_{B}\check{\psi}_{n,s}(x)dx\Big|=o\Big(n^{-(s-2)/2}\Big)\nonumber \\
\leq\; & \Big|\mathbf{P}\Big(\hat{H}_n\in B\Big)- \mathbf{P}\Big(\sqrt{n}\big(H(\bar{Y}_n+d_nZ)-H(\bm{0}))\in B\Big)\Big| + o\Big(n^{-(s-2)/2}\Big)\nonumber\\
\leq \;& \mathbf{P}\Big(\|R_n\|\neq o\big(n^{-(s-2)/2}\big)\Big) + 2\mathbf{P}\Big(\sqrt{n}\big(H(\bar{Y}_n+d_nZ)-H(\bm{0}))\in B^{n,s, \delta}\Big)+ o\Big(n^{-(s-2)/2}\Big)\nonumber\\
= \; & 2\mathbf{P}\Big(\sqrt{n}\big(H(\bar{Y}_n+d_nZ)-H(\bm{0}))\in B^{n,s, \delta}\Big)+ o\Big(n^{-(s-2)/2}\Big)\nonumber\\
=\; & 2\int_{B^{n,s,\delta}} \check{\psi}_{n,s}(x)dx + o\Big(n^{-(s-2)/2}\Big)
\end{align}
for any $\delta>0$. Now calculations at page 213 of BR(86) and arguments at page 58 of Lahiri(1989) imply that for any $B\in \mathcal{A}_m$, $$\int_{B^{n,s,\delta}} \check{\psi}_{n,s}(x)dx \leq C_{21}(s)\sup\limits_{B \in  \mathcal{A}_m}\; \mathbf{\Phi}\big(B^{n,s, \delta}\big)+o\Big(n^{-(s-2)/2}\Big)=o\Big(n^{-(s-2)/2}\Big),$$ since $\delta >0$ is arbitrary. Therefore (\ref{eqn:21}) follows from (\ref{eqn:22}).\\

\begin{lem}\label{lem:8}
Let $\bm{A}$ and $\bm{B}$ be positive definite matrices of same order. Then for some given matrix $\bm{C}$, the solution of the equation $\bm{AX}+\bm{XB}=\bm{C}$ can be expressed as $$\bm{X}=\int_{0}^{\infty}e^{-t\bm{A}}\bm{C}e^{-t\bm{B}}dt,$$
where $e^{-t\bm{A}}$ and $e^{-t\bm{B}}$ are defined in the usual way.
\end{lem}
Proof of Lemma \ref{lem:8}. This lemma is an easy consequence of Theorem VII.2.3 in Bhatia (1996).\\

\begin{lem}\label{lem:9}
Let $W_1,\dots,W_n$ be n independent mean $0$ random variables with average variance $s_n^2=n^{-1}\sum_{i=1}^{n}\mathbf{E}W_i^2$ and $\mathbf{P}\big(\max\{|W_j|:i \in \{1,\dots,n\}\} \leq C_{30}\big)=1$ for some positive constant $C_{30}$ and integer $s\geq 3$. $\bar{\chi}_{\nu,n}$ is the average $\nu$th cumulant. Recall the polynomial $\tilde{P}_r$ for any non-negative integer $r$, as defined in the beginning of this section. Then there exists two constants $0<C_{31}(s)<1$ and $C_{32}(s)>0$ such that whenever $|t|\leq C_{31}(s)\sqrt{n}\min\{C_{30}^{-2}s_n, C_{30}^{-s/(s-2)}s_n^{s/(s-2)}\}$, we have $$\Big|\prod_{j=1}^{n}\mathbf{E}\Big(e^{in^{-1/2}tW_j}\Big)-\sum_{r=0}^{s-2}n^{-r/2}\tilde{P}_r\big(it:\{\bar{\chi}_{\nu,n}\}\big)e^{-(s_n^2t^2)/2}\Big|\leq C_{32}(s)C_{30}^{s}s_n^{-s} n^{-(s-2)/2}\Big[(s_nt)^s + (s_n t)^{3(s-2)}\Big]e^{-(s_n^2t^2)/4}$$
\end{lem}

Proof of Lemma \ref{lem:9}. In view of Theorem 9.9 of BR(86), enough to show that for any $j \in \{1,\dots,n\}$, $\Big|\mathbf{E}\big(e^{its_n^{-1}n^{-1/2}W_j}\big)-1\Big|\leq 1/2$ whenever $|t|\leq C_{31}(s)\sqrt{n}\min\{C_{30}^{-2}s_n, C_{30}^{-s/(s-2)}s_n^{s/(s-2)}\}$. This is indeed the case due to the fact that $$\Big|\mathbf{E}\big(e^{itn^{-1/2}W_j}\big)-1\Big|\leq \dfrac{t^2\mathbf{E}W_j^2}{2n s_n^2}.$$

\begin{lem}\label{lem:10}
Assume the setup of Theorem \ref{thm:negboot} and let $X_i=y_ix_i$, $i\in \{1,\dots,n\}$. Define $\sigma_n^2=n^{-1}\sum_{i=1}^{n}Var(X_i)$ and $\bar{\chi}_{\nu,n}$ as the $\nu$th average cumulant of $\{(X_1-E(X_1)),\dots, (X_n-E(X_n))\}$.  $P_r\big(-\Phi_{\sigma_n^2}:\{\bar{\chi}_{\nu,n}\}\big)$ is the finite signed measure on $\mathcal{R}$ whose density is $\tilde{P}_r\big(-D: \{\bar{\chi}_{\nu,n}\}\big)\phi_{\sigma_n^2}(x)$. Let $S_0(x)=1$ and $S_1(x)=x-1/2$. Suppose $\sigma_n^2$ is bounded away from both $0$ \& $\infty$ and assumptions (C.1)-(C.3) of Theorem \ref{thm:negboot} hold. Then we have
\begin{align}\label{eqn:23}
\sup_{x\in \mathcal{R}}\Big|&\mathbf{P}\Big(n^{-1/2}\sum_{i=1}^{n}\big(X_i-E(X_i)\big)\leq x\Big)- \sum_{r=0}^{1}n^{-r/2}(-1)^rS_r(n\mu_n+n^{1/2}x)\dfrac{d^r}{dx^r}\Phi_{\sigma_n^2}(x)\nonumber\\
&-n^{-1/2}P_1\big(-\Phi_{\sigma_n^2}:\{\bar{\chi}_{\nu,n}\}\big)(x)\Big| = o\big(n^{-1/2}\big),
\end{align}
where $P_r\big(-\Phi_{\sigma_n^2}:\{\bar{\chi}_{\nu,n}\}\big)(x)$ is the $P_r\big(-\Phi_{\sigma_n^2}:\{\bar{\chi}_{\nu,n}\}\big)-$measure of the set $(-\infty,x]$.

\end{lem}

Proof of Lemma \ref{lem:10}. See Section 1 in the Supplementary material file.

\begin{lem}\label{lem:11}
Let $\breve{W}_1,\dots,\breve{W}_n$ be iid mean $\bm{0}$ non-degenerate random vectors in $\mathcal{R}^{l}$ for some natural number $l$, with finite fourth absolute moment and $\limsup_{\|\bm{t}\|\rightarrow \infty}\big|\mathbf{E}e^{i\bm{t}^\prime\breve{W}_1}\big|<1$ (i.e. Cramer's condition holds). Suppose $\breve{W}_i=(\breve{W}_{i1}^{\prime},\dots,\breve{W}_{im}^{\prime})^\prime$ where $\breve{W}_{ij}$ is a random vector in $\mathcal{R}^{l_j}$ and $\sum_{j=1}^{m}l_j=l$, $m$ being a fixed natural number. Consider the sequence of random variables $\tilde{W}_1,\dots,\tilde{W}_n$ where $\tilde{W}_i=(c_{i1}\breve{W}_{i1}^\prime,\dots,c_{im}\breve{W}_{im}^\prime)^\prime$. $\{c_{ij}:i\in \{1,\dots,n\}, j \in \{1,\dots,m\}\}$ is a collection of real numbers such that for any $j\in \{1,\dots,m\}$, $\Big\{n^{-1}\sum_{i=1}^{n}|c_{ij}|^4\Big\}=O(1)$ and $\liminf_{n\rightarrow \infty}n^{-1}\sum_{i=1}^{n}c_{ij}^2 > 0$. Also assume that $\tilde{\bm{V}}_n = Var(\tilde{W}_i)$ converges to some positive definite matrix and $\bar{\chi}_{\nu,n}$ denotes the average $\nu$th cumulant of $\tilde{W}_1,\dots,\tilde{W}_n$. Then we have
\begin{align}\label{eqn:29}
\sup_{\bm{B}\in \mathcal{A}_l}\Big|\mathbf{P}\Big(n^{-1/2}\sum_{i=1}^{n}\tilde{W}_i \in \bm{B}\Big)-\int_{\bm{B}}\Big[1+n^{-1/2}\tilde{P}_r\big(-D: \{\bar{\chi}_{\nu,n}\}\big)\Big]\phi_{\tilde{\bm{V}}_n}(\bm{t})d\bm{t}\Big|=o\big(n^{-1/2}\big),
\end{align}
where the collection of sets $\mathcal{A}_l$ is as defined in section \ref{sec:main}.
\end{lem}

Proof of Lemma \ref{lem:11}. See Section 1 in the Supplementary material file.

\subsection{Proof of Main Results}
Proof of Theorem \ref{thm:nor}. Recall that the studentized pivot is
\begin{align*}
\tilde{\mathbf{H}}_n=\sqrt{n}\hat{\bm{L}}_n^{1/2}\big(\hat{\bm{\beta}}_n-\bm{\beta}\big),
\end{align*}
where $\hat{\bm{L}}_n=n^{-1}\sum_{i=1}^{n}\bm{x}_i\bm{x}_i^\prime e^{\bm{x}_i^\prime \hat{\bm{\beta}}_n}\big(1+e^{\bm{x}_i^\prime}\hat{\bm{\beta}}_n\big)^{-2}$. $\hat{\bm{\beta}}_n$ is the solution of (\ref{eqn:def}). By Taylor's theorem, from (\ref{eqn:def}) we have
\begin{align}\label{eqn:1.1}
\bm{L}_n\big(\hat{\bm{\beta}}_n-\bm{\beta}\big)=n^{-1}\sum_{i=1}^{n}(y_i-p(\bm{\beta}|\bm{x}_i))\bm{x}_i - (2n)^{-1}\sum_{i=1}^{n}\bm{x}_ie^{z_i}(1-e^{z_i})(1+e^{z_i})^{-3}\big[\bm{x}_i^\prime(\hat{\bm{\beta}_n}-\bm{\beta})\big]^2,
\end{align}
where $|z_i-\bm{x}_i^\prime \bm{\beta}|\leq|\bm{x}_i^\prime (\hat{\bm{\beta}}_n-\bm{\beta})|$ for all $i\in \{1,\dots,n\}$. Now due to the assumption $n^{-1}\sum_{i=1}^{n}\|\bm{x}_i\|^3=O(1)$, by Lemma \ref{lem:1} (with $t=3$) we have
\begin{align}\label{eqn:1.2}
\mathbf{P}\Big(\big|n^{-1}\sum_{i=1}^{n}(y-p(\bm{\beta}|\bm{x}_i))x_{ij}\big|\leq C_{40}(p) n^{-1/2}(\log n)^{1/2}\Big)=o\big(n^{-1/2}\big),
\end{align}
for any $j\in \{1,\dots,p\}$. Again by assumption $\bm{L}_n$ converges to some positive definite matrix $\bm{L}$. Moreover,
\begin{align*}
\big\|(2n)^{-1}\sum_{i=1}^{n}\bm{x}_ie^{z_i}(1-e^{z_i})(1+e^{z_i})^{-3}\big[\bm{x}_i^\prime(\hat{\bm{\beta}}-\bm{\beta})\big]^2\big\| \leq \Big(n^{-1}\sum_{i=1}^{n}\|\bm{x}_i\|^3\Big)\|\hat{\bm{\beta}}_n-\bm{\beta}\|^2.  
\end{align*}
Hence (\ref{eqn:1.1}) can be rewritten as 
$$(\hat{\bm{\beta}}_n-\bm{\beta})= f_n(\hat{\bm{\beta}}_n-\bm{\beta}),$$ 
where $f_n$ is a continuous function from $\mathcal{R}^p$ to $\mathcal{R}^p$ satisfying $\mathbf{P}\Big(\|f_n\big(\hat{\bm{\beta}}_n-\bm{\beta}\big)\|\leq C_{40}n^{-1/2}(\log n)^{1/2}\Big) = 1 - o\big(n^{-1/2}\big)$ whenever $\|(\hat{\bm{\beta}}_n-\bm{\beta})\|\leq C_{40}n^{-1/2}(logn)^{1/2}$. Therefore, part (a) of Theorem \ref{thm:nor} follows by Brouwer's fixed point theorem. Now we are going to prove part (b). Note that from (\ref{eqn:1.1}) and the fact that $\bm{L}_n$ converges to some positive definite matrix $\bm{L}$, we have for large enough $n$,
\begin{align}\label{eqn:1.3}
\tilde{\mathbf{H}}_n = \hat{\bm{L}}_n^{1/2}\big[\bm{L}_n^{-1}\Lambda_n+ R_{1n}\big].
\end{align}
Here $\Lambda_n=n^{-1/2}\sum_{i=1}^{n}(y-p(\bm{\beta}|\bm{x}_i))\bm{x}_i$ and $R_{1n}= - \bm{L}_n^{-1}\dfrac{1}{2\sqrt{n}}\sum_{i=1}^{n}\bm{x}_ie^{z_i}(1-e^{z_i})(1+e^{z_i})^{-3}\big[\bm{x}_i^\prime(\hat{\bm{\beta}_n}-\bm{\beta})\big]^2$ with $|z_i-\bm{x}_i^\prime \bm{\beta}|\leq |\bm{x}_i^\prime (\hat{\bm{\beta}}_n-\bm{\beta})|$ for all $i\in \{1,\dots,n\}$. $\bm{L}_n$ and $\hat{\bm{L}}_n$ are as defined earlier. Now applying part (a) we have $\mathbf{P}\Big(\|R_{1n}\|=O\Big(n^{-1/2}(log n )\Big)\Big)=1- o\Big(n^{-1/2}\Big)$. Again by Taylor's theorem we have
\begin{align}\label{eqn:1.3}
\hat{\bm{L}}_n-\bm{L}_n = n^{-1}\sum_{i=1}^{n}\bm{x}_i\bm{x}_i^\prime e^{\bm{x}_i^\prime \bm{\beta}}\big(1-e^{\bm{x}_i^\prime \bm{\beta}}\big)\big(1+e^{\bm{x}_i^\prime \bm{\beta}}\big)^{-3}\big[\bm{x}_i^\prime(\hat{\bm{\beta}}_n-\bm{\beta})\big] + \bm{L}_{1n},
\end{align}
where by part (a), we have $\mathbf{P}\Big(\|\bm{L}_{1n}\|=O\Big(n^{-1}(log n )\Big)\Big)=1- o\Big(n^{-1/2}\Big)$.
Hence using Lemma \ref{lem:6}, part (a) and Taylor's theorem, one can show that $\mathbf{P}\Big(\|\hat{\bm{L}}_n^{1/2}-\bm{L}_n^{1/2}\|=O\Big(n^{-1/2}(log n )^{1/2}\Big)\Big)=1- o\Big(n^{-1/2}\Big)$. Therefore (\ref{eqn:1.2}) and (\ref{eqn:1.3}) will imply that $$\tilde{\mathbf{H}}_n=\bm{L}_n^{-1/2}\Lambda_n + R_{2n},$$ where $\mathbf{P}\Big(\|R_{2n}\|=O\Big(n^{-1/2}(log n )\Big)\Big)=1- o\Big(n^{-1/2}\Big)$. Hence for any set $B \in \mathcal{A}_p$, there exists a constant $C_{41}(p)>0$ such that
\begin{align*}
& \Big|\mathbf{P}\Big(\tilde{\mathbf{H}}_n \in B\Big)-\Phi(B)\Big|\nonumber\\
\leq \; &  \Big|\mathbf{P}\Big(\tilde{\mathbf{H}}_n \in B\Big)-\mathbf{P}\Big(\bm{L}_n^{-1/2}\Lambda_n \in B\Big)\Big| + \Big|\mathbf{P}\Big(\bm{L}_n^{-1/2}\Lambda_n \in B\Big)-\Phi(B)\Big|\nonumber \\
 \leq \; & \mathbf{P}\Big(\|R_{2n}\|> C_{41}(p) n^{-1/2}(log n)\Big) +  2\mathbf{P}\Big(\bm{L}_n^{-1/2}\Lambda_n \in (\partial B)^{C_{41}(p) n^{-1/2} (\log n)}\Big) + \Big|\mathbf{P}\Big(\bm{L}_n^{-1/2}\Lambda_n \in B\Big)-\Phi(B)\Big|\nonumber\\
 = \; & O\Big(n^{-1/2}(\log n)\Big).
\end{align*}
The last equality is a consequence of Lemma \ref{lem:5} and the bound on $\|R_{2n}\|$. Therefore part (b) is proved.\\

Proof of Theorem \ref{thm:posboot}. By applying Taylor's theorem, it follows from (\ref{eqn:defpb}) that
\begin{align}\label{eqn:3.1}
\hat{\bm{L}}_n\big(\hat{\bm{\beta}}_n^*-\hat{\bm{\beta}}_n\big)=n^{-1}\sum_{i=1}^{n}(y-\hat{p}(\bm{x}_i))\bm{x}_i - (2n)^{-1}\sum_{i=1}^{n}\bm{x}_ie^{z_i^*}(1-e^{z_i^*})(1+e^{z_i^*})^{-3}\big[\bm{x}_i^\prime(\hat{\bm{\beta}}_n^*-\hat{\bm{\beta}}_n)\big]^2,
\end{align}
where $|z_i^*-\bm{x}_i^\prime \bm{\beta}|\leq|\bm{x}_i^\prime (\hat{\bm{\beta}}_n^*-\hat{\bm{\beta}}_n)|$ for all $i\in \{1,\dots,n\}$. Now rest of part (a) of Theorem \ref{thm:posboot} follows exactly in the same line as the proof of part (a) of Theorem \ref{thm:nor}. To establish part (b), assume that $W_i = \Big(y_i\bm{x}^\prime_i,\big[y_i^2-\mathbf{E}y_i^2\big]\bm{z}^\prime_i\Big)^\prime$ and $W_i^*=\Big(\hat{Y}_i\big[(G_i^*-\mu_{G^*})\mu_{G^*}^{-1}\big]\bm{x}^\prime_i, \hat{Y}_i^2\big[\mu_{G^*}^{-2}(G_i^*-\mu_{G^*})^2-1\big]\bm{z}_i^\prime\Big)^\prime$. Here $y_i=(y_i-p(\bm{\beta}|\bm{x}_i))$ and $\hat{Y}_i=(y_i-\hat{p}(\bm{x}_i))$. First we are going to show that $$\check{\mathbf{H}}_n=  \sqrt{n}\Big(H\big(\bar{W}_n+n^{-1/2}b_nZ\big)\Big) + R_n\;\;\; \text{and}\;\;\; \check{\mathbf{H}}_n^*=  \sqrt{n}\Big(\hat{H}\big(\bar{W}_n^*+n^{-1/2}b_nZ\big)\Big) + R_n^*,$$ for some functions $H, \hat{H}:\mathcal{R}^k \rightarrow \mathcal{R}^p$ where $k=p+q$ with $q=\dfrac{p(p+1)}{2}$. $H(\cdot), \hat{H}(\cdot)$ have continuous partial derivatives of all orders, $H(\bm{0})=\hat{H}(\bm{0})=\bm{0}$ and $\mathbf{P}\Big(\|R_n\|=o\big(n^{-1/2}\big)\Big)=1-o\big(n^{-1/2}\big)$ \& $\mathbf{P}_*\Big(\|R_n^*\|=o\big(n^{-1/2}\big)\Big)=1-o_p\big(n^{-1/2}\big)$. Next step is to apply Lemma \ref{lem:3}, Lemma \ref{lem:4} and Lemma \ref{lem:7} to claim that suitable Edgeworth expansions exist for both $\check{\mathbf{H}}_n$ and $\check{\mathbf{H}}_n^*$. The last step is to conclude SOC of Bootstrap by comparing the Edgeworth expansions. Now (\ref{eqn:1.1}) and part (a) of Theorem \ref{thm:nor} imply that 
\begin{align}\label{eqn:3.2}
\sqrt{n}\big(\hat{\bm{\beta}}_n-\bm{\beta}\big)=\bm{L}_n^{-1}\Big[\Lambda_n-\xi_n/2\Big]+R_{3n},
\end{align}
where $\mathbf{P}\Big(\|R_{3n}\|\leq C_{42}(p)n^{-1} (\log n)^{3/2}\big)\Big)=1-o\big(n^{-1/2}\big)$. Here $\Lambda_n=n^{-1/2}\sum_{i=1}^{n}y_i\bm{x}_i$ and $\xi_n = n^{-3/2}\sum_{i=1}^{n}\bm{x}_ie^{\bm{x}_i^\prime \bm{\beta}}\big(1-e^{\bm{x}_i^\prime \bm{\beta}}\big)\big(1+e^{\bm{x}_i^\prime \bm{\beta}}\big)^{-3}\Big[\bm{x}_i^{\prime}\big(\bm{L}_n^{-1}\Lambda_n\big)\Big]^2$. Clearly, $\mathbf{P}\Big(\|\xi_n\|\leq C_{43}(p)n^{-1/2} (\log n)\big)\Big)=1-o\big(n^{-1/2}\big)$. Therefore, by Taylor's theorem we have 
\begin{align}\label{eqn:3.3}
\sqrt{n}\big(\hat{\bm{L}}_n-\bm{L}_n\big)\big(\hat{\bm{\beta}}_n-\bm{\beta}\big)= \xi_n + R_{4n},
\end{align}
where $\mathbf{P}\Big(\|R_{4n}\|\leq C_{44}(p)n^{-1} (\log n)^2\big)\Big)=1-o\big(n^{-1/2}\big)$. Again noting (\ref{eqn:3.3}), by equation (5) at page 52 of Turnbull (1929) we have 
\begin{align}\label{eqn:3.4}
\hat{\bm{M}}_n^{-1/2}-\bm{L}_n^{-1/2}=-\bm{L}_n^{-1/2}\bm{Z}_{1n}\bm{L}_n^{-1/2}+\bm{Z}_{2n},
\end{align}
where $\big(\hat{\bm{M}}_n-\bm{L}_n\big)=\bm{L}_n^{1/2}\bm{Z}_{1n}+\bm{Z}_{1n}\bm{L}_n^{1/2}$. Also easy to show that $$\mathbf{P}\Big(\|\hat{\bm{M}}_n-\bm{M}_n\|\leq C_{45}(p)n^{-1}(\log n)\Big)=1-o\big(n^{-1/2}\big),$$ where $\bm{M}_n=n^{-1}\sum_{i=1}^{n}y_i^2\bm{x}_i\bm{x}_i^\prime$. Hence using Lemma \ref{lem:6} we have  $\mathbf{P}\Big(\|\bm{Z}_{2n}\|\leq C_{46}(p)n^{-1} (\log n)^2\big)\Big)=1-o\big(n^{-1/2}\big)$. Therefore from (\ref{eqn:3.2})-(\ref{eqn:3.4}), Lemma \ref{lem:8} and the fact that $b_n=O(n^{-d})$ (for some $d>0$) will imply that
\begin{align}\label{eqn:3.5}
\check{\mathbf{H}}_n = \bm{L}_n^{-1/2}\Big[\Lambda_n+b_nZ+\xi_n/2\Big]    - \bm{L}_n^{-1/2}\Big[\int_{0}^{\infty}e^{-t\bm{L}_n^{1/2}}\big(\bm{M}_n-\bm{L}_n\big)e^{-t\bm{L}_n^{1/2}}dt\Big]\bm{L}_n^{-1/2}\Lambda_n +R_{5n},
\end{align}
where $\mathbf{P}\Big(\|R_{5n}\|\leq C_{47}(p)n^{-1/2} (\log n)^{-1}\Big)=1-o\big(n^{-1/2}\big)$. Now writing $W_i=(W_{i1}^\prime,W_{i2}^\prime)^\prime$ and $\bar{W}_n=n^{-1}\sum_{i=1}^{n}W_i=(\bar{W}_{n,1}^\prime,\bar{W}_{n2}^\prime)^\prime$ with $W_{i1}$ has first $p$ components of $W_i$ for all $i\in \{1,\dots,n\}$, we have
\begin{align*}
&\Lambda_n + b_n Z = \sqrt{n}\big(\bar{W}_{n1}+n^{-1/2}b_n Z\big)\\
& \xi_n = n^{-1/2}\sum_{i=1}^{n}\bm{x}_ie^{\bm{x}_i^\prime \bm{\beta}}\big(1-e^{\bm{x}_i^\prime \bm{\beta}}\big)\big(1+e^{\bm{x}_i^\prime \bm{\beta}}\big)^{-3}\Big[\bar{W}_{n1}^\prime \bm{L}_n^{-1}\bm{x}_i \bm{x}_i^\prime \bm{L}_n^{-1} \bar{W}_{n1}\Big]^2\\
&\;\;\; = \sqrt{n}\Big(\bar{W}_{n1}^\prime \tilde{\bm{M}}_1 \bar{W}_{n1},\dots, \bar{W}_{n1}^\prime \tilde{\bm{M}}_p \bar{W}_{n1}\Big)^\prime,
\end{align*}
where $\tilde{\bm{M}}_k = n^{-1}\sum_{i=1}^{n}x_{ik}e^{\bm{x}_i^\prime \bm{\beta}}\big(1-e^{\bm{x}_i^\prime \bm{\beta}}\big)\big(1+e^{\bm{x}_i^\prime \bm{\beta}}\big)^{-3}\Big(\bm{L}_n^{-1}\bm{x}_i \bm{x}_i^\prime \bm{L}_n^{-1}\Big)$ for $k\in \{1,\dots,p\}$. Hence writing $\tilde{W}_{n1}=\bar{W}_{n1}+n^{-1/2}b_n Z$ we have
\begin{align}\label{eqn:3.6}
\bm{L}_n^{-1/2}\Big[\Lambda_n+b_nZ+\xi_n/2\Big]  = \sqrt{n}\bigg[\bm{L}_n^{-1/2}\tilde{W}_n + \Big(\tilde{W}_{n1}^\prime \breve{\bm{M}}_1 \tilde{W}_{n1},\dots, \tilde{W}_{n1}^\prime \breve{\bm{M}}_p \tilde{W}_{n1}\Big)^\prime\bigg],
\end{align}
since $b_n = O(n^{-d})$ and $\|\tilde{\bm{M}}_k\|=O(1)$ for any $k\in \{1,\dots,p\}$. Here $\breve{\bm{M}}_k=\sum_{j=1}^{p}L_{kjn}^{-1/2}\tilde{\bm{M}}_k$, $k \in \{1,\dots,p\}$, with $L_{kjn}^{-1/2}$ being the $(k,j)$th element of $\bm{L}_n^{-1/2}$. Again the $j$th row of $\big(\bm{M}_n-\bm{L}_n\big)$ is $\bar{W}^\prime_{n2}\bm{E}_{jn}$ where $\bm{E}_{jn}$ is a matrix of order $q\times p$ with $\|\bm{E}_{jn}\|\leq q$, $j\in\{1,\dots,p\}$. Therefore from (\ref{eqn:1.1}) and (\ref{eqn:1.3}) we have
\begin{align}\label{eqn:3.7}
\bm{L}_n^{-1/2}\Big[\int_{0}^{\infty}e^{-t\bm{L}_n^{1/2}}\big(\bm{M}_n-\bm{L}_n\big)e^{-t\bm{L}_n^{1/2}}dt\Big]\bm{L}_n^{-1/2}\Lambda_n =  \sqrt{n} \Big(\tilde{W}_{n2}^\prime \check{\bm{M}}_1 \tilde{W}_{n2}, \dots, \tilde{W}_{n2}^\prime \check{\bm{M}}_p \tilde{W}_{n2}\Big)^\prime +R_{6n}, 
\end{align}where $\tilde{W}_{n2}=\bar{W}_{n2}+n^{-1/2}b_nZ_1$ with $Z_1\sim N_{q}\big(\bm{0},\bm{I}_q\big)$, independent of $Z$ \& $\{y_1,\dots,y_n\}$. $\bar{\bm{M}}_k=\bigintss_{0}^{\infty}\Big[\sum_{j=1}^{p}m_{kjn}(t)\check{\bm{M}}_j(t)\big]dt$ where $m_{kjn}(t)$ is the $(k,j)$th element of the matrix $\bm{L}_n^{-1/2}e^{t\bm{L}_n^{1/2}}$ and $\check{\bm{M}}_j(t)=\bm{E}_{jn}e^{t\bm{L}_n^{1/2}}\bm{L}_n^{-1/2}$, $k,j\in \{1,\dots,p\}$.  Moreover, $\mathbf{P}\Big(\|R_{6n}\|\leq C_{48}n^{-1/2}(\log n)^{-1}\big)\Big)=1-o\big(n^{-1/2}\big)$. Now define the $(p+q)\times (p+q)$ matrices $\{\bm{M}_1^\dagger, \dots, \bm{M}_p^\dagger\}$ where $\bm{M}_k^\dagger=\begin{bmatrix}\bm{\breve{M}}_k & \bm{0} \\ \bar{\bm{M}}_k  & \bm{0}\end{bmatrix}$.
 Therefore from (\ref{eqn:3.5})-(\ref{eqn:3.7}) we have
\begin{align}\label{eqn:3.8}
\check{\mathbf{H}}_n & = \sqrt{n}\bigg[\Big(\bm{L}_n^{-1/2}\;\; \bm{0}\Big)\tilde{W}_n + \Big(\tilde{W}_n^\prime \bm{M}_1^\dagger \tilde{W}_n,\dots, \tilde{W}_n^\prime \bm{M}_p^\dagger \tilde{W}_n\Big)^\prime\bigg]+ R_n\nonumber\\
& = \sqrt{n}H\big(\tilde{W}_n\big) + R_{n},
\end{align}
where the function $H(\cdot)$ has continuous partial derivatives of all orders, $\tilde{W}_n=\big(\tilde{W}_{n1}^\prime,\tilde{W}_{n2}^\prime\big)^\prime$ and $R_n =R_{5n}+R_{6n}$.\\

Through the same line of arguments, writing $\bar{W}_{n1}^*=n^{-1}\sum_{i=1}^{n}W_{i1}^*=n^{-1}\sum_{i=1}^{n}\hat{Y}_i\mu_{G^*}^{-1}(G_i^*-\mu_{G^*})\bm{x}_i$ and $\bar{W}_{n2}^*=n^{-1}\sum_{i=1}^{n}W_{i2}^*=n^{-1}\sum_{i=1}^{n}\hat{Y}_i^2\big[\mu_{G^*}^{-2}(G_i^*-\mu_{G^*})^2-1\big]\bm{z}_i$, it can be shown that
\begin{align}\label{eqn:3.9}
\check{\mathbf{H}}_n^* & = \sqrt{n}\bigg[\Big(\hat{\bm{M}}_n^{-1/2}\;\;\bm{0}\Big)\tilde{W}_n^* + \Big(\tilde{W}_n^{*\prime} \bm{M}^{*\dagger}_1 \tilde{W}_n^*,\dots, \tilde{W}_n^{*\prime}\bm{M}^{*\dagger}_p \tilde{W}_n^*\Big)^\prime\bigg]+ R_n^*\nonumber\\
& = \sqrt{n}\hat{H}\big(\tilde{W}_n^*\big) + R_{n}^*,
\end{align}
where $\tilde{W}_n^*=\big(\tilde{W}_{n1}^{*\prime},\tilde{W}_{n2}^{*\prime}\big)^\prime$ with $\tilde{W}_{n1}^*=\bar{W}_{n1}^*+n^{-1/2}b_n Z^*$ and $\tilde{W}_{n2}^*=\bar{W}_n^*+n^{-1/2}b_nZ_1^*$, $Z_1^*$ being a $N_q\big(\bm{0},\bm{I}_q\big)$ distributed random vector independent of $\{G_1^*,\dots,G_n^*\}$ and $Z^*$. $\bm{M}^{*\dagger}_j=\begin{bmatrix}\bm{\breve{M}}^*_k & \bm{0} \\ \bar{\bm{M}}^*_k  & \bm{0}\end{bmatrix}$ where $\breve{\bm{M}}_j^*=\sum_{j=1}^{p}\hat{M}_{kjn}^{-1/2}\tilde{\bm{M}}_j^*$ with $\hat{M}_{kjn}^{-1/2}$ being the $(k,j)$th element of $\hat{\bm{M}}_n^{-1/2}$, $\tilde{\bm{M}}_j^*$ being same as $\tilde{\bm{M}}_j$ after replacing $\bm{\beta}$ by $\hat{\bm{\beta}}_n$. $\bar{\bm{M}}_j^*=\bigintss_{0}^{\infty}\big[\sum_{j=1}^{p}m_{kjn}^*\check{\bm{M}}_j(t)\big]dt$ where $m^*_{kjn}(t)$ is the $(k,j)$th element of the matrix $\hat{\bm{M}}_n^{-1/2}e^{-t\hat{\bm{M}}_n^{1/2}}$ and $\check{\bm{M}}^*_j(t)=\bm{E}_{jn}e^{-t\hat{\bm{M}}_n^{1/2}}\hat{\bm{M}}_n^{-1/2}$.  Also $\mathbf{P}_*\big(\|R_n^*\|\leq C_{49}n^{-1/2}(\log n)^{-1}\big)=1-o_p\big(n^{-1/2}\big)$. Now by applying Lemma \ref{lem:3}, Lemma \ref{lem:4} and Lemma \ref{lem:7} with $s=3$, Edgeworth expansions of the densities of $\check{\mathbf{H}}_n$ and $\check{\mathbf{H}}_n^*$ can be found uniformly over the class $\mathcal{A}_p$ upto an error $o\big(n^{-1/2}\big)$ and $o_p\big(n^{-1/2}\big)$ respectively. Call those Edgeworth expansions $\tilde{\psi}_{n,3}(\cdot)$ and $\tilde{\psi}_{n,3}^*(\cdot)$ respectively. Now if $\tilde{\psi}_{n,3}(\cdot)$ is compared with $\check{\psi}_{n,3}(\cdot)$ of Lemma \ref{lem:4}, then $\check{\bm{M}}_n= I_p$. Similarly for $\tilde{\psi}_{n,3}^*(\cdot)$ also $\check{\bm{M}}_n=I_p$. Therefore, $\tilde{\psi}_{n,3}(\cdot)$ and $\tilde{\psi}_{n,3}^*(\cdot)$ have the forms
\begin{align*}
&\tilde{\psi}_{n,3}(\bm{x}) = \Big[1 + n^{-1/2}q_1(\bm{\beta}, \mu_{W},\bm{x}) + \sum_{j=1}^{m_2-1}b_n^{2j}q_{2j}(\bm{\beta},\bm{L}_n,\bm{x})\Big]\phi(\bm{x})\\
&\tilde{\psi}_{n,3}^*(\bm{x}) = \Big[1 + n^{-1/2}q_1(\hat{\bm{\beta}}_n, \hat{\mu}_{W}, \bm{x}) + \sum_{j=1}^{m_2-1}b_n^{2j}q_{2j}(\hat{\bm{\beta}}_n,\hat{\bm{M}}_n,\bm{x})\Big]\phi(\bm{x}),
\end{align*}
where $m_2=\inf\{j:b_n^{2j}=o(n^{-1/2})\}$, $\mu_{W}$ is the vector of $\{n^{-1}\sum_{i=}^{n}\mathbf{E}(y_i-p(\bm{\beta}|\bm{x}_i))^2x_{ij}^{l_1}x_{ij^\prime}^{l_2}:j,j^\prime\in\{1,\dots,p\}, l_1,l_2\in \{0,1,2\}, l_1+l_2=2\}$ and $\{n^{-1}\sum_{i=}^{n}\mathbf{E}(y_i-p(\bm{\beta}|\bm{x}_i))^3x_{ij}^{l_1}x_{ij^\prime}^{l_2}x_{ij^{\prime\prime}}^{l_3}:j,j^\prime, j^{\prime\prime}\in\{1,\dots,p\}, l_1,l_2,l_3\in \{0,1,2,3\},  l_1+l_2+l_3=3\}$. $\hat{\mu}_{W}$ is the vector of $\{n^{-1}\sum_{i=}^{n}(y_i-\hat{p}(\bm{x}_i))^2x_{ij}^{l_1}x_{ij^\prime}^{l_2}:j,j^\prime\in\{1,\dots,p\}, l_1,l_2\in \{0,1,2\}, l_1+l_2=2\}$ and $\{n^{-1}\sum_{i=}^{n}(y_i-\hat{p}(\bm{x}_i))^3x_{ij}^{l_1}x_{ij^\prime}^{l_2}x_{ij^{\prime\prime}}^{l_3}:j,j^\prime,j^{\prime\prime}\in\{1,\dots,p\}, l_1,l_2\in \{0,1,2,3\},  l_1+l_2+l_3=3\}$.  $q_1(\bm{a},\bm{b},\bm{c})$ is a polynomial in $\bm{c}$ whose coefficients are continuous functions of $(\bm{a},\bm{b})^\prime$. $q_{2j}(\bm{a},\bm{b},\bm{c})$ are polynomials in $\bm{c}$ whose coefficients are continuous functions of $\bm{a}$ and $\bm{b}$.
Now Theorem \ref{thm:posboot} follows by comparing $\tilde{\psi}_{n,3}(\cdot)$ and $\tilde{\psi}_{n,3}^*(\cdot)$ and due to part (a) of Theorem \ref{thm:nor}.\\

Proof of Theorem \ref{thm:negboot}. Recall that here $p=1$ and hence $q=1$. 
Define, $\bm{B}_n=\sqrt{n}H(\bm{E}_n\times \mathcal{R})$ with $\bm{E}_n=(-\infty,z_n]$ and $z_n=\Big(\dfrac{3}{4n}-\mu_n\Big)$. Here $\mu_n=n^{-1}\sum_{i=1}^{n}x_ip(\beta|x_i)$. Note that $B_n$ is an interval, as argued in section \ref{sec:main} just after the description of Theorem \ref{thm:negboot}. The function $H(\cdot)$ is defined in (\ref{eqn:3.8}). We are going to show that there exists a positive constant $M_2$ such that $$\lim_{n\rightarrow \infty}\mathbf{P}\Big(\sqrt{n}\Big|\mathbf{P}_{*}\big(\mathbf{H}_n^* \in \bm{B}_n\big)-\mathbf{P}\big(\mathbf{H}_n \in \bm{B}_n\big)\Big|\geq M_3\Big)=1.$$
Define the set $\bm{Q}_n=\Big\{\big|\hat{\bm{\beta}}_n-\bm{\beta}\big| = o\big(n^{-1/2} (\log n)\big)\Big\} \cap \Big\{\big|n^{-1}\sum_{i=1}^{n}\big[(y_i-p(\beta|x_i))^2-\mathbf{E}(y_i-p(\beta|x_i))^2\big]x_i^2\big| =o\big(n^{-1/2}(log n)\big)\Big\} \cap \Big\{\big|n^{-1}\sum_{i=1}^{n}\big[(y_i-p(\beta|x_i))^3-(y_i-p(\beta|x_i))^3\big]x_i^3\big| =o(1)\Big\}$. Now due to a stronger version of (\ref{eqn:1.2}), it is easy to see that $\mathbf{P}\Big(\big|\hat{\bm{\beta}}_n-\bm{\beta}\big| = o\big(n^{-1/2} (\log n)\Big)=1$ for all but finitely many $n$, upon application of Borel-Cantelli lemma and noting that $\max\{|x_i|:i\in \{1,\dots,n\}\}=O(1)$. Again by applying Lemma \ref{lem:1}, it is easy to show that $\mathbf{P}\Big(\Big\{\big|n^{-1}\sum_{i=1}^{n}\big[(y_i-p(\beta|x_i))^2-(y_i-p(\beta|x_i))^2\big]\big| =o\big(n^{-1/2}(log n)\big)\Big\} \cap \Big\{\big|n^{-1}\sum_{i=1}^{n}\big[(y_i-p(\beta|x_i))^3-(y_i-p(\beta|x_i))^3\big]\big| =o(1)\Big\}\Big)=1$ for large enough $n$. Hence $\mathbf{P}\big(\bm{Q}_n\big)=1$ for large enough $n$. Similarly define the Bootstrap version of $\bm{Q}_n$ as $\bm{Q}_n^*=\Big\{\big|\hat{\bm{\beta}}_n^*-\hat{\bm{\beta}}_n\big|= o\big(n^{-1/2} (\log n)\big)\Big\} \cap \Big\{\big|n^{-1}\sum_{i=1}^{n}\big[(y_i-\hat{p}(x_i))^2\big(\mu_{G^*}^{-2}(G_i^*-\mu_{G^*})^2-1\big)\big]x_i^2\big| =o\big(n^{-1/2}(log n)\big)\Big\} \cap \Big\{\big|n^{-1}\sum_{i=1}^{n}\big[(y_i-\hat{p}(x_i))^3\big(\mu_{G^*}^{-3}(G_i^*-\mu_{G^*})^3-1\big)\big]x_i^3\big| =o\big(1)\big)\Big\}$. Through the same line, it is easy to establish that $\mathbf{P}\Big(\mathbf{P}^*\big(\bm{Q}^*_n\big)=1\Big)=1$ for large enough $n$.
Hence enough to show 
\begin{align}\label{eqn:2.1}
\lim_{n\rightarrow \infty}\mathbf{P}\Big(\Big\{\sqrt{n}\Big|\mathbf{P}_{*}\big(\big\{\mathbf{H}_n^* \in \bm{B}_n\}\cap \bm{Q}_n^*\big)-\mathbf{P}\big(\big\{\mathbf{H}_n \in \bm{B}_n\big\}\cap \bm{Q}_n\big)\Big|\geq M_2\Big\}\cap \bm{Q}_n\Big)=1.
\end{align}
Recall the definitions of $\bar{W}_n$ and $\bar{W}_n^*$ from the proof of Theorem \ref{thm:posboot}. Similar to (\ref{eqn:3.8}) and (\ref{eqn:3.9}), it is easy to observe that
\begin{align}\label{eqn:2.2}
\mathbf{H}_n=\sqrt{n}H(\bar{W}_n)+R_{n}\;\; \text{and}\;\; \mathbf{H}^*_n=\sqrt{n}\hat{H}(\bar{W}^*_n)+R^*_{n},
\end{align}
where $\big\{|R_{n}|=O(n^{-1/2}(\log n)^{-1})\big\}\subseteq \bm{Q}_n$ and $\big\{|R_{n}^*|=O(n^{-1/2}(\log n)^{-1})\big\}\subseteq \bm{Q}_n^*$. To prove (\ref{eqn:2.1}), first we are going to show for large enough $n$, 
\begin{align}\label{eqn:2.3}
&\bigg\{\Big\{\sqrt{n}\Big|\mathbf{P}_{*}\big(\big\{\mathbf{H}_n^* \in \bm{B}_n\}\cap \bm{Q}_n^*\big)-\mathbf{P}\big(\big\{\mathbf{H}_n \in \bm{B}_n\big\}\cap \bm{Q}_n\big)\Big|\geq M_2\Big\}\cap \bm{Q}_n\bigg\}\nonumber\\
&\supseteq \bigg\{\Big\{\sqrt{n}\Big|\mathbf{P}_{*}\big(\big\{\sqrt{n}\hat{H}(\bar{W}_n^*) \in \bm{B}_n\}\cap \bm{Q}_n^*\big)-\mathbf{P}\big(\big\{\sqrt{n}H(\bar{W}_n) \in \bm{B}_n\big\}\cap \bm{Q}_n\big)\Big|\geq 2M_2\Big\}\cap \bm{Q}_n\bigg\}.
\end{align}
Now due to (\ref{eqn:2.2}), we have 
\begin{align*}
\Big|\mathbf{P}\Big(\mathbf{H}_n \in \bm{B}_n\Big) - \mathbf{P}\Big(\sqrt{n}H(\bar{W}_n) \in \bm{B}_n\Big)\Big| \leq\; &  \mathbf{P}\Big(\sqrt{n}H(\bar{W}_n) \in \big(\partial\bm{B}_n\big)^{(n \log n)^{-1/2}}\Big)\nonumber\\ 
& + \mathbf{P}\Big(|R_{n}|\neq o(n^{-1/2}(\log n)^{-1})\Big)\nonumber\\
 \Big|\mathbf{P}_{*}\Big(\mathbf{H}_n^* \in \bm{B}_n\Big) - \mathbf{P}_{*}\Big(\sqrt{n}\hat{H}(\bar{W}^*_n) \in \bm{B}_n\Big)\Big| \leq\; & \mathbf{P}_{*}\Big(\sqrt{n}\hat{H}(\bar{W}^*_n) \in \big(\partial\bm{B}_n\big)^{(n \log n)^{-1/2}}\Big)\nonumber\\
 &+ \mathbf{P}_*\Big(|R_{n}^*|\neq o(n^{-1/2}(\log n)^{-1})\Big)
\end{align*}
To establish (\ref{eqn:2.3}), enough to show $\mathbf{P}\Big(\sqrt{n}\hat{H}(\bar{W}_n) \in \big(\partial\bm{B}_n\big)^{(n \log n)^{-1/2}}\Big)=o\big(n^{-1/2}\big)$ and $\mathbf{P}\Big(\Big\{\mathbf{P}_{*}$ $\Big(\sqrt{n}\hat{H}(\bar{W}^*_n) \in \big(\partial\bm{B}_n\big)^{(n \log n)^{-1/2}}\Big)=o\big(n^{-1/2}\big)\Big\}\cap \bm{Q}_n\Big)=1$ for large enough $n$. An Edgeworth expansion of $\sqrt{n}\bar{W}_n^*$ with an error $o(n^{-1/2})$ (in almost sure sense) can be established using Lemma \ref{lem:11}. Then we can use transformation technique of Bhattacharya and Ghosh (1978) to find an Edgeworth expansion $\hat{\eta}_n(\cdot)$ of the density of $\sqrt{n}\hat{H}(\bar{W}_n^*)$ with an error $o(n^{-1/2})$ (in almost sure sense). Now the calculations similar to page 213 of BR(86) will imply that 
$\mathbf{P}\Big(\Big\{\mathbf{P}_{*}\Big(\sqrt{n}\hat{H}(\bar{W}^*_n) \in \big(\partial\bm{B}_n\big)^{(n \log n)^{-1/2}}\Big)=o\big(n^{-1/2}\big)\Big\}\cap \bm{Q}_n\Big)=1$, since $\bm{B}_n$ is an interval. Next we are going to show that $\mathbf{P}\Big(\sqrt{n}\hat{H}(\bar{W}_n) \in \big(\partial\bm{B}_n\big)^{(n \log n)^{-1/2}}\Big)=0$ for large enough $n$ and to show that we need to utilize the form of $\bm{B}_n$, as Edgeworth expansion of $\sqrt{n}H(\bar{W}_n)$ similar to $\sqrt{n}\hat{H}(\bar{W}_n^*)$ does not exist due to the lattice nature of $W_1,\dots,W_n$. To this end define $k_n(\bm{x})=\big(\sqrt{n}H(\bm{x}/\sqrt{n}),x_2\big)^\prime$ where $\bm{x}=(x_1,x_2)^\prime$. Note that $k_n(\cdot)$ is a diffeomorphism (cf. proof of lemma 3.2 in Lahiri (1989)). Hence $k_n(\cdot)$ is a bijection and $k_n(\cdot)$ \& $k_n^{-1}(\cdot)$ have derivatives of all orders. Therefore, arguments given between (2.15) and (2.18) at page 444 of Bhattacharya and Ghosh (1978) with $g_n$ replaced by $k_n^{-1}(\cdot)$ will imply that 
\begin{align*}
\Big|\mathbf{P}\Big(\mathbf{H}_n \in \bm{B}_n\Big) - \mathbf{P}\Big(\sqrt{n}H(\bar{W}_n) \in \bm{B}_n\Big)\Big| &\leq \mathbf{P}\Big(\big(\sqrt{n}\bar{W}_n \in \big(\partial k_n^{-1}(\bm{B}_n\times\mathcal{R})\big)^{d_n(n \log n)^{-1/2}}\Big) +o\big(n^{-1/2}\big)\\
 & = \mathbf{P}\Big(\sqrt{n}\bar{W}_{n1} \in \big(\partial \bm{E}_n\big)^{d_n(n \log n)^{-1/2}}\Big) +o\big(n^{-1/2}\big), 
\end{align*}
where $d_n \leq \max\big\{|det\big(Grad\big[k_n(x)\big]\big)|^{-1}:|x|=O(\sqrt{\log n})\big\}$. Now by looking into the form of $H(\cdot)$ in (\ref{eqn:1.2}), it is easy to see that $d_n=O(1)$, say $d_n\leq C_{44}$ for some positive constant $C_{44}$. Now note that 
\begin{align*}
&\mathbf{P}\Big(\sqrt{n}\bar{W}_{n1} \in \big(\partial \bm{E}_n\big)^{C_{44}(n \log n)^{-1/2}}\Big)\\
& = \mathbf{P}\Big(\Big[n^{-1/2}\sum_{i=1}^{n}y_ix_i - \sqrt{n}\mu_n\Big] \in \big(z_n-C_{44}(n \log n)^{-1/2},z_n+C_{44} (n \log n)^{-1/2}\big)\Big)\\
& = \mathbf{P}\Big(\sum_{i=1}^{n}y_ix_i \in \big(3/4-C_{44}(\log n)^{-1/2},3/4+C_{44} (\log n)^{-1/2}\big)\Big)\\
& = 0,
\end{align*}
for large enough $n$, since $\sum_{i=1}^{n}y_ix_i$ can take only integer values. Therefore (\ref{eqn:2.3}) is established.
Now recalling that $\hat{\eta}_n(\cdot)$ is the Egeworth expansion of the density of $\sqrt{n}\hat{H}(\bar{W}_n^*)$ with an almost sure error $o(n^{-1/2})$, we have for large enough $n$,
\begin{align}\label{eqn:2.4}
\mathbf{P}\Big(\sqrt{n}\Big|\mathbf{P}_*\Big(\sqrt{n}\hat{H}(\bar{W}_n^*) \in \bm{B}_n\Big)-\int_{\bm{B}_n}\hat{\eta}_n(x)dx\Big|=o(1)\Big)=1.    
\end{align}
Now define $U_i=\Big(\big(y_i-p(\beta|x_i)\big)x_iV_i,\big(y_i-p(\beta|x_i)\big)^2x_i^2\big[V_i^2-1\big]\Big)^\prime$, $i\in\{1,\dots,n\}$, where $V_1,\dots,V_n$ are iid continuous random variables which are independent of $\{y_1,\dots,y_n\}$. Also $\mathbf{E}(V_1)=0$, $\mathbf{E}(V_1^2)=\mathbf{E}(V_1)^3=1$ and $\mathbf{E}V_1^8< \infty$. An immediate choice of the distribution of $V_1$ is that of $(G_1^*-\mu_{G^*})\mu_{G^*}^{-1}$. Other choices of $\{V_1,\dots,V_n\}$ can be found in Liu(1988), Mammen (1993) and Das et al. (2019). Now since $\max\{|x_i|:i\in \{1,\dots,n\}\}=O(1)$, there exists a natural number $n_0$ and constants $0<\delta_2\leq\delta_1<1$ such that $\sup_{n\geq n_0}p(\beta|x_n)\leq\delta_1$ and $\inf_{n\geq n_0}p(\beta|x_n)\geq \delta_2$. Again $V_1,\dots,V_n$ are iid continuous random variables. Hence writing $p_n=p(\beta|x_n)$, for any $b>0$ we have 
\begin{align*}
& \sup_{n\geq n_0}\sup_{\|\bm{t}\|>b}\Big|\mathbf{E}e^{i\bm{t}^\prime U_n}\Big|\\
& \leq \sup_{n\geq n_0}\bigg[p_n \sup_{\|\bm{t}\|>b(1-\delta_1)^2}\Big|\mathbf{E}e^{it_1(1-p_n)V_1+it_2(-p_n)^2[V_1^2-1]}\Big|+(1-p_n) \sup_{\|\bm{t}\|>b\delta_2^2}\Big|\mathbf{E}e^{it_1(-p_n)V_1+t_2(-p_n)^2[V_1^2-1]}\Big|\bigg]\\
& <1,
\end{align*}
i.e. uniform Cramer's condition holds. Also the minimum eigen value condition of Theorem 20.6 of BR(86) holds due to $\max\{|x_i|:i\in \{1,\dots,n\}\}=O(1)$ and $\liminf_{n\rightarrow \infty}n^{-1}\sum_{i=1}^{n}x_i^6>0$. Hence applying Theorem 20.6 of BR(86) and then applying transformation technique of Bhattacharya and Ghosh (1978) we have
\begin{align}\label{eqn:2.5}
\Big|\mathbf{P}\Big(\sqrt{n}H(\bar{U}_n) \in \bm{B}_n\Big)-\int_{\bm{B}_n}\eta_n(x)dx\Big|=o\big(n^{-1/2}\big),
\end{align}
where $\bar{U}_n=n^{-1}\sum_{i=1}^{n}U_i$. Note that in both the expansions $\eta_n(\cdot)$ and $\hat{\eta}_n(\cdot)$ the variances corresponding to normal terms are $1$. Also $\hat{H}(\cdot)$ can be obtained from $H(\cdot)$ first replacing $\bm{L}_n$ by $\hat{M}_n$ and then $\beta$ by $\hat{\beta}_n$ (cf. (\ref{eqn:3.8}) and (\ref{eqn:3.9})). Hence we can conclude that for any Borel set $\bm{C}$,
$$\mathbf{P}\Big(\Big\{\sqrt{n}\Big|\int_{\bm{C}}\eta_n(x)dx-\int_{\bm{C}}\hat{\eta}_n(x)dx\Big|=o\big(1\big)\Big\}\cap \bm{Q}_n\Big)=1$$ Hence from (\ref{eqn:2.4}) and (\ref{eqn:2.5}), we have
\begin{align}\label{eqn:2.6}
\mathbf{P}\bigg(\Big\{\sqrt{n}\Big|\mathbf{P}_*\Big(\sqrt{n}\hat{H}(\bar{W}_n^*) \in \bm{B}_n\Big)-\mathbf{P}\Big(\sqrt{n}H(\bar{U}_n) \in \bm{B}_n\Big)\Big|=o(1)\Big\}\cap \bm{Q}_n\bigg)=1,  
\end{align}
for large enough $n$. To establish (\ref{eqn:2.1}), in view of (\ref{eqn:2.3}) and (\ref{eqn:2.6}) it is enough to find a positive constant $M_3$ such that 
\begin{align*}
\sqrt{n}\Big|\mathbf{P}\Big(\sqrt{n}H(\bar{W}_n) \in \bm{B}_n\Big)-\mathbf{P}\Big(\sqrt{n}H(\bar{U}_{n}) \in \bm{B}_n\Big)\Big|=  \sqrt{n}\Big|\mathbf{P}\Big(\sqrt{n}\bar{W}_{n1} \in \bm{E}_n\Big)-\mathbf{P}\Big(\sqrt{n}\bar{U}_{n1} \in \bm{E}_n\Big)\Big|\geq 4M_3.
\end{align*}
Note that since $\mathbf{E}V_i^2=\mathbf{E}V_i^3=1$ for all $i\in \{1,\dots,n\}$, the first three average moments of $\{W_{11},\dots,W_{n1}\}$ are same as that of $\{U_{11},\dots,U_{n1}\}$. However $\{W_{11},\dots,W_{n1}\}$ are independent lattice random variables whereas $\{U_{11},\dots,U_{n1}\}$ are independent random variables for which uniform Cramer's condition holds. Therefore by Lemma \ref{lem:10} and Theorem 20.6 of BR(86) we have
\begin{align}\label{eqn:2.7}
&\sup_{x\in \mathcal{R}}\Big|\mathbf{P}\Big(\sqrt{n}\bar{W}_{n1}\leq x\Big)-\Phi_{\sigma_n^2}(x)-n^{-1/2}P_1\big(-\Phi_{\sigma_n^2}:\{\bar{\chi}_{\nu,n}\}\big)(x)\nonumber\\
&\;\;\;\;\;\;\;\;+ n^{-1/2}\big(n\mu_n+\sqrt{n}x-1/2)\dfrac{d}{dx}\Phi_{\sigma_n^2}(x)\Big| = o\big(n^{-1/2}\big)\nonumber\\
\text{and}\;\;& \sup_{x\in \mathcal{R}}\Big|\mathbf{P}\Big(\sqrt{n}\bar{U}_{n1}\leq x\Big)- \Phi_{\sigma_n^2}(x)
-n^{-1/2}P_1\big(-\Phi_{\sigma_n^2}:\{\bar{\chi}_{\nu,n}\}\big)(x)\Big| = o\big(n^{-1/2}\big),
\end{align}
where $P_1\big(-\Phi_{\sigma_n^2}:\{\bar{\chi}_{\nu,n}\}\big)(x)$ is as defined in Lemma \ref{lem:10}. Recall that $\bm{E}_n=(-\infty,z_n]$ where $z_n=\Big(\dfrac{3}{4n}-\mu_n\Big)$. Therefore for some positive constants $C_{46}, C_{47}, C_{48}$ we have
\begin{align*}
\sqrt{n}\Big|\mathbf{P}\Big(\sqrt{n}\bar{W}_{n1} \in \bm{E}_n\Big)-\mathbf{P}\Big(\sqrt{n}\bar{U}_{n1} \in \bm{E}_n\Big)\Big|
=\;&\sqrt{n}\Big|\mathbf{P}\Big(\sqrt{n}\bar{W}_{n1} \leq \sqrt{n}z_n\Big)-\mathbf{P}\Big(\sqrt{n}\bar{U}_{n1} \leq \sqrt{n}z_n\Big)\Big|\\
\geq \; & \big(n\mu_n+nz_n-1/2\big)\big(\sqrt{2\pi}\sigma_n\big)^{-1}e^{-(nz_n^2)/(2\sigma_n^2)} - o(1)\\
=\; & \big(4\sqrt{2\pi}\sigma_n\big)^{-1}e^{-(nz_n^2)/(2\sigma_n^2)}-o(1)\\
\geq \; & C_{46} \exp\Big\{-C_{47}n^{-1}\Big(\dfrac{9}{16}+n^2\mu_n^2-\dfrac{3n\mu_n}{2}\Big)\Big\} - o(1)\\
\geq \;& C_{48} \exp\Big\{-C_{47}M_1^2\Big\}. 
\end{align*}
The first inequality follows due to (\ref{eqn:2.7}). Second one is due to $\max\{|x_i|:i\in \{1,\dots,n\}\}=O(1)$ and the last one is due to the assumption $\sqrt{n}|\mu_n|<M_1$. Taking $4M_2=C_{48} \exp\Big\{-C_{47}M_1^2\Big\}$, the proof of Theorem \ref{thm:negboot} is now complete.

\section{Conclusion} \label{sec_con}  In this paper we consider the studentized version of the logistic regression estimator and proposed a novel Bootstrap method called PEBBLE. The rate of convergence of the studentized version to normal distribution is found to be sub-optimal with respect to the classical Berry-Esseen rate $O\big(n^{-1/2}\big)$. %A novel perturbation Bootstrap technique is developed for approximating the distribution of the logistic regression estimator to improve the rate of convergence. 
We observe that the usual studentization also fails  significantly in improving the error rate in the Bootstrap approximation due to the underlying lattice structure. Therefore, a novel modification is proposed in the form of studentized pivots to achieve SOC by the Bootstrap in approximating the distribution of the studentized logistic regression estimator. The proposed Bootstrap method can be used in practical purposes to draw inferences about the regression parameter which will be more accurate than that based on asymptotic normality. PEBBLE is shown perform better than other existing method performance-wise, in general, via simulation experiments. Specifically for larger $p$, smaller $n$ settings, PEBBLE outperforms other methods by a large margin. The proposed method is used to find the middle, upper and lower CIs for the covariates in a real data application concerning the dependency of the type of delivery on several related clinical variables. As a future extension, the SOC of Bootstrap in the generalized linear model (GLM) can be explored. %\rd{(Probably not required: One can expect to obtain interesting results under different sub-models.   For example when the underlying setup is linear it is already known that the usual Bootstrap achieves SOC, whereas in logistic regression it requires smoothing to achieve SOC.)}
Additionally, one can also explore the high dimensional structure in GLM, that is when dimension $p$ grows with $n$, by adding suitable penalty terms in the underlying objective function.

\section{Supplementary Material}

\subsection{Supplementary Proof Details}

%\setcounter{section}{0}
%\setcounter{page}{1}
%\section{Supplementary Materials : Additional Simulation Results and Proofs}\label{sec:sm}
%\renewcommand{\theequation}{\thesection.\arabic{equation}}
%\renewcommand{\theequation}{S.\arabic{equation}}
%\setcounter{equation}{0}
\par
\subsubsection{Notation}
Suppose, $\mathbf{\Phi_V}$ and $\mathbf{\phi_V}$ respectively denote the normal distribution and its density with mean $\mathbf{0}$ and covariance matrix $\bm{V}$. We will write $\mathbf{\Phi_V} = \mathbf{\Phi}$ and $\phi_{\mathbf{V}}=\phi$ when the dispersion matrix $\mathbf{V}$ is the identity matrix.
$C, C_1, C_2,\cdots$ denote generic constants that do not depend on the variables like $n, x$, and so on. $\bm{\nu}_1$, $\bm{\nu}_2$ denote vectors in $\mathscr{R}^p$, sometimes with some specific structures (as mentioned in the proofs). $(\mathbf{e_1},\ldots,\mathbf{e_p})'$ denote the standard basis of $\mathcal{R}^p$. For a non-negative integral vector $\bm{\alpha} = (\alpha_1, \alpha_2,\ldots,\alpha_l)'$ and a function $f = (f_1,f_2,\ldots,f_l):\ \mathcal{R}^l\ \rightarrow \ \mathcal{R}^l$, $l\geq 1$, let $|\bm{\alpha}| = \alpha_1 +\ldots+ \alpha_l$, $\bm{\alpha}! = \alpha_1!\ldots \alpha_l!$, $f^{\bm{\alpha}} = (f_1^{\alpha_1})\ldots(f_l^{\alpha_l})$, $D^{\bm{\alpha}}f_1 = D_1^{\alpha_1}\cdots D_l^{\alpha_l}f_1$, where $D_jf_1$ denotes the partial derivative of $f_1$ with respect to the $j$th component of $\bm{\alpha}$, $1\leq j \leq l$. We will write $D^{\bm{\alpha}}=D$ if $\bm{\alpha}$ has all the component equal to 1. For $\mathbf{t} =(t_1,t_2,\cdots t_l)'\in \mathcal{R}^l$ and $\mathbf{\alpha}$ as above, define $\mathbf{t}^{\bm{\alpha}} = t_1^{\alpha_1}\cdots t_l^{\alpha_l}$. For any two vectors $\bm{\alpha}, \bm{\beta} \in \mathcal{R}^k$, $\bm{\alpha} \leq \bm{\beta}$ means that each of the component of $\bm{\alpha}$ is smaller than that of $\bm{\beta}$. For a set $\bm{A}$ and real constants $a_1,a_2$, $a_1\bm{A}+a_2=\{a_1y+a_2:y\in \bm{A}\}$, $\partial A$ is the boundary of $A$ and $A^{\epsilon}$ denotes the $\epsilon-$neighbourhood of $A$ for any $\epsilon>0$. $\mathcal{N}$ is the set of natural numbers. $C(\cdot),C_1(\cdot),\dots$ denote generic constants which depend on only their arguments. Given two probability measures $P_1$ and $P_2$ defined on the same space $(\Omega,\mathcal{F})$, $P_1*P_2$ defines the measure on $(\Omega,\mathcal{F})$ by convolution of $P_1$ \& $P_2$ and $\|P_1-P_2\|=|P_1-P_2|(\Omega)$, $|P_1-P_2|$ being the total variation of $(P_1-P_2)$. For a function $g:\mathcal{R}^k\rightarrow \mathcal{R}^m$ with $g=(g_1,\dots,g_m)^\prime$, $$Grad [g(\bm{x})]=\Big(\Big(\frac{\partial g_i(\bm{x})}{\partial x_j}\Big)\Big)_{m\times k}.$$

For any natural number $m$, the class of sets $\mathcal{A}_m$ is the collection of Borel subsets of $\mathcal{R}^m$ satisfying 
\begin{equation}\label{eqn:setdef}
\sup\limits_{B \in  \mathcal{A}_m}\; \mathbf{\Phi}((\delta B)^{\epsilon}) = O(\epsilon) \;\;\; \text{as}\; \;\epsilon \downarrow 0.
\end{equation} 

For Lemma \ref{lem:3} below, define $\xi_{1,n,s}(\bm{t})=\Big(1+\sum_{i=1}^{s-2}n^{-r/2}\tilde{P}_r\big(i\bm{t}:\{\bar{\chi}_{\nu,n}\}\big)\Big)\exp\Big\{-\bm{t}^\prime \bm{E}_n \bm{t}/2\Big\}$ where $\bm{E}_n=n^{-1}\sum_{i=1}^{n}Var(Y_i)$ and $\bar{\chi}_{\nu,n}$ is the average $\nu$th cumulant of $Y_1,\dots,Y_n$. Define $\bar{\rho}_l=n^{-1}\sum_{i=1}^{n}$ $\mathbf{E}\|Y_i\|^{l}$, the average $l$th absolute moment of $\{Y_1,\dots,Y_n\}$. The polynomials $\tilde{P}_r\big(\bm{z}:\{\bar{\chi}_{\nu,n}\}\big)$ are defined on the pages of $51-53$ of Bhattacharya and Rao (1986). Define the identity $$\xi_{1,n,s}(\bm{t})\Big(\sum_{j=0}^{\infty}(-\|\bm{t}\|^2b_n^2)^j/j!\Big)=\xi_{n,s}(\bm{t})+ o\big(n^{-(s-2)/2}\big),$$ uniformly in $\|\bm{t}\|<1$, where $c_n$ is defined in Lemma \ref{lem:3}. $\psi_{n,s}(\cdot)$ is the Fourier inverse of $\xi_{n,s}(\cdot)$.

%Note that, $\check{\bm{b}}_n=O_p(n^{-\delta_1})$, by Lemma \ref{lem:betahat} and \ref{lem:Sigma}, described below. Suppose, $r_1=\min\{a\in \mathcal{N}:\|\check{\bm{b}}_n\|^{a+1}=o_p(n^{-1/2})\}$, $\mathcal{N}$ being the set of natural numbers. Define the conditional Lebesgue density of two-term Edgeworth expansion of $\bm{R}_{n}^*$ as
%\begin{align*}
%\xi_n^*(\bm{x})=&\phi(\bm{x})\Bigg[1+\sum_{k=1}^{r_1}\dfrac{1}{k!}\big{\{}\sum_{|\bm{\alpha}|=k}\check{\bm{b}}_n^{\alpha}H_{\bm{\alpha}}(\bm{x})\big{\}}+\dfrac{1}{\sqrt{n}}\bigg[ \dfrac{1}{6}\sum_{|\bm{\alpha}|=3}\bm{t}^{\bm{\alpha}}\bar{\bm{\xi}}_n^{*(1)}({\bm{\alpha}})H_{\bm{\alpha}}(\bm{x})\\ 
%&-\dfrac{1}{2\hat{\sigma}_n^2}\Big{\{}\sum_{|\bm{\alpha}|=1}\bm{t}^{\bm{\alpha}}\bar{\bm{\xi}}_n^{*(3)}({\bm{\alpha}})H_{\bm{\alpha}}(\bm{x})
%+\sum_{|\bm{\alpha}|=1}\sum_{|\bm{\zeta}|=2}\bm{t}^{\bm{\alpha}+\bm{\zeta}}\bar{\bm{\xi}}_n^{*(3)}({\bm{\alpha}})\bar{\bm{\xi}}_n^{*(1)}({\bm{\zeta}})H_{\bm{\alpha}+\bm{\zeta}}(\bm{x})\Big{\}}\bigg]\Bigg],
%\end{align*}
%where $x\in \mathcal{R}^q$, $\bar{\bm{\xi}}_{n}^{*(j)}(\bm{\alpha})=n^{-1}\sum_{i=1}^{n}\Big(\check{\bm{\xi}}_i^{(0)}\hat{\epsilon}_i^j\Big)^{\bm{\alpha}}$, $j=0,1,\ldots$ and $H_{\bm{\alpha}}(\bm{x})=(-D)^{\bm{\alpha}}\phi(\bm{x})$, where $\phi(\bm{\cdot})$ is the standard normal density on $\mathcal{R}^q$. 

\subsection{Proofs of Lemma \ref{lem:3}, \ref{lem:10} and \ref{lem:11}}

\begin{customlemma}{3}\label{lem:3}
Suppose $Y_1,\dots,Y_n$ are mean zero independent random vectors in $\mathcal{R}^k$ with $\bm{E}_n=n^{-1}$ $ \sum_{i=1}^{n}Var(Y_i)$ converging to some positive definite matrix $V$. Let $s\geq 3$ be an integer and $\bar{\rho}_{s
+\delta} =O(1)$ for some $\delta>0$. Additionally assume $Z$ to be a $N(\bm{0},\bm{I}_k)$ random vector which is independent of $Y_1,\dots,Y_n$ and the sequence $\{c_n\}_{n\geq 1}$ to be such that $c_n=O(n^{-d})$ \& $n^{-(s-2)/\tilde{k}}\log n = o(c_n^2)$ where $\tilde{k}=\max\{k+1,s+1\}$ \& $d>0$ is a constant. Then for any Borel set $B$ of $\mathcal{R}^k$,
\begin{align}\label{eqn:1}
\Big|\mathbf{P}\big(\sqrt{n}\bar{Y}+c_nZ\in B\big)- \int_{B}\psi_{n,s}(x)dx\Big|=o\Big(n^{-(s-2)/2}\Big),
\end{align}
where $\psi_{n,s}(\cdot)$ is defined above.
\end{customlemma}

Proof of Lemma \ref{lem:3}. Define $V_i=Y_iI\Big(\|Y_i\|\leq \sqrt{n}\Big)$ and $W_i=V_i-EV_i$. Suppose $\bar{\tilde{\chi}}_{\nu,n}$ is the average cumulant of $W_1,\dots,W_n$ and $D_n=n^{-1}\sum_{i=1}^{n}Var(W_i)$. Let $\tilde{\xi}_{1,n,s}$, $\tilde{\xi}_{n,s}$ and $\tilde{\psi}_{n,s}$ are respectively obtained from $\xi_{1,n,s}$, $\xi_{n,s}$ and $\psi_{n,s}$ with $\bar{\chi}_{\nu,n}$ replaced by $\bar{\tilde{\chi}}_{\nu,n}$ and $E_n$ replaced by $D_n$. For any Borel set $B\in \mathcal{R}^k$, define $B_n=B-n^{-1/2}\sum_{i=1}^{n}EV_i$. Then we have
\begin{align}\label{eqn:2}
&\Big|\mathbf{P}\big(\sqrt{n}\bar{Y}_n+c_nZ\in B\big)- \int_{B}\psi_{n,s}(x)dx\Big|\nonumber\\
\leq & \Big|\mathbf{P}\big(\sqrt{n}\bar{Y}_n+c_nZ\in B\big)- \mathbf{P}\Big(\sqrt{n}\bar{V}_n+c_nZ\in B\Big)\Big|\nonumber\\
&+ \Big|\mathbf{P}\big(\sqrt{n}\bar{W}_n+c_nZ\in B_n\big)- \int_{B_n}\tilde{\psi}_{n,s}(x)dx\Big| +\Big|\int_{B_n}\tilde{\psi}_{n,s}(x)dx- \int_{B}\psi_{n,s}(x)dx\Big|\nonumber\\
=& I_1 +I_2 +I_3\;\;\;\; \text{(say)}.
\end{align}

First we are going to show that $I_1=o\Big(n^{-(s-2)/2}\Big)$. Now writing $G_j$ and $G_j^\prime$ to be distributions of $n^{-1/2}Y_j$ and $n^{-1/2}V_j$, $j\in \{1,\dots,n\}$, we have
\begin{align}\label{eqn:3}
I_1 &\leq \sum_{j=1}^{n}\|G_j-G_j^\prime\|\nonumber\\
&= 2 \sum_{j=1}^{n}P\Big(\|Y_j\|>n^{1/2}\Big)\nonumber\\
&=o\Big(n^{-(s-2)/2}\Big),
\end{align}
due to the fact that $n^{-1}\sum_{j=1}^{n}E\|Y_j\|^{s+\delta}=O(1)$. Next we are going to show $I_3=o\Big(n^{-(s-2)/2}\Big)$. \\
Define $m_1=\inf \{j: b_n^{2j}=o\big(n^{-(s-2)/2}\big)\}$. Again note that the eigen values of $D_n$ are bounded away from 0, due to (14.18) in corollary 14.2 of Bhattacharya and Rao (1986) and the fact that $E_n$ converges to some positive definite matrix. Therefore we have
\begin{align}\label{eqn:4}
I_3=\Big|\int_{B_n}\tilde{\psi}_{n,s}^{m_1}(x)dx- \int_{B}\psi_{n,s}^{m_1}(x)dx\Big|+o\Big(n^{-(s-2)/2}\Big)=I_{31}+o\Big(n^{-(s-2)/2}\Big)\;\;\;\; (say),
\end{align}
uniformly for any Borel set $B$ of $\mathcal{R}^k$, where 
$$\psi_{n,s}^{m_1}(x)=\bigg\{\Big[\sum_{r=0}^{s-2}n^{-r/2}\tilde{P}_r\big(-D:\big\{\bar{\chi}_{\nu,n}\big\}\big)\Big]\Big[\sum_{j=0}^{m_1-1}2^{-j}(j!)^{-1}c_n^{2j}(D^\prime D)^j\Big]\bigg\}\phi_{E_n}(x)\;\;\; \text{and}$$
$$\tilde{\psi}_{n,s}^{m_1}(x)=\bigg\{\Big[\sum_{r=0}^{s-2}n^{-r/2}\tilde{P}_r\big(-D:\big\{\bar{\tilde{\chi}}_{\nu,n}\big\}\big)\Big]\Big[\sum_{j=0}^{m_1-1}2^{-j}(j!)^{-1}c_n^{2j}(D^\prime D)^j\Big]\bigg\}\phi_{D_n}(x).$$
Now writing $l(\bm{u})=\|\bm{u}\|/2,\; \bm{u}\in \mathcal{R}^k,$ and $a_n=n^{-1/2}\sum_{i=1}^{n}EV_i$, from (\ref{eqn:3}) we have
\begin{align}\label{eqn:5}
I_{31} \leq & \sum_{r=0}^{s-2}\sum_{j=0}^{m_1-1}n^{-r/2} b_n^{2j}\bigg[\int_{B_n}\Big|\Big\{\tilde{P}_r\big(-D:\big\{\bar{\chi}_{\nu,n}\big\}\big)\dfrac{l(-D)}{j!}\Big\}\phi_{E_n}(x)-\Big\{\tilde{P}_r\big(-D:\big\{\bar{\tilde{\chi}}_{\nu,n}\big\}\big)\dfrac{l(-D)}{j!}\Big\}\phi_{D_n}(x)\Big|dx\nonumber\\
&+\int_{B}\Big|\Big\{\tilde{P}_r\big(-D:\big\{\bar{\chi}_{\nu,n}\big\}\big)\dfrac{l(-D)}{j!}\Big\}\phi_{E_n}(x)-\Big\{\tilde{P}_r\big(-D:\big\{\bar{\chi}_{\nu,n}\big\}\big)\dfrac{l(-D)}{j!}\Big\}\phi_{E_n}(x-a_n)\Big|dx\bigg]\nonumber \\
& + o\Big(n^{-(s-2)/2}\Big)\nonumber\\
= & I_{311}+I_{312}+o\Big(n^{-(s-2)/2}\Big)\;\;\;\; \text{(say)}.
\end{align}
Now assume $E_n=I_k$, the $k\times k$ identity matrix. Then following the proof of Lemma 14.6 of Bhattacharya and Rao (86), it can be shown that $I_{311}+I_{312}=o\Big(n^{-(s-2)/2}\Big)$. Main ingredients of the proof are (14.74), (14.78), (14.79) and bounds similar to (14.80) and (14.86) in Bhattacharya and Rao (86). The general case when $E_n$ converges to a positive definite matrix, will follow essentially through the same line. Hence from (\ref{eqn:4}) and (\ref{eqn:5}), we have $I_3=o\big(n^{-(s-2)/2}\big)$. The last step is to show $I_2=o\Big(n^{(s-2)/2}\Big)$. Now let us write $\Gamma_n=\sqrt{n}\bar{W}_n+c_nZ$. Then recall that $$I_2=\Big|\mathbf{P}\big(\Gamma_n\in B_n\big)- \int_{B_n}\tilde{\psi}_{n,s}(x)dx\Big|.$$ By Theorem 4 of chapter 5 of Feller(2014), we can say that $\Gamma_n$ has density with respect to the Lebesgue measure. Let us call that density by $q_n(\cdot)$. Then we have
\begin{align}\label{eqn:6}
I_2 \leq \int \big|q_n(x)-\tilde{\psi}_{n,s}(x)\big|dx\leq \int \big|q_n(x)-\tilde{\psi}_{n,(\tilde{k}-1)}(x)\big|dx + \int \big|\tilde{\psi}_{n,s}(x)-\tilde{\psi}_{n,(\tilde{k}-1)}(x)\big|dx,
\end{align}
where $\tilde{k}=\max\{k+1,s+1\}$. Note that $\int \|x\|^j\big|q_n(x)-\tilde{\psi}_{n,(\tilde{k}-1)}(x)\big|dx< \infty$ for any $j\in \mathcal{N}$, since $\tilde{\psi}_{n,(\tilde{k}-1)}(x)$ has negative exponential term and $\bar{W}_n$ is bounded. Therefore by Lemma 11.6 of Bhattacharya and Rao (86) we have
\begin{align}\label{eqn:7}
I_2 &\leq  C(k)\bigg[\max_{|\bm{\beta}|\in \{0,\dots,(k+1)\}}\int \Big|D^{\bm{\beta}}\Big(\hat{q}(t)-\tilde{\xi}_{n,(\tilde{k}-1)}(t)\Big)\Big|dt\bigg] +   \int \big|\tilde{\psi}_{n,s}(x)-\tilde{\psi}_{n,(\tilde{k}-1)}(x)\big|dx \nonumber \\
& = I_{21} + I_{22}\;\;\;\; \text{(say)}.
\end{align}
Here $\hat{q}_n(\cdot)$ is the Fourier transform of the density $q(\cdot)$. Clearly $I_{22}=o\Big(n^{-(s-2)/2}\Big)$ by looking into the definition of $\tilde{\psi}_{n,s}(\cdot)$. Now define $$\breve{\xi}_{n,(\tilde{k}-1)}(t)=\bigg[\sum_{r=0}^{\tilde{k}-3}n^{-r/2}\tilde{P}_r\Big(i \bm{t}:\Big\{\bar{\tilde{\chi}}_{\nu,n}\Big\}\Big)\bigg]\exp{\Big(\dfrac{-\bm{t}^\prime\bm{D}_n\bm{t}-c_n^2\|\bm{t}\|^2}{2}\Big)}.$$ Then we have
\begin{align}\label{eqn:8}
I_{21}   &\leq C(k)\max_{|\bm{\beta}|\in \{0,\dots,(k+1)\}}\bigg[\int \Big|D^{\bm{\beta}}\Big(\hat{q}_n(t)-\breve{\xi}_{n,(\tilde{k}-1)}(t)\Big)\Big|dt +  \int \Big|D^{\bm{\beta}}\Big(\breve{\xi}_{n,(\tilde{k}-1)}(t)-\tilde{\xi}_{n,(\tilde{k}-1)}(t)\Big)\Big|dt\bigg] \nonumber\\
&= I_{211}+I_{212}\;\;\;\; \text{(say)}
\end{align}
First, we are going to show that $I_{212}=o\Big(n^{-(s-2)/2}\Big)$. Note that
\begin{align*}
\breve{\xi}_{n,(\tilde{k}-1)}(t)-\tilde{\xi}_{n,(\tilde{k}-1)}(t) = \bigg[\sum_{r=0}^{\tilde{k}-3}n^{-r/2}\tilde{P}_r\Big(i \bm{t}:\Big\{\bar{\tilde{\chi}}_{\nu,n}\Big\}\Big)\bigg]\exp{\Big(\dfrac{-\bm{t}^\prime\bm{D}_n\bm{t}}{2}}\Big)\sum_{j=m_2}^{\infty}\dfrac{c_n^{2j}\|\bm{t}\|^{2j}
(-1)^j}{2^j j!},
\end{align*}
where $m_2=m_2(r)=(s-2)^{-1}m_1(\tilde{k}-3-r)$. Therefore for any $\bm{\beta}\in \mathcal{N}^k$ with $|\bm{\beta}|\in \{0,\dots,k+1\}$ we have
\begin{align}\label{eqn:9}
&D^{\bm{\beta}}\Big(\breve{\xi}_{n,(\tilde{k}-1)}(t)-\tilde{\xi}_{n,(\tilde{k}-1)}(t)\Big)\nonumber\\
& = \sum^{*}\sum_{r=0}^{\tilde{k}-3}\sum_{j=m_2}^{\infty}C_1(\bm{\alpha},\bm{\beta},\bm{\gamma})\dfrac{n^{-r/2}(-1)^jc_n^{2j}}{2^j j!}\bigg[D^{\bm{\alpha}}\bigg(\tilde{P}_r\Big(i \bm{t}:\Big\{\bar{\tilde{\chi}}_{\nu,n}\Big\}\Big)\bigg)\bigg]\bigg[D^{\bm{\gamma}}\bigg(\exp{\Big(\dfrac{-\bm{t}^\prime\bm{D}_n\bm{t}}{2}}\Big)\bigg)\bigg]D^{\bm{\beta}-\bm{\alpha}-\bm{\gamma}}\Big(\|\bm{t}\|^{2j}\Big),
\end{align}
where $\sum^{*}$ is over $\bm{\alpha}, \bm{\gamma}\in \mathcal{N}^k$ such that $0\leq \bm{\alpha}, \bm{\gamma}\leq \bm{\beta}$. Since the degree of the polynomial $\tilde{P}_r\Big(i\bm{t}:\{\bar{\tilde{\chi}}_{\nu,n}\}\Big)$ is $3r$, $D^{\bm{\alpha}}\bigg(\tilde{P}_r\Big(i \bm{t}:\Big\{\bar{\tilde{\chi}}_{\nu,n}\Big\}\Big)\bigg)=0$ if $|\bm{\alpha}|>3r$.  When $|\bm{\alpha}|\leq 3r$, then recalling that $n^{-1}\sum_{i=1}^{n}E\|Y_i\|^{s}=O(1)$ and by Lemma 9.5 \& Lemma 14.1(v) of Bhattacharya and Rao (1986) we have
\begin{align}\label{eqn:10}
\bigg|D^{\bm{\alpha}}\bigg(\tilde{P}_r\Big(i \bm{t}:\Big\{\bar{\tilde{\chi}}_{\nu,n}\Big\}\Big)\bigg)\bigg| \leq 
\begin{dcases}
    C_2(\bm{\alpha},r)\big(\bar{\rho}_{s}\big)^{r/(s-2)}\Big(1+\big(\bar{\rho}_2\big)^{r(s-3)/(s-2)}\Big)(1+\|\bm{t}\|^{3r-|\bm{\alpha}|}),& \text{if } 0\leq r\leq (s-2)\\
    C_3(\bm{\alpha},r)n^{(r+2-s)/2}\bar{\rho}_{s}\Big(1+\big(\bar{\rho}_2\big)^{r-1}\Big)\Big(1+\|\bm{t}\|^{3r-|\bm{\alpha}|}\Big),              & \text{if } r> (s-2).
\end{dcases}
\end{align}

Again note that
\begin{align}\label{eqn:11}
\bigg|D^{\bm{\gamma}}\bigg(\exp{\Big(\dfrac{-\bm{t}^\prime\bm{D}_n\bm{t}}{2}}\Big)\bigg)\bigg| \leq C_4(\bm{\gamma})\Big(1+\|\bm{t}\|\Big)^
{|\bm{\gamma}|}\|\bm{D}_n\|^{|\bm{\gamma}|}\bigg(\exp{\Big(\dfrac{-\bm{t}^\prime\bm{D}_n\bm{t}}{2}}\Big)\bigg)
\end{align}
\begin{align}\label{eqn:12}
\text{and}\;\;\;\sum_{j=m_2}^{\infty}\bigg|\dfrac{c_n^{2j}D^{\bm{\beta}-\bm{\alpha}-\bm{\gamma}}\Big(\|\bm{t}\|^{2j}\Big)}{2^j j!}\bigg|\leq  C_5(\bm{\alpha}, \bm{\beta}, \bm{\gamma})c_n^{2m_3}\Big[e^{c_n^2/2}+\|\bm{t}\|^{m_3}\exp(c_n^2\|\bm{t}\|^2/2)\Big],
\end{align}
where $m_3=m_3(\bm{\alpha}, \bm{\beta}, \bm{\gamma},r)=\max\{m_2,|\bm{\beta}-\bm{\alpha}-\bm{\gamma}|/2\}$. Now combining (\ref{eqn:10})-(\ref{eqn:12}), from (\ref{eqn:9}) we have $I_{212}=o\Big(n^{-(s-2)/2}\Big)$. Last step is to show $I_{211}=o\Big(n^{-(s-2)/2}\Big)$. Recall that 
\begin{align}\label{eqn:13}
I_{211}&= C(k)\max_{|\bm{\beta}|\in \{0,\dots,(k+1)\}}\bigg[\int \Big|D^{\bm{\beta}}\Big(\hat{q}_n(t)-\breve{\xi}_{n,(\tilde{k}-1)}(t)\Big)\Big|dt\bigg]\nonumber \\
& \leq C(k)\max_{|\bm{\beta}|\in \{0,\dots,(k+1)\}}\bigg[\int_{A_n} \Big|D^{\bm{\beta}}\Big(\hat{q}_n(t)-\breve{\xi}_{n,(\tilde{k}-1)}(t)\Big)\Big|dt + \int_{A_n^c} \Big|D^{\bm{\beta}}\Big(\hat{q}(t)-\breve{\xi}_{n,(\tilde{k}-1)}(t)\Big)\Big|dt\bigg]\nonumber \\
& = I_{2111} + I_{2112}\;\;\;\; \text{(say)},
\end{align}
where $$A_n = \Bigg\{\bm{t}\in \mathcal{R}^k:\|\bm{t}\|\leq C_6(k)\lambda_n^{-1/2}\bigg(\dfrac{n^{1/2}}{\eta_{\tilde{k}}^{1/(\tilde{k}-2)}}\bigg)^{(\tilde{k}-2)/\tilde{k}}\Bigg\},$$
with $C_6(k)$ being some fixed positive constant, $\lambda_n$ being the largest eigen value of $\bm{D}_n$, $\eta_{\tilde{k}}=n^{-1}\sum_{i=1}^{n}E\|\bm{B}_nW_i\|^{\tilde{k}}$ and $\bm{B}_n^2=\bm{D}_n^{-1}$. Note that
\begin{align}\label{eqn:14}
&D^{\bm{\beta}}\Big(\hat{q}_n(t)-\breve{\xi}_{n,(\tilde{k}-1)}(t)\Big)\nonumber\\& =\sum_{\bm{0}\leq \bm{\alpha}\leq \bm{\beta}}C_7(\bm{\alpha},\bm{\beta})D^{\bm{\alpha}}\bigg(E\Big(e^{i\sqrt{n}\bm{t}^\prime}\bar{W}_n\Big)-\exp{\Big(\dfrac{-\bm{t}^\prime\bm{D}_n\bm{t}}{2}\Big)}\sum_{r=0}^{\tilde{k}-3}n^{-r/2}\tilde{P}_r\Big(i \bm{t}:\Big\{\bar{\tilde{\chi}}_{\nu,n}\Big\}\Big)\bigg)D^{\bm{\beta}-\bm{\alpha}}\bigg(\exp{\Big(\dfrac{-c_n^2\|\bm{t}\|^2}{2}\Big)}\bigg),
\end{align}
where $$\bigg|D^{\bm{\beta}-\bm{\alpha}}\bigg(\exp{\Big(\dfrac{-c_n^2\|\bm{t}\|^2}{2}\Big)}\bigg)\bigg|\leq C_8(\bm{\alpha},\bm{\beta})c_n^{2|\bm{\beta}-\bm{\alpha}|}\|\bm{t}\|^{|\bm{\beta}-\bm{\alpha}|}\exp{\Big(\dfrac{-c_n^2\|\bm{t}\|^2}{2}\Big)}\;\;\; \text{and}$$ 
by Theorem 9.11 and the following remark of Bhattacharya and Rao (86) we have
\begin{align}\label{eqn:15}
& \bigg|D^{\bm{\alpha}}\bigg(E\Big(e^{i\sqrt{n}\bm{t}^\prime\bar{W}_n}\Big)-\exp{\Big(\dfrac{-\bm{t}^\prime\bm{D}_n\bm{t}}{2}\Big)}\sum_{r=0}^{\tilde{k}-3}n^{-r/2}\tilde{P}_r\Big(i \bm{t}:\Big\{\bar{\tilde{\chi}}_{\nu,n}\Big\}\Big)\bigg)\bigg|\nonumber\\
&\leq C_9(k)\lambda_n^{|\bm{\alpha}|/2}\eta_{\tilde{k}}n^{-(\tilde{k}-2)/2}\Big[(\bm{t}^\prime D_n\bm{t})^{(\tilde{k}-|\bm{\alpha}|/2)}+(\bm{t}^\prime D_n\bm{t})^{(3(\tilde{k}-2)+|\bm{\alpha}|)/2}\Big]\exp{\Big(\dfrac{-\bm{t}^\prime\bm{D}_n\bm{t}}{4}\Big)}.
\end{align}
Now note that $\bar{\rho}_{s+\delta}=O(1)$ and $E_n$ converges to a positive definite matrix $\bm{E}$. Hence applying Lemma 14.1(v) (with $s^\prime=\tilde{k}$) and corollary 14.2 of Bhattacharya and Rao (86), from (\ref{eqn:14}) we have
$I_{2111}=o\Big(n^{-(s-2)/2}\Big)$. Again applying Lemma 14.1(v) and corollary 14.2 of Bhattacharya and Rao (86) we have $\eta_{\tilde{k}}\leq C_{10}(\tilde{k},s)n^{(\tilde{k}-s)/2}\bar{\rho}_s$ for large enough $n$ and $\lambda_n$ being converged to some positive number. Therefore we have for large enough $n$, $$A_n^c \subseteq B_n\;\; \text{where}\;\; B_n=\Big\{\bm{t}\in \mathcal{R}^k:\|\bm{t}\| > C_{11}(k,\bm{E})n^{(s-2)/{2\tilde{k}}}\Big\},$$ implying
\begin{align}\label{eqn:16}
I_{2112} &\leq C(k)\max_{|\bm{\beta}|\in \{0,\dots,(k+1)\}} \int_{B_n} \Big|D^{\bm{\beta}}\Big(\hat{q}_n(t)-\breve{\xi}_{n,(\tilde{k}-1)}(t)\Big)\Big|dt\nonumber\\
&\leq C(k)\max_{|\bm{\beta}|\in \{0,\dots,(k+1)\}} \bigg[\int_{B_n} \Big|D^{\bm{\beta}}\Big(\hat{q}_n(t)\Big)\Big|d\bm{t}+\int_{B_n}\Big|D^{\bm{\beta}}\Big(\breve{\xi}_{n,(\tilde{k}-1)}(t)\Big)\Big|dt\bigg]\nonumber\\
&= I_{21121} + I_{21122}\;\;\; \text{(say)},
\end{align}
for large enough $n$. To establish $I_{2112}=o\Big(n^{-(s-2)/2}\Big)$, first we are going to show $I_{21122}=o\Big(n^{-(s-2)/2}\Big)$. Note that 
\begin{align*}
D^{\bm{\beta}}\Big(\breve{\xi}_{n,(\tilde{k}-1)}(t)\Big)=\sum_{\bm{0}\leq \bm{\alpha}\leq \bm{\beta}}C_{12}(\bm{\alpha},\bm{\beta})D^{\bm{\alpha}}\bigg(\sum_{r=0}^{\tilde{k}-3}n^{-r/2}\tilde{P}_r\Big(i \bm{t}:\Big\{\bar{\tilde{\chi}}_{\nu,n}\Big\}\Big)\bigg)D^{\bm{\beta}-\bm{\alpha}}\bigg(\exp{\Big(\dfrac{-\bm{t}^\prime\tilde{\bm{D}}_n\bm{t}}{2}\Big)}\bigg),
\end{align*}
where $\tilde{\bm{D}}_n=\bm{D}_n+c_n^2\bm{I}_k$. We are going to use bounds (\ref{eqn:10}) and (\ref{eqn:11}) with $\bm{D}_n$ being replaced by $\tilde{\bm{D}}_n$. Note that by Corollary 14.2 of Bhattacharya and Rao (86) and the fact that $c_n=O(n^{-d})$, $\tilde{\bm{D}}_n$ converges to the positive definite matrix $\bm{E}$, which is the limit of $\bm{E}_n$. Hence those bounds will imply that for large enough $n$,
\begin{align}
I_{21122}&=C(k)\max_{|\bm{\beta}|\in \{0,\dots,(k+1)\}} \int_{B_n}\Big|D^{\bm{\beta}}\Big(\breve{\xi}_{n,(\tilde{k}-1)}(t)\Big)\Big|dt\nonumber\\
&\leq C_{13}(k,\bm{E}) n^{(\tilde{k}+1-s)/2}\int_{B_n} \Big(1+\|\bm{t}\|^{3(\tilde{k}-1)}\Big)\exp\Big(-C_{14}(\bm{E})\|\bm{t}\|^2/2\Big)d\bm{t}\nonumber\\
& \leq C_{15}(k,\bm{E}) n^{(\tilde{k}+1-s)/2}\int_{B_n} \exp\Big(-C_{14}(\bm{E})\|\bm{t}\|^2/4\Big)d\bm{t}.
\end{align}
Now apply Lemma 2 of the main paper to conclude that $I_{21122}=o\Big(n^{-(s-2)/2}\Big)$. Only remaining thing to show is $I_{21121}=o\Big(n^{-(s-2)/2}\Big)$. Note that  
\begin{align}\label{eqn:18}
D^{\bm{\beta}}\Big(\hat{q}_n(\bm{t})\Big)=\sum_{\bm{0}\leq \bm{\alpha}\leq \bm{\beta}}C_{16}(\bm{\alpha},\bm{\beta})D^{\bm{\alpha}}\bigg(E\Big(e^{i\sqrt{n}\bm{t}^\prime\bar{W}_n}\Big)\bigg)D^{\bm{\beta}-\bm{\alpha}}\bigg(\exp{\Big(\dfrac{-c_n^2\|\bm{t}\|^2}{2}\Big)}\bigg),
\end{align}
where 
\begin{align*}
&\Big|D^{\bm{\alpha}}\Big(E\Big(e^{i\sqrt{n}\bm{t}^\prime\bar{W}_n}\Big)\Big)\Big|\leq \Big|D^{\bm{\alpha}}\bigg(\prod_{i=1}^{n}E\Big(e^{i\bm{t}^\prime W_i/{\sqrt{n}}}\Big)\Big)\bigg|\\  \text{and}\;\;\;& \bigg|D^{\bm{\beta}-\bm{\alpha}}\bigg(\exp{\Big(\dfrac{-c_n^2\|\bm{t}\|^2}{2}\Big)}\bigg)\bigg| \leq C_{17}(\bm{\alpha},\bm{\beta})\Big(1+\|\bm{t}\|^{|\bm{\beta}-\bm{\alpha}|}\Big)\exp\Big(\dfrac{-c_n^2\|\bm{t}\|^2}{2}\Big).
\end{align*}
Now by Leibniz's rule of differentiation, $D^{\bm{\alpha}}\Big(E\Big(e^{i\sqrt{n}\bm{t}^\prime\bar{W}_n}\Big)\Big)$ is the sum of $n^{|\bm{\alpha}|}$ terms. A typical term is of the form
$$\prod_{i\not \in C_r}E\Big(e^{i\bm{t}^\prime W_{i}/{\sqrt{n}}}\Big)\prod_{l=1}^{r}D^{\bm{\beta}_l}\Big(E\Big(e^{i\bm{t}^\prime W_{i_l}/{\sqrt{n}}}\Big)\Big),$$
where $C_r=\{i_1,\dots,i_r\}\subset \{1,\dots,n\}$, $1\leq r\leq |\bm{\alpha}|$. $\bm{\beta}_1,\dots,\bm{\beta}_r$ are non-negative integral vectors satisfying $|\bm{\beta}_j|\geq 1$ for all $j\in \{1,\dots,r\}$ and $\sum_{j=1}^{r}\bm{\beta}_i=\bm{\alpha}$. Note that $\Big|D^{\bm{\beta}_l}\Big(E\Big(e^{i\bm{t}^\prime W_{i_l}/{\sqrt{n}}}\Big)\Big)\Big|\leq n^{-|\bm{\beta}_l|/2}E\|W_{i_l}\|^{|\bm{\beta}_l|}$ and $W_{j_l}\leq 2\sqrt{n}$, which imply that 
\begin{align*}
&\bigg|\prod_{i\not \in C_r}E\Big(e^{i\bm{t}^\prime W_{i}/{\sqrt{n}}}\Big)\prod_{l=1}^{r}D^{\bm{\beta}_l}\Big(E\Big(e^{i\bm{t}^\prime W_{i_l}/{\sqrt{n}}}\Big)\Big)\bigg| \leq 2^{\sum_{l=1}^{r}|\bm{\beta}_l|}=2^{|\bm{\alpha}|}\nonumber \\
\Rightarrow\;\; & \Big|D^{\bm{\alpha}}\Big(E\Big(e^{i\sqrt{n}\bm{t}^\prime\bar{W}_n}\Big)\Big)\Big|\leq (2n)^{|\bm{\alpha}|}.
\end{align*}
Let $K_n=C_{11}(k,\bm{E})n^{(s-2)/{2\tilde{k}}}$. Therefore from (\ref{eqn:18}), for large enough $n$ we have
\begin{align}
I_{21121}&\leq \Big[\max_{|\bm{\beta}|\in \{0,\dots,(k+1)\}}\sum_{\bm{0}\leq \bm{\alpha \leq \bm{\beta}}}C_{16}(\bm{\alpha},\bm{\beta})\Big](2n)^{k+1}\Big[\int_{B_n}\Big(1+\|\bm{t}\|^{k+1}\Big)\exp\Big(\dfrac{-c_n^2\|\bm{t}\|^2}{2}\Big)\Big]\nonumber\\
& \leq C_{18}(k)(2n)^{k+1}\int_{r\geq K_n}r^{k-1}\Big(1+r^{k+1}\Big)e^{-c_n^2r^2/2}dr\nonumber\\
& \leq C_{19}(k)(2n)^{k+1}c_n^{-1}\int_{r\geq K_n}\dfrac{1}{2\sqrt{\pi} c_n^{-1}}e^{-c_n^2r^2/4}dr\nonumber\\
&\leq C_{20}(k)n^{k+d+1} \int_{c_nK_n/{\sqrt{2}}}^{\infty}\dfrac{1}{\sqrt{2\pi} }e^{-z^2/2}dr\nonumber\\
& = o\Big(n^{-(s-2)/2}\Big).
\end{align}
The second inequality follows by considering polar transformation. Third inequality follows due to the assumptions that $n^{-(s-2)/{\tilde{k}}}(\log n)=o(c_n^2)$ and $c_n=O(n^{-d})$.
%$r^{2k}e^{-c_n^2r^2/4}\leq 1$ is required, i.e. $c_n^2r^2-8k \log r \geq 0$ is required. This is true if $c_n^2r^2 \geq 4k$, i.e. $r\geq \sqrt{4k}c_n^{-1}$. But actually we have $r\geq K_n$. Hence enough to have $K_n\geq \sqrt{4k}c_n^{-1}$ which is true due to $n^{-(s-2)/{\tilde{k}}}(\log n)=o(c_n^2)$.
The last equality is the implication of Lemma 2 presented in main paper. Therefore the proof of Lemma \ref{lem:3} is now complete.

\begin{customlemma}{10}\label{lem:10}
Assume the setup of Theorem 3 and let $X_i=y_ix_i$, $i\in \{1,\dots,n\}$. Define $\sigma_n^2=n^{-1}\sum_{i=1}^{n}Var(X_i)$ and $\bar{\chi}_{\nu,n}$ as the $\nu$th average cumulant of $\{(X_1-E(X_1)),\dots, (X_n-E(X_n))\}$.  $P_r\big(-\Phi_{\sigma_n^2}:\{\bar{\chi}_{\nu,n}\}\big)$ is the finite signed measure on $\mathcal{R}$ whose density is $\tilde{P}_r\big(-D: \{\bar{\chi}_{\nu,n}\}\big)\phi_{\sigma_n^2}(x)$. Let $S_0(x)=1$ and $S_1(x)=x-1/2$. Suppose $\sigma_n^2$ is bounded away from both $0$ \& $\infty$ and assumptions (C.1)-(C.3) of Theorem 3 hold. Then we have
\begin{align}\label{eqn:23}
\sup_{x\in \mathcal{R}}\Big|&\mathbf{P}\Big(n^{-1/2}\sum_{i=1}^{n}\big(X_i-E(X_i)\big)\leq x\Big)- \sum_{r=0}^{1}n^{-r/2}(-1)^rS_r(n\mu_n+n^{1/2}x)\dfrac{d^r}{dx^r}\Phi_{\sigma_n^2}(x)\nonumber\\
&-n^{-1/2}P_1\big(-\Phi_{\sigma_n^2}:\{\bar{\chi}_{\nu,n}\}\big)(x)\Big| = o\big(n^{-1/2}\big),
\end{align}
where $P_r\big(-\Phi_{\sigma_n^2}:\{\bar{\chi}_{\nu,n}\}\big)(x)$ is the $P_r\big(-\Phi_{\sigma_n^2}:\{\bar{\chi}_{\nu,n}\}\big)-$measure of the set $(-\infty,x]$.

\end{customlemma}

Proof of Lemma \ref{lem:10}. For any integer $\alpha$, define $p_n(x)=\mathbf{P}\big(\sum_{i=1}^{n}X_i = \alpha\big)$ and $x_{\alpha,n}= n^{-1/2}(\alpha-n \mu_n)$. Also define $\tilde{X}_n=n^{-1/2}\sum_{i=1}^{n}\big(X_i-E(X_i)\big)$ and $q_{n,3}(x)=n^{-1/2}\sum_{r=0}^{1}n^{-r/2}\tilde{P}_r\big(-D: \{\bar{\chi}_{\nu,n}\}\big)\phi_{\sigma_n^2}(x)$. Note that
\begin{align}\label{eqn:24}
\sup_{x\in \mathcal{R}}\Big|&\mathbf{P}\Big(n^{-1/2}\sum_{i=1}^{n}\big(X_i-E(X_i)\big)\leq x\Big)- \sum_{r=0}^{1}n^{-r/2}(-1)^rS_r(n\mu_n+n^{1/2}x)\dfrac{d^r}{dx^r}\Phi_{\sigma_n^2}(x)\nonumber\\
&-n^{-1/2}P_1\big(-\Phi_{\sigma_n^2}:\{\bar{\chi}_{\nu,n}\}\big)(x)\Big|\nonumber\\
\leq   \sup_{x\in \mathcal{R}}\big|&\mathbf{P}\Big(\tilde{X}_n\leq x\Big)-Q_{n,3}(x)\big| +\sup_{x\in \mathcal{R}}\big|Q_{n,3}(x)-\sum_{r=0}^{1}n^{-r/2}(-1)^rS_r(n\mu_n+n^{1/2}x)\dfrac{d^r}{dx^r}\Phi_{\sigma_n^2}(x)\nonumber\\
&-n^{-1/2}P_1\big(-\Phi_{\sigma_n^2}:\{\bar{\chi}_{\nu,n}\}\big)(x)\big|\\\nonumber
= J_1 + & J_2\;\;\;\; \text{(say)},
\end{align}
where $Q_{n,3}(x)=\sum_{\{\alpha: x_{\alpha,n}\leq x\}}q_{n,3}(x_{\alpha,n})$. Now the fact that $J_2=o\big(n^{-1/2}\big)$ follows from Theorem A.4.3 of Bhattacharya and Rao (86) and dropping terms of order $n^{-1}$. Now we are going to show $J_1=O\big(n^{-1}\big)$. Note that $$J_1 \leq \sum_{\alpha \in \Theta}\big|p_n(x_{\alpha,n})-q_{n,3}(x_{\alpha,n})\big| = J_{3}\;\;\;\; \text{(say)},$$ where $\Theta$ has cardinality $\leq C_{33}n$, since $\mathbf{P}\big(\big|n^{-1}\sum_{i=1}^{n}X_i\big| \leq C_{33}\big)=1$ for some constant $C_{33}>0$, due to the assumption that $\max\{|x_j|^5:j\in\{1,\dots,n\}\}=O(1)$. Hence $n^{-1}J_{3}\leq C_{33} \sup_{\alpha \in \Theta}\big|p_n(x_{\alpha,n})-q_{n,3}(x_{\alpha,n})\big| = C_{33}\sup_{\alpha \in \Theta}J_4(\alpha)$ (say). Hence enough to show $\sup_{\alpha \in \Theta}J_4(\alpha)=O\big(n^{-2}\big)$. Now define $g_j(t)=\mathbf{E}\big(e^{itX_j}\big)$ and $f_n(t)=\mathbf{E}(it\tilde{X}_n)$. Then we have 
\begin{align*}
f_n\big(\sqrt{n}t\big) = \sum_{\alpha \in \Theta}p_n\big(x_{\alpha,n}\big)e^{i \sqrt{n}t x_{\alpha,n}}.
\end{align*}
Hence by Fourier inversion formula for lattice random variables (cf. page 230 of Bhattacharya and Rao (86)), we have
\begin{align}\label{eqn:25}
p_n\big(x_{\alpha,n}\big) & = (2\pi)^{-1}\int_{\mathcal{F}^*} e^{-i \sqrt{n}t x_{\alpha,n}}f_n\big(\sqrt{n}t\big)dt\nonumber\\
& = (2\pi)^{-1}n^{-1/2}\int_{\sqrt{n}\mathcal{F}^*} e^{-i t x_{\alpha,n}}f_n\big(t\big)dt,
\end{align}
where $\mathcal{F}^*=(-\pi,\pi)$, the fundamental domain corresponding to the lattice distribution of $\sum_{i=1}^{n}X_i$.

Again note that 
\begin{align}\label{eqn:26}
q_{n,3}(x_{\alpha,n})= (2\pi)^{-1}n^{-1/2}\int_{\mathcal{R}}e^{-itx_{\alpha,n}}\sum_{r=0}^{1}n^{-r/2}\tilde{P}_r\big(it:\{\bar{\chi}_{\nu,n}\}\big)e^{-\sigma_n^2t^2/2}dt.
\end{align}
Now defining the set $E=\Big\{t\in \mathcal{R} :|t|\leq C_{31}(s)\sqrt{n}\min\big\{C_{33}^{-2}\sigma_n, C_{33}^{-5/3}\sigma_n^{5/3}\big\}\Big\}$, from (\ref{eqn:25}) \& (\ref{eqn:26}) we have
\begin{align}\label{eqn:27}
sup_{\alpha \in \Theta}J_4(\alpha) \leq &(2\pi)^{-1}n^{-1/2}\bigg[\int_{E}\Big|f_n(t)-\sum_{r=0}^{1}n^{-r/2}\tilde{P}_r\big(it:\{\bar{\chi}_{\nu,n}\}\big)e^{-\sigma_n^2t^2/2}\Big|dt\nonumber\\
& +\int_{\sqrt{n}\mathcal{F}^*\cap E^c}|f_n(t)|dt+\int_{\mathcal{R}\cap (\sqrt{n}\mathcal{F}^*)^c}\Big|\sum_{r=0}^{1}n^{-r/2}\tilde{P}_r\big(it:\{\bar{\chi}_{\nu,n}\}\big)e^{-\sigma_n^2t^2/2}\Big|dt\bigg]\nonumber\\
= & (2\pi)^{-1}n^{-1/2}\big(J_{41}+J_{42}+J_{43}\big)\;\;\;\;\; \text{(say)}.
\end{align}
Note that $J_{41}=O\big(n^{-3/2}\big)$ by applying Lemma 9 of the main paper with $s=5$. $J_{43}=O\big(n^{-3/2}\big)$ due to the presence of the exponential term in the integrand and the form of the set $E$. Moreover noting the form of the set $\mathcal{F}^*$, we can say that there exists constants $C_{34}>0$, $0<C_{35},C_{36}<\pi$ such that
\begin{align}\label{eqn:28}
J_{42} \leq C_{34} \sup_{t \in \sqrt{n}\mathcal{F}^* \cap E^c}\prod_{i=1}^{n}\big|g_j(n^{-1/2}t)\big| \leq   C_{34} \sup_{C_{35}\leq |t|\leq C_{36}} \big|\mathbf{E}(e^{ity_{i_1}})\big|^m  \leq C_{34} \delta^m,
\end{align}
for some $0<\delta<1$. Recall that $x_{i_j}=1$ for all $j\in \{1,\dots,m\}$. The last inequality is due to the fact that there is no period of $\mathbf{E}(e^{ity_{i_1}})$ in the interval $[C_{35},C_{36}]\cup [-C_{36},-C_{35}]$. Now $J_{42}=O(n^{-3/2})$ follows from (\ref{eqn:28}) since $m\geq (\log n)^2$. Therefore the proof is complete.\\

\begin{customlemma}{11}\label{lem:11}
Let $\breve{W}_1,\dots,\breve{W}_n$ be iid mean $\bm{0}$ non-degenerate random vectors in $\mathcal{R}^{l}$ for some natural number $l$, with finite fourth absolute moment and $\limsup_{\|\bm{t}\|\rightarrow \infty}\big|\mathbf{E}e^{i\bm{t}^\prime\breve{W}_1}\big|<1$ (i.e. Cramer's condition holds). Suppose $\breve{W}_i=(\breve{W}_{i1}^{\prime},\dots,\breve{W}_{im}^{\prime})^\prime$ where $\breve{W}_{ij}$ is a random vector in $\mathcal{R}^{l_j}$ and $\sum_{j=1}^{m}l_j=l$, $m$ being a fixed natural number. Consider the sequence of random variables $\tilde{W}_1,\dots,\tilde{W}_n$ where $\tilde{W}_i=(c_{i1}\breve{W}_{i1}^\prime,\dots,c_{im}\breve{W}_{im}^\prime)^\prime$. $\{c_{ij}:i\in \{1,\dots,n\}, j \in \{1,\dots,m\}\}$ is a collection of real numbers such that for any $j\in \{1,\dots,m\}$, $\Big\{n^{-1}\sum_{i=1}^{n}|c_{ij}|^4\Big\}=O(1)$ and $\liminf_{n\rightarrow \infty}n^{-1}\sum_{i=1}^{n}c_{ij}^2 > 0$. Also assume that $\tilde{\bm{V}}_n = Var(\tilde{W}_i)$ converges to some positive definite matrix and $\bar{\chi}_{\nu,n}$ denotes the average $\nu$th cumulant of $\tilde{W}_1,\dots,\tilde{W}_n$. Then we have
\begin{align}\label{eqn:29}
\sup_{\bm{B}\in \mathcal{A}_l}\Big|\mathbf{P}\Big(n^{-1/2}\sum_{i=1}^{n}\tilde{W}_i \in \bm{B}\Big)-\int_{\bm{B}}\Big[1+n^{-1/2}\tilde{P}_r\big(-D: \{\bar{\chi}_{\nu,n}\}\big)\Big]\phi_{\tilde{\bm{V}}_n}(\bm{t})d\bm{t}\Big|=o\big(n^{-1/2}\big),
\end{align}
where the collection of sets $\mathcal{A}_l$ is as defined in (\ref{eqn:setdef}).
\end{customlemma}

Proof of Lemma \ref{lem:11}. First note that $\tilde{W}_1,\dots, \tilde{W}_n$ is a sequence of independent random variables. Hence (\ref{eqn:29}) follows by Theorem 20.6 of Bhattacharya and Rao (1986), provided there exists $\delta_4 \in (0,1)$, independent of $n$, such that for all $\upsilon \leq \delta_4$,
\begin{align}\label{eqn:30}
n^{-1}\sum_{i=1}^{n}\mathbf{E}\big\|\tilde{W}_i\big\|^{3}\mathbf{1}\Big(\big\|\tilde{W}_i\big\|> \upsilon \sqrt{n}\Big)=o(1)
\end{align}
and
\begin{align}\label{eqn:31}
\max_{|\bm{\alpha}|\leq l+2}\int_{\|\bm{t}\|\geq \upsilon\sqrt{n}}\Big|D^{\bm{\alpha}}\mathbf{E}\exp(i\bm{t}^\prime\bm{R}_{1n}^{\dagger})\Big|d\bm{t}=o\Big(n^{-1/2}\Big)
\end{align}
where $\bm{R}_{1n}^{\dagger}=n^{-1/2}\sum_{i=1}^{n}\big(\bm{Z}_{i}-\mathbf{E}\bm{Z}_i\big)$ with $$\bm{Z}_{i}=\tilde{W}_i\mathbf{1}\Big(\big\|\tilde{W}_i\big\|\leq \upsilon \sqrt{n}\Big).$$

First consider (\ref{eqn:30}). Note that 
$\max\Big\{|c_{ij}|: i\in\{1,\dots,n\}, j \in \{1,\dots,m\}\Big\} = O\big(n^{1/4}\big).$
Therefore, we have for any $\upsilon >0$,
\begin{align*}
&n^{-1}\sum_{i=1}^{n}\mathbf{E}\big\|\tilde{W}_i\big\|^{3}\mathbf{1}\Big(\big\|\tilde{W}_i\big\|> \upsilon \sqrt{n}\Big)\\
\leq& n^{-1}\sum_{i=1}^{n}\mathbf{E}\Big(\sum_{j=1}^{m}c_{ij}^2\big\|\breve{W}_{ij}\big\|^2\Big)^{3/2}\mathbf{1}\Big(\sum_{j=1}^{m}c_{ij}^2\big\|\breve{W}_{ij}\big\|^2> \upsilon^2 n\Big)\\
\leq& n^{-1}\sum_{i=1}^{n}\Big(1+\sum_{j=1}^{m}c_{ij}^2\Big)^2\mathbf{E}\bigg[\big\|\breve{W}_{1}\big\|^3\mathbf{1}\Big(\big\|\breve{W}_{1}\big\|^2>C_{37}\upsilon^2 n^{1/2}\Big)\bigg]\\
=& o(1).
\end{align*}
Now consider (\ref{eqn:31}). Note that for any $|\bm{\alpha}|\leq l+2$,  $|D^{\bm{\alpha}}\mathbf{E}\exp(i\bm{t}^\prime\bm{R}_{1n}^{\dagger})|$ is bounded above by a sum of $n^{|\alpha|}$-terms, each of which is bounded above by
\begin{align}\label{eqn:32}
C_{38}(\alpha) \cdot n^{-|\bm{\alpha}|/2} \max\{\mathbf{E}\|\bm{Z}_{i}-\mathbf{E}\bm{Z}_i\|^{|\bm{\alpha}|}: k \in \bm{I}_n\} \cdot \prod_{i \in \bm{I}^{\mathsf{c}}_n}|\mathbf{E}\exp(i\bm{t}^\prime\bm{Z}_{i}/\sqrt{n})|
\end{align}
where $\bm{I}_n\subset \{1,\dots, n\}$ is of size $|\bm{\alpha}|$ and $\bm{I^{\mathsf{c}}_n}=\{1,\dots,n\}\backslash \bm{I}_n$. Now for any $\omega>0$ and $\bm{t}\in \mathcal{R}^{l_j}$, define the set $$\bm{B}_n^{(j)}(\bm{t},\omega) = \Big\{i:1\leq i \leq n  \; \text{and}\; |c_{ij}|\|\bm{t}\| > \omega\Big\}.$$ Hence for any $\bm{t}\in \mathcal{R}^{l}$ writing $\bm{t}=\big(\bm{t}_1^\prime,\dots,\bm{t}_{m}^\prime\big)^\prime$, $\bm{t}_j$ is of length $l_j$, we have
\begin{align*}
&\sup\Bigg{\{}\prod_{i \in \bm{I}^{\mathsf{c}}_n}|\mathbf{E}\exp(i\bm{t}^\prime\bm{Z}_{k}/\sqrt{n})|:\|\bm{t}\|\geq \upsilon \sqrt{n} \Bigg{\}} \\
=&\sup\Bigg{\{}\prod_{i \in \bm{I}^{\mathsf{c}}_n}|\mathbf{E}\exp(i\bm{t}^\prime\bm{Z}_{k})|:\|\bm{t}\|^2\geq \upsilon^2  \Bigg{\}} \\
\leq & \max\Bigg\{\sup\bigg{\{}\prod_{i \in \bm{I}^{\mathsf{c}}_n\cap \bm{B}_n^{(j)}\Big(\dfrac{\bm{t}_j}{\|\bm{t}_j\|}, \upsilon/\sqrt{2}\Big)}\Big[|\mathbf{E}\exp\Big(ic_{ij}\bm{t}_j^\prime\breve{W}_{1j}\Big)|+\mathbf{P}\Big(\|\breve{W}_1\|> C_{37}\upsilon^2 n^{1/2}\Big)\Big]\\
&:\|\bm{t}_j\|\geq \upsilon/\sqrt{2}  \bigg{\}}:j\in \big\{1,\dots,m\big\}\Bigg\}
\end{align*}
Now since $\big|\bm{I}^{\mathsf{c}}_n\big|\geq \Big|\bm{I}^{\mathsf{c}}_n\cap \bm{B}_n^{(j)}\Big(\dfrac{\bm{t}_j}{\|\bm{t}_j\|}, \upsilon/\sqrt{2}\Big)\Big|\geq \big|\bm{B}_n^{(j)}\Big(\dfrac{\bm{t}_j}{\|\bm{t}_j\|}, \upsilon/\sqrt{2}\Big)\big|-|\alpha|$, due to Cramer's condition we have
\begin{align}\label{eqn:33}
&\sup\bigg{\{}\prod_{i \in \bm{I}^{\mathsf{c}}_n\cap \bm{B}_n^{(j)}\Big(\dfrac{\bm{t}_j}{\|\bm{t}_j\|}, \upsilon/\sqrt{2}\Big)}\Big[|\mathbf{E}\exp\Big(ic_{ij}\bm{t}_j^\prime\breve{W}_{1j}\Big)|+\mathbf{P}\Big(\|\breve{W}_1\|> C_{37}\upsilon^2 n^{1/2}\Big)\Big]
:\|\bm{t}_q\|\geq \upsilon/\sqrt{2}  \bigg{\}}\nonumber\\
&\leq \theta^{\Big|\bm{B}_n^{(j)}\Big(\dfrac{\bm{t}_j}{\|\bm{t}_j\|}, \upsilon/\sqrt{2}\Big)\Big|-\big|\alpha\big|}
\end{align}

Next note that $\liminf_{n\rightarrow \infty}n^{-1}\sum_{i=1}^{n}c_{ij}^2>0$ for all $j\in \{1,\dots,m\}$. Therefore for any $j\in \{1,\dots,m\}$, $\bm{u}\in \mathcal{R}^{l_j}$ with $|u|=1$, there exists $0<\delta_5<1$ such that for sufficiently large $n$ we have 
\begin{align*}
\dfrac{n\delta_5}{2}\leq & \sum_{i=1}^{n}\big|u c_{ij}\big|^2\\
\leq & \max\Big\{\big|c_{ij}\big|^2:1\leq i \leq n\Big\}\cdot |\bm{B}_n^{(j)}(u,\omega)|+\Big(n-|\bm{B}_n^{(j)}(u,\omega)|\Big)\cdot \omega^2\\
\leq & C_{38}\cdot n^{1/2}\cdot |\bm{B}_n^{(j)}(\bm{u},\omega)| + n\omega^2
\end{align*}
which implies $|\bm{B}_n^{(j)}(\bm{u},\omega)|\geq C_{39}\cdot n^{1/2}$ whenever $\omega < \sqrt{\delta_5/2}$. Therefore taking $\delta_4=\sqrt{\delta_5/3}$, (\ref{eqn:31}) follows from (\ref{eqn:32}) and (\ref{eqn:33}).\\

\subsection{Supplementary Simulation Details}\label{sec_supp_sim}
In this section we present expanded forms of the pivots and the forms of the confidence intervals obtained based on our proposed Bootstrap method. Code details for the reproduction of the results of Section 6 and 7 of the main manuscript can be supplied if required.

Recall that our model is 
\begin{align*}
    y_i =\; & 1, \; w.p. \; p(\bm{\beta}|\bm{x}_i),\\
    = \;& 0, \; w.p. \; [1 - p(\bm{\beta}|\bm{x}_i)],
\end{align*}
where  $p(\bm{\beta}|\bm{x}_i)  = \frac{\exp{(x_i^T \beta)}}{1+\exp{(x_i^T \beta)}}$, $i\in \{1,\dots,n\}$. Here $y_1,\ldots,y_n$ are independent binary responses and $\bm{x}_1,\ldots,\bm{x}_n$ are known non-random design vectors. $\bm{\beta}=(\beta_{1,n},\ldots, \beta_{p,n})$ is the $p$-dimensional vector of regression parameters. For the rest of this section $\bm{x}_{i,\bm{A}}$ denotes the sub vector of $\bm{x}_i$ comprising of only components belonging to the set $\bm{A}$ where $\bm{A}\subseteq \{1,\dots,p\}$. For any vector $\bm{\gamma}$ of length $p$, $\bm{\gamma}_{\bm{A}}$ is the sub vector of $\bm{\gamma}$ comprising of only components belonging to the set $\bm{A}$.

The logistic regression estimator $\bm{\beta}_n$ of $\bm{\beta}$ is defined as $$\hat{\bm{\beta}}_n = \mbox{Argmax}_{\bm{\beta}} L(\bm{\beta}|y_1,\dots,y_n, \bm{x}_1,\dots,\bm{x}_n),$$ where $L(\bm{\beta}|y_1,\dots,y_n, \bm{x}_1,\dots,\bm{x}_n)=\prod_{i=1}^{n}p(\bm{x}_i)^{y_i}(1-p(\bm{x}_i))^{1-y_i}$ is the likelihood. The Bootstrap version [hereafter referred to as PEBBLE] $\bm{\hat{\beta}}_n^*$ of $\hat{\bm{\beta}}_n$ is defined as
\begin{align*}
\hat{\bm{\beta}}_n^* = \operatorname*{arg\,max}_{\bm{t}}\Bigg[\sum_{i=1}^{n}\Big\{(y_i - \hat{p}(\bm{x}_i))\bm{x}_i^{\prime}\bm{t}\Big\}(G^*_i-\mu_{G^*}) + \mu_{G^*}\sum_{i=1}^{n}\Big\{\hat{p}(\bm{x}_i)(\bm{x}_i^{\prime}\bm{t})-\log (1+e^{\bm{x}_i^{\prime}\bm{t}})\Big\}\Bigg],
\end{align*}
where $G_1,\dots,G_n^*$ are iid copies of a non-negative \& non-degenerate random variable $G^*$ with $Var(G^*)=\mu_{G^*}^2$ and $\mathbf{E}(G^*-\mu_{G^*})^3=\mu_{G^*}^3$. One example of the distribution of $G^*$ is $Beta(1/2,3/2)$.

%\subsection{Forms of the confidence Sets}
\subsubsection{Form of the confidence region for the parameter}
The original studentized pivot for the parameter vector $\bm{\beta}$ is
$$\check{\mathbf{H}}_n=\hat{\bm{M}}_n^{-1/2}\hat{\bm{L}}_n\big[\sqrt{n}\big(\hat{\bm{\beta}}_n-\bm{\beta}\big)\big] +\hat{\bm{M}}_n^{-1/2}b_n\bm{Z},$$
where $\hat{\bm{L}}_n=n^{-1}\sum_{i=1}^{n}\bm{x}_i\bm{x}_i^{\prime}e^{\bm{x}_i^{\prime}\hat{\bm{\beta}}_n}(1+e^{\bm{x}_i^{\prime}\hat{\bm{\beta}}_n})^{-2}$, $\hat{\bm{M}}_n=n^{-1}\sum_{i=1}^{n}\big(y_i-\hat{p}(\bm{x}_i)\big)^2\bm{x}_i\bm{x}_i^\prime$, $\hat{p}(\bm{x}_i)  = \frac{\exp{(x_i^T \hat{\bm{\beta}}_n)}}{1+\exp{(x_i^T \hat{\bm{\beta}}_n)}}$. $\bm{Z}$ is distributed as $N\big(\bm{0},\bm{D}\big)$ where $\bm{D}$ is a $p\times p$ diagonal matrix, independent of $y_1,\dots,y_n$. $\{b_n\}_{n\geq 1}$ is a sequence of real numbers such that $b_n=O(n^{-d})$ and $n^{-1/p_1}\log n = o(b_n^2)$ where $d>0$ is a constant and $p_1=\max\{p+1,4\}$. Corresponding PEBBLE version of the studentized pivot is defined as
$$\check{\mathbf{H}}_n^*=\hat{\bm{M}}_n^{*-1/2}\bm{L}_n^*\big[\sqrt{n}\big(\hat{\bm{\beta}}_n^*-\hat{\bm{\beta}}_n\big)\big] +\hat{\bm{M}}_n^{*-1/2}b_n\bm{Z}^*,$$
where $\bm{L}_n^*=n^{-1}\sum_{i=1}^{n}\bm{x}_i\bm{x}_i^{\prime}e^{\bm{x}_i^{\prime}\hat{\bm{\beta}}_n^*}\big(1+e^{\bm{x}_i^{\prime}\hat{\bm{\beta}}_n^*}\big)^{-2}$ and $\hat{\bm{M}}_n^*=n^{-1}\sum_{i=1}^{n}\big(y_i-\hat{p}(\bm{x}_i)\big)^2\bm{x}_i\bm{x}_i^\prime\mu_{G^*}^{-2}(G_i^*-\mu_{G^*})^2$. $\bm{Z}^*$ has the same distribution as $\bm{Z}$, independent of $y_1,\dots,y_n$ and $G_1^*,\dots,G_n^*$.

For some $\alpha \in (0,1)$, let $\Big(\|\check{\mathbf{H}}_n^*\|\Big)_{\alpha}$ be the $\alpha$th quantile of the Bootstrap distribution of $\|\check{\mathbf{H}}_n^*\|$. Then the $100(1-\alpha) \%$ confidence region of $\bm{\beta}$ is given by
$$\bigg\{\bm{\beta}:\|\check{\mathbf{H}}_n\|\leq \Big(\|\check{\mathbf{H}}_n^*\|\Big)_{(1-\alpha)}\bigg\}.$$

\subsubsection{Form of the confidence intervals for the components of the parameter}
The pivotal quantity for the $j$th component of $\bm{\beta}$ is formulated as
\begin{align*}
    \check{H}_{j,n} = \hat{\Sigma}_{j,n}^{-\frac{1}{2}}\bigg(\sqrt{n}(\hat{\beta}_{j,n}-\beta_j)+b_n\Big(\hat{\bm{L}}_n^{-1}\Big)_{j\cdot}^\prime \bm{Z}\bigg),
\end{align*}
where $\hat{\beta}_{j,n}$ \& $\beta_{j,n}$ are respectively the $j$th component of $\hat{\bm{\beta}}_n$ and $\bm{\beta}$, $j\in \{1,\dots,p\}$. $\hat{\Sigma}_{j,n}$ is the $(j,j)$-th element of $\hat{\bm{\Sigma}}_n$ where $\hat{\bm{\Sigma}}_n = \hat{\bm{L}}_n^{-1} \bm{\hat{M}}_n \hat{\bm{L}}_n^{-1}$ and $\Big(\hat{L}_n^{-1}\Big)_{j\cdot}^\prime$ is the $j$-th row of $\hat{\bm{L}}_n^{-1}$. Similarly the Bootstrap version corresponding to $\check{H}_{j,n}$ is defined as
\begin{align*}
\check{H}^*_{j,n} = \Sigma_{j,n}^{*-\frac{1}{2}}\bigg(\sqrt{n}(\hat{\beta}^*_{j,n}-\hat{\beta}_{j,n})+b_n\Big(\bm{L}_n^{*-1}\Big)^\prime_{j\cdot} \bm{Z}^*\bigg),
\end{align*}
where $\hat{\beta}_{j,n}^*$ is the $j$th component of the vector $\hat{\bm{\beta}}_n^*$, $j\in \{1,\dots,p\}$. $\Sigma_{j,n}^*$ is the $(j,j)$-th element of $\bm{\Sigma}_n^*$ where $\bm{\Sigma}^*_n = \bm{L}_n^{*-1}\hat{\bm{M}}_n^*\bm{L}_n^{*-1}$ and $\Big(\bm{L}_n^{*-1}\Big)_{j\cdot}^\prime$ is the $j$-th row of $\hat{\bm{L}}_n^{*-1}$.  

Now define $\big(\check{H}^*_{j,n}\big)_{\alpha}$ to be the $\alpha$th quantile of the Bootstrap distribution of $\check{H}^*_{j,n}$ for some $\alpha \in (0,1)$. Then $100(1-\alpha)\%$ two-sided confidence interval of $\beta_{j}$ is given by 
$$\Bigg[\bigg\{\hat{\beta}_{j,n}-\frac{\hat{\Sigma} _{j,n}^{1/2}u_{1j}^*}{\sqrt{n}}\bigg\}, \bigg\{\hat{\beta}_{j,n}-\frac{\hat{\Sigma}_{j,n}^{1/2}l_{1j}^*}{\sqrt{n}}\bigg\}\Bigg],$$
where $l_{1j}^* = \Big[(\check{H}_{j,n}^*)_{\alpha/2} - b_n\hat{\Sigma}_{j,n}^{-\frac{1}{2}}\Big(\hat{L}_n^{-1}\Big)^\prime_{j\cdot}\bm{Z}\Big]$ and $u_{1j}^*=\Big[(\check{H}_{j,n}^*)_{(1-\alpha)/2} - b_n\hat{\Sigma}_{j,n}^{-\frac{1}{2}}\Big(\hat{L}_n^{-1}\Big)^\prime_{j\cdot}\bm{Z}\Big]$. Again $100(1-\alpha)\%$ lower and upper confidence intervals of $\beta_{j}$ are respectively given by
$$\Bigg(-\infty, \bigg\{\hat{\beta}_{j,n}-\frac{\hat{\Sigma}_{j,n}^{1/2}l_{2j}^*}{\sqrt{n}}\bigg\}\Bigg]\;\;\;\; \text{and}\;\;\;\; \Bigg[\bigg\{\hat{\beta}_{j,n}-\frac{\hat{\Sigma} _{j,n}^{1/2}u_{2j}^*}{\sqrt{n}}\bigg\}, \infty\Bigg),$$
 where $l_{2j}^* = \Big[(\check{H}_{j,n}^*)_{\alpha} - b_n\hat{\Sigma}_{j,n}^{-\frac{1}{2}}\Big(\hat{L}_n^{-1}\Big)^\prime_{j\cdot}\bm{Z}\Big]$ and $u_{2j}^*=\Big[(\check{H}_{j,n}^*)_{(1-\alpha)} - b_n\hat{\Sigma}_{j,n}^{-\frac{1}{2}}\Big(\hat{L}_n^{-1}\Big)^\prime_{j\cdot}\bm{Z}\Big]$.

\end{document}